\documentclass[10pt, a4paper]{amsart}
\usepackage[margin=2cm]{geometry}

\usepackage[T1]{fontenc}
\usepackage[utf8]{inputenc}
\usepackage{lmodern}
\usepackage[british]{babel}

\usepackage{csquotes}
\usepackage[
sorting=nyt,
    doi=true,
    url=true,
    maxbibnames=4,
    backend=biber,
    giveninits=true,
    ]{biblatex}
\bibliography{literature}

\AtEveryBibitem{\clearfield{month}}
\AtEveryBibitem{\clearfield{day}}

\newcommand{\citep}[1]{(\citeauthor{#1}, \citeyear{#1})}

\usepackage{graphicx}
\usepackage{subcaption}
\usepackage[dvipsnames]{xcolor}
  \definecolor{darkred}{RGB}{139,0,0}
  \definecolor{mediumblue}{RGB}{0,0,205}
  \definecolor{forestgreen}{RGB}{34,139,34}

\usepackage[colorlinks=true]{hyperref}
\hypersetup{linkcolor=darkred, urlcolor=forestgreen, citecolor=mediumblue}
\usepackage{cleveref}

\usepackage{amsmath, amsfonts, amssymb}
\usepackage{mathtools}
\usepackage{commath}

\usepackage{bm}
\usepackage{calrsfs}
\usepackage{stmaryrd}
\DeclareMathAlphabet{\pazocal}{OMS}{zplm}{m}{n}
\DeclareMathAlphabet{\pazocalbf}{OMS}{cmsy}{b}{n}
\usepackage{etoolbox}
\preto\subequations{\ifhmode\unskip\fi}
\usepackage{nicefrac}

\usepackage{amsthm}
\theoremstyle{remark}
\newtheorem{rmk}{Remark}[section]

\usepackage{enumitem}
  \setlist[enumerate]{nosep, topsep=0pt, wide = 1em, leftmargin=*}
\usepackage{booktabs}
\AtBeginEnvironment{tabular}{\small}
\usepackage{siunitx}
\usepackage[dvipsnames]{xcolor}

\newcommand{\eb}{\bm{e}}

\newcommand{\nb}{\bm{n}}

\newcommand{\ub}{\bm{u}}

\newcommand{\mb}{\bm{m}}

\newcommand{\ubu}{\mathbf{u}}
\newcommand{\Mbu}{\mathbf{M}}
\newcommand{\Fbu}{\mathbf{F}}
\newcommand{\Gbu}{\mathbf{G}}

\newcommand{\rhobar}{\overline{\rho}}
\newcommand{\qbar}{\overline{q}}
\newcommand{\pbar}{\overline{p}}
\newcommand{\Tbar}{\overline{T}}

\renewcommand{\O}{\Omega}

\DeclareMathOperator{\Id}{Id}
\DeclareMathOperator{\spann}{span}

\newcommand{\lprod}[3]{\big(#1,#2\big)_{#3}}

\newcommand{\jump}[1]{\llbracket#1\rrbracket}
\newcommand{\mean}[1]{\{\!\{#1\}\!\}}
\newcommand{\restr}[2]{{\left.\kern-\nulldelimiterspace#1\right|_{#2}}}

\newcommand{\pref}{p_\text{ref}}

\newcommand{\Tref}{T_\text{ref}}
\newcommand{\Lref}{L_\text{ref}}
\newcommand{\eref}{e_\text{ref}}

\newcommand{\HC}{\pazocal{H}}

\newcommand{\Th}{\mathcal{T}_h}
\newcommand{\Fh}{\mathcal{F}_h}
\newcommand{\Fhi}{\Fh^i}
\newcommand{\Fhb}{\Fh^b}

\newcommand{\PP}{\mathbb{P}}
\newcommand{\QQ}{\mathbb{Q}}

\newcommand{\dt}{\Delta t}

\begin{document}

\title[A discontinuous Galerkin approach for atmospheric flows]{A discontinuous Galerkin approach for atmospheric flows with implicit condensation}

\author[S.~Hittmeir]{Sabine Hittmeir}
\address{Faculty of Mathematics, University of Vienna, Austria}
\email{sabine.hittmeir@univie.ac.at}

\author[P.~L.~Lederer]{Philip L. Lederer}
\address{Department of Applied Mathematics, University of Twente, Netherlands}
\email{p.l.lederer@utwente.nl}

\author[J.~Sch\"oberl]{Joachim Sch\"oberl}
\address{Institute for Analysis and Scientific Computing, TU Wien, Austria}
\email{joachim.schoeberl@tuwien.ac.at}

\author[H.~v.~Wahl]{Henry von Wahl}
\address{Faculty of Mathematics, University of Vienna, Austria}
\email{henry.wahl@univie.ac.at}

\begin{abstract}
We present a discontinuous Galerkin method for moist
atmospheric dynamics, with and without warm rain. By considering a
combined density for water vapour and cloud water, we avoid the need
to model and compute a source term for condensation. We recover
the vapour and cloud densities by solving a pointwise non-linear
problem each time step. Consequently, we enforce the requirement
for the water vapour not to be supersaturated implicitly. Together
with an explicit time-stepping scheme, the method is highly
parallelisable and can utilise high-performance computing hardware.
Furthermore, the discretisation works on structured and unstructured
meshes in two and three spatial dimensions. We illustrate the
performance of our approach using several test cases in two and three
spatial dimensions. In the case of a smooth, exact solution, we
illustrate the optimal higher-order convergence rates of the method.
\end{abstract}

\maketitle

\keywords{
Atmospheric flow, Discontinuous Galerkin, Matrix-free, High-order, compressible Euler equations with source terms, moisture, implicit condensation, Hyperbolic conservation laws}

\section{Introduction}\label{sec.intro}

Precipitation still causes one of the largest uncertainties in weather forecasts and climate models. This is evident, as the modelling has to build a bridge across many scales to arrive from the processes in and around precipitation particles at the evolutionary dynamics of clouds and cloud systems on the mesoscales. Thus, cloud models exist with very different ranges of complexity. To model the moisture balances, we implement a bulk microphysical warm cloud model corresponding to the Kessler scheme, which is widely used in meteorology. Then, compared to the dry dynamics, the thermodynamic equation gets much more complex, as it provides a strong coupling to the moisture balances. This additional complexity is also due to the different heat capacities for vapour and liquid water and the dependence of latent heat on the temperature, which are often neglected. In this work, we retain all these thermodynamic details, as they have been demonstrated to be essential, e.g. in the case of deep convective clouds via asymptotic analysis in~\cite{HK17}.

While the equations of dry air are generally accepted, the equations governing cloudy air are still actively debated. Therefore, it is necessary to develop numerical methods to simulate these complex equations accurately and efficiently. Typically, numerical schemes are established for individual model reductions obtained by scale analysis. For an overview of numerical methods for numerical weather prediction, see, for example, \cite{MKM15} and the references therein.

On the other hand, increasing, massively-parallel computing power allows for ever more fine-scale computations for weather and climate models. However, numerical methods have for the most part been based on structured grids and discretisations exploiting the structure of the computational grids. While the orthogonality of the grid can lead to a number of desirable properties, see, for example, \cite{ST11}, there are drawbacks in the setting of massively parallel high-performance computing hardware. In particular, load balancing can be an issue when the discretisation is limited to column-based subdivisions of the mesh. Furthermore, for global discretisations, structured grids on the globe lead to strong clustering of grid points at the poles, creating an additional bottleneck. To avoid these issues, it is desirable to consider discretisation methods that work on unstructured meshes. This brings the additional benefit that steep topographical changes can be modelled by unstructured meshes easily and accurately, which would cause severe distortion of structured meshes. In the context of dry atmospheric flows, methods based on unstructured meshes were fist considerer in \cite{AVH10,SSW13}.

While the equations modelling moist air are usually given in a non-conservative form, they can be written in the form of a non-linear hyperbolic conservation law system with source terms. This form of the equation lends itself to a discontinuous Galerkin (DG) discretisation, which we consider in this work. Additionaly, DG methods work both on structured and unstructured meshes. The system of equations is essentially the compressible Euler equations with multiple densities and some additional constraints and source terms modelling the phase changes of water. The different densities model the different water phases. While discontinuous Galerkin methods for the compressible Euler equations have been studied very well~\cite{BR97a,BO99,HH02a}, only a few approaches consider the system with moisture and even fewer include rain dynamics.
Among the (discontinuous) element-based approaches for atmospheric motion, we mention~\cite{TMQ22} in the context of the Nonhydrostatic Unified Model of the Atmosphere (NUMA), where both continuous and discontinuous elements are used on unstructured grids and rain dynamics are included, and~\cite{SBB14} in the context of the COSMO (Consortium for Small-scale Modelling) model. In the latter work, rain dynamics were not included. Furthermore, we mention~\cite{BGS20}, where a compatible finite element discretisation (partially discontinuous) was considered, based on the ENDGame model by the UK Met Office. While this work included rain dynamics, the discretisation relied on vertically structured grids. Finally, we mention ClimateMachine~\cite{STM22}, where discontinuous elements using vertical columns on Cartesian grids were considered and implemented both on CPUs and GPUs.

A particular challenge for discretisations is the realisation of the
source term, which models the phase change between water vapour and cloud
water, i.e., \emph{condensation}. The main cause of problems is that
the source term is only non-zero when the atmosphere is fully
saturated by water vapour. However, once the saturation threshold is
reached, the condensation source term is not given as an explicit
function of the other variables, but rather as the right-hand side
term balancing another transport equation. To avoid this, we take an
approach used, for example, in~\cite{STM22,DAB14,Ooy90}, and consider a
single density for cloud water and water vapour, such that the phase
change between these two does not have to be modelled. We then
reconstruct the individual densities after every time-step by solving an algebraic non-linear problem.

While individual aspects of our approach have been used in the literature, our combination and realisation of these techniques is novel and goes beyond what is considered in most cases in related works.
For example, while \cite{STM22} essentially consider the same set of equations and use a similar approach to deal with condensation, their work is based on structured grids, does not consider rain and only presents a convergence study for fourth-order polynomials and in a dry test-case over a very short time.
While \cite{TMQ22} also considered a DG approach, including unstructured meshes and rain, the equations used are different and simplified. In particular, the momentum and potential temperature only take dry mass into account. Furthermore, a constant artificial viscosity term is added, which modifies the problem being solved, condensation is implemented through a source term and correction scheme, and no convergence results for the scheme is presented in \cite{TMQ22}.
In \cite{GGD12}, the same set of equations with an artificial viscosity is used together with spectral elements. While higher-order polynomials up to order ten are used, this is only applied on structured quadrilateral meshes, and a diagonal mass matrix is only achieved through inexact integration. While the temporal convergence of the utilised IMEX scheme is studied, the spatial convergence of the scheme is not presented. A similar set of equations to those in \cite{TMQ22} is also used in \cite{SBB14}. A convergence study for linear and quadratic elements is shown in \cite{SBB14}. However, rain was not included in that work.
Finally, we mention \cite{BGS20}, where again, a different set of equations, compared to our work, is used. Convergence for the moist case is also shown here for the piecewise constant and linear cases, and rain is included in further examples. We also note that \cite{BGS20} require structured meshes and only use discontinuous elements for parts of the discretisation. However, this then leads to a compatible finite element method. 

The remainder of this work is structured as follows. In
\Cref{sec.eqn}, we present the equations under consideration in this
paper. This includes the microphysics parametrisation used and the
derivation of the implicit condensation equations. Our discontinuous
Galerkin approach is presented in \Cref{sec.discr}, and we present a
number of numerical examples taken from the literature and realised
with this approach in \Cref{sec.numex}.
This includes a high-order convergence study for multiple polynomial orders in the moist case without rain.
The results in
\Cref{sec.numex} are based on the open-source finite element library
NGSolve and are fully reproducible through openly accessible python
scripts implementing the presented examples. Details of the
thermodynamic constants are stated in \Cref{appendix:constants}, and
additional details on the numerical set-ups are given in
\Cref{appendix.setups}.

\section{Governing Equations}\label{sec.eqn}
The equations of motion for cloudy air with warm rain 
are given in a conservative form as
\begin{subequations}
\begin{align}
  \partial_t \rho_d + \nabla\cdot(\rho_d \ub) &= 0, \label{eqn.dry-density}\\
  \partial_t \rho_v + \nabla\cdot(\rho_v \ub) &= (S_{ev} - S_{cd}), \label{eqn.vapour-density}\\
  \partial_t \rho_c + \nabla\cdot(\rho_c \ub) &= (S_{cd} - S_{au} - S_{ac}), \label{eqn.cloud-density}\\
  \partial_t \rho_r + \nabla\cdot(\rho_r \ub - \rho_r v_r \eb_z) &= (S_{au} + S_{ac} - S_{ev}),
    \label{eqn.rain-density}\\
  \partial_t(\rho \ub) + 2 \bm{\Omega}\times\rho\ub +  \nabla\cdot(\rho \ub\otimes \ub - \rho_r v_r \ub\otimes \eb_z + \Id p)
    &= -\rho g \eb_z,\label{eqn.momentum}\\
    \partial_t E + \nabla\cdot((E + p) \ub - (c_l(T - \Tref) + \ub^2/2)\rho_rv_r \eb_z)
      &= - \rho g \eb_z \cdot \ub, \label{eqn.energy}
\end{align}
\end{subequations}
see, e.g.~\cite{CBvdH11,DAB14,Ooy01}. The variables $\rho_d, \rho_v,
\rho_c$, and $\rho_r$ are the densities of dry air, water vapour,
cloud water and rain, respectively. While $\ub$ is related to the
velocity of air, $v_r$ is the terminal rain speed, $T$ is
the temperature, $\Tref$ is the reference temperature and $E$ is the
sum of the internal and kinetic energy densities. 
The total density $\rho$ is given as the sum of the component densities
\begin{equation}
  \rho = \rho_d + \rho_v + \rho_c + \rho_r.
\end{equation}
The vector $\bm{\Omega}$ is the Coriolis vector due to the Earth's
rotation. As our simulations are on a time-scale where the effects from this
are negligible, we will not consider this term further below.
To close the system we have the pressure given by the equation of state
\begin{equation}
  p = (\rho_d R_d + \rho_v R_v) T,
\end{equation}
with the gas constants $R_d, R_v$ for dry air and water vapour,
respectively. The temperature can be
recovered from the equation for the sum of the internal and kinetic energy densities
\begin{equation}\label{eqn.mass-total-energy}
  E = (c_{vd}\rho_d + c_{vv}\rho_v + c_l(\rho_c + \rho_r))(T - \Tref)
    + \rho_v (\Lref - R_v\Tref)
    + \rho \frac{1}{2} \ub\cdot \ub,
\end{equation}
where $c_{vd}, c_{vv}$ are the specific heats of dry air and
water vapour, both at constant volume, respectively, while $c_l$ is the
specific heat of liquid water at the reference temperature. The constant
$\Lref$ is the latent heat of vaporization at the reference temperature.
The values used for these quantities are summarized in the appendix in
\Cref{tab.parameters}.
Finally, we note that while this form of the equations is rarely used 
for atmospheric modelling, it is the form that best lends itself to a 
discontinuous Galerkin discretisation. We refer to \cite{GR08} for an
overview of other common forms of the equations used in mesoscale
atmospheric modelling.

\subsection{Microphysics parametrisation}
\label{sec.eqn.subsec.micro-physics}
The source terms modelling the phase-changes of water are $S_{ev},
S_{cd}, S_{au}, S_{ac}$ which describe evaporation, condensation,
auto-conversion of cloud water into rain droplets and the collection
of cloud droplets by raindrops (accretion), respectively. The closure
of the moist dynamics is based on the microphysics closure due to
Kessler~\cite{Kes69}.

We take the specific form of the microphysics parametrisation from the
COSMO model, as described in~\cite[Section~5.4]{DFH21}. The evaporation,
auto-conversion and collecting-rain source terms are defined as
\begin{subequations}\label{eqn.source-terms}
\begin{align}
  S_{ev} &\coloneqq (3.86\times10^{-3} - 9.41\times10^{-5}(T - \Tref)) \left(1 + 9.1\rho_r^{3/16} \right)(\rho_{vs} - \rho_v)\rho_r^{1/2},\\
  S_{au} &\coloneqq 0.001\max\{\rho_c - q_{au}\rho, 0\},\\
  S_{ac} &\coloneqq 1.72 \rho_c \rho_r^{7 / 8},
\end{align}
\end{subequations}
with the auto-conversion threshold $q_{au}$, which is chosen as $q_{au} = 0$.

This leaves the condensation source term open.
Let us now denote by $q_i = \rho_i / \rho_d$, the mixing ratios of the
different water-phase densities with respect to the dry air density.
The condensation source term is defined implicitly by
\eqref{eqn.vapour-density} via two assumptions:
\begin{enumerate}[label=(\alph*)]
  \item Water vapour is saturated in the presence of cloud water, i.e.,
  \begin{equation*}
    \text{if }q_c > 0,\text{ then }q_v = q_{vs} \text{ and } S_{ev} = 0.
  \end{equation*}
  \item Cloud water evaporates instantaneously in undersaturated regions,
  i.e.,
  \begin{equation*}
    \text{if } q_v < q_{vs},\text{ then }q_c = 0 \text{ and } S_{cd}=0.
  \end{equation*}
\end{enumerate}
In oversaturated regions, the source term due to condensation
$S_{cd}$ is then defined through the relationship
\begin{equation*}
  \partial_t q_{vs} + u\cdot\nabla q_{vs} = -S_{cd},
\end{equation*}
where the saturation vapour mixing ratio is approximated in a saturated
atmosphere in terms of the saturation vapour pressure $e_s$ by
\begin{equation}\label{eqn.saturation-ratio}
  q_{vs} = \frac{\epsilon e_s(T)}{p - e_s(T)}.
\end{equation}
The saturation vapour pressure is recovered from the
Clausius–Clapeyron relation~\cite{CBvdH11}, and by using a linear
approximation of the latent heat of condensation under the assumption of
constant specific heats. This leads to the saturation vapour pressure
\begin{equation}\label{eqn.saturation-vapour-pre}
  e_s(T) = \eref \left(\frac{T}{\Tref}\right)^\frac{c_{pv} - c_l}{R_v}
    \exp \left[\frac{\Lref - (c_{pv} - c_l)\Tref }{R_v}\left(\frac{1}{\Tref}
      - \frac{1}{T}\right) \right],
\end{equation}
where $\eref\coloneqq e_s(\Tref)$ and $c_{pv}$ is the specific
heat of water vapour at constant pressure, see also~\cite{Rom08}.
The specific values used for these constants are again shown
in \Cref{tab.parameters}.

Finally, again following the warm rain scheme of the COSMO model, the
mean terminal velocity of rain is defined as
\begin{equation*}
  v_r = (\pi \rho_w N_0^r)^{-1/8} \frac{v_0^r\Gamma(4.5)}{6}\rho_r^{1/8},
\end{equation*}
with the distribution parameter set as $N_0^r = 
8\times10^6\,\unit{\per\meter\tothe{4}}$,
$v_0^r = 130\,\unit{\metre\tothe{1/2}\per\second}$ and the water density
$\rho_w = \rho_v + \rho_c + \rho_r$,~\cite{DFH21}.

\subsection{Implicit condensation} For the numerical method, we want all
the source terms to be given explicitly. To avoid the modelling of the
condensation source term, we follow~\cite{DAB14} and remove the
condensation source term by adding equations \eqref{eqn.vapour-density}
and \eqref{eqn.cloud-density}. This defines a new \emph{moist} density
$\rho_m = \rho_v + \rho_c$ of water which is only transported by the
velocity of dry air. This leads to the equation
\begin{equation}\label{eqn.moist-density2}
  \partial_t \rho_m + \nabla\cdot(\rho_m \ub) = S_{ev} - S_{au} - S_{ac}.
\end{equation}
To recover $\rho_c, \rho_v$, we then need to compute the vapour density
at saturation. This can be obtained from the saturation pressure
\eqref{eqn.saturation-vapour-pre} and the equation of state for the
vapour pressure $p_v = \rho_v R_v T$. Since this depends on the
temperature, we need to recover the temperature from the 
internal energy contribution in the energy equation \eqref{eqn.mass-total-energy}
, which in turn depends on all densities separately. We assume
that all cloud water evaporates instantaneously in undersaturated
regions. Using the moist density, the internal energy density
$\rho e = E - \rho \ub \cdot \ub / 2$, where $e$ denotes the internal
energy, and the saturation requirement, we then
arrive at the non-linear system
\begin{subequations}\label{eqn.implicit_condensation}
  \begin{align}
    \rho e &= (c_{vd}\rho_d + c_{vv}\rho_v + c_l(\rho_c + \rho_r))(T - \Tref) + \rho_v (\Lref - R_v\Tref),\\
    \rho_v &= \min\left(\frac{e_s(T)}{R_v T}, \rho_m\right),\\
    \rho_c &= \rho_m - \rho_v,
  \end{align}
\end{subequations}
which we can use to solve for the temperature, vapour density and
cloud density. The saturation vapour pressure $e_s(T)$ is defined in
\eqref{eqn.saturation-vapour-pre}. We note that while~\cite{DAB14}
also considered a single moist density, rain dynamics were not
included. For a detailed discussion of our solution strategy
of system \ref{eqn.implicit_condensation}, see \Cref{sec.discr:subsec.reconstruction} below.

\subsection{Perturbation formulation}
To facilitate the numerical approximation of the system
\eqref{eqn.dry-density}, \eqref{eqn.moist-density2},
\eqref{eqn.rain-density}, \eqref{eqn.momentum} and \eqref{eqn.energy},
we split the densities, pressure and energy density into the
hydrostatic part (time-independent and only dependent on the spatial
$z$-direction) plus a perturbation from the hydrostatic state:
\begin{subequations}
  \begin{align}
    \rho_i(x,y,z,t) &= \rhobar_i(z) + \rho_i'(x,y,z,t),  \qquad i \in \{d,v,c,m,r\},\\
    E (x,y,z,t) &= \overline{E}(z) + {E}'(x,y,z,t),\\
    p (x,y,z,t) &= \pbar(z) + p'(x,y,z,t),
  \end{align}
\end{subequations}
for which we have the relation
\begin{equation}
  \nabla \overline{p} = -\overline{\rho} g \eb_z.
\end{equation}
This then leads to the set of equations we use in our numerical method
\begin{subequations}\label{eqn.fullsysterm.perturbation}
\begin{align}
  \partial_t \rho_d' + \nabla\cdot((\rhobar_d + \rho_d') \ub) &= 0, \\
  \partial_t \rho_m' + \nabla\cdot((\rhobar_m + \rho_m') \ub) &= S_{ev} - S_{au} - S_{ac},\\
  \partial_t \rho_r' + \nabla\cdot((\rhobar_r + \rho_r')( \ub - v_r \eb_z)) &= S_{au} + S_{ac} - S_{ev},\\
  \partial_t(\rho \ub) + \nabla\cdot(\rho \ub\otimes \ub - (\rhobar_r + \rho_r') v_r \ub\otimes \eb_z + \Id p') &= -\rho' g \eb_z,\\
  \partial_t {E}' + \nabla \cdot ((\overline{E} + {E}' + \overline{p} + p') \ub - (c_l(T - \Tref) + \ub^2/2)(\rhobar_r + \rho_r')v_r \eb_z) &= - (\rhobar + \rho') g \eb_z \cdot \ub,
\end{align}
\end{subequations}
where the vapour, cloud densities and temperature (perturbations) are
again recovered via system~\eqref{eqn.implicit_condensation}.

\section{Discretisation}\label{sec.discr}

Our system of equations is a hyperbolic balance law, i.e., we may write it in short as
\begin{equation}\label{eqn.hyperbolic-equation}
  \partial_t U + \nabla\cdot F(U) = G(U),
\end{equation}
with $U = (\rho_d', \rho_m', \rho_r', \mb, E') \in [L^2(\O)]^{4+d}$, $d\in\{2,3\}$ the spatial dimension, and
$\mb=\rho \ub$, which we aim to solve using a discontinuous Galerkin scheme.

\subsection{Preliminaries and notation}
For the discretisation in space, let $\Th$ be a quasi-uniform and
shape-regular simplicial or hexahedral mesh of $\O$. 
In the case of a simplicial mesh, we define our discrete space as
\begin{equation*}
  W_h = [\PP^k(\Th)]^{4+d} \coloneqq \{V_h\in L^2(\Omega) : \restr{V_h}{K}\in\PP^k(K), \text{ for all } K \in\Th\}^{4+d},
\end{equation*}
and in the case of a hexahedral mesh, we set
\begin{equation*}
  W_h = [\QQ^k(\Th)]^{4+d} \coloneqq \{V_h\in L^2(\Omega) : \restr{V_h}{K}\in\QQ^k(K), \text{ for all } K \in\Th\}^{4+d}.
\end{equation*}
The dimension $4+d$ is due to the presence of four scalar variables
in \eqref{eqn.fullsysterm.perturbation} plus the number of momentum variables given by $d$.
The polynomial and tensor product spaces are
\begin{equation*}
  \PP^k = \spann\left\{\prod_{i=1}^d x_i^{\alpha_i} : \alpha_i \in\mathbb{N}_0, \sum_{i=1}^d {\alpha_i} \leq k\right\},
  \quad\text{and}\quad
  \QQ^k = \spann\left\{\prod_{i=1}^d x_i^{\alpha_i} : \alpha_i \in\mathbb{N}_0, 0\leq {\alpha_i} \leq k \right\},
\end{equation*}
for a given order $k\geq 1$. We note that in \texttt{NGSolve}, we use a basis of Legendre polynomials \cite{Zag06} leading to a diagonal mass-matrix.
 Furthermore, let $\Fh \coloneqq \{ F =
\partial K_1\cap \partial K_2 \;|\; \forall K_1,K_2\in\Th\}$ be the
set of facets of the mesh. We divide this into the set of boundary
facets $\Fhb \coloneqq \{F\in\Fh\;|\;F\subset\partial\O\}$ and
interior facets $\Fhi\coloneqq \Fh\setminus\Fhb$. On a facet
$F\in\Fh$, we define a fixed unit normal vector as
\begin{equation*}
  \nb = \begin{cases}
    \nb_{K_1} & \parbox[t]{11cm}{the unit normal vector to $F$ at $x$ pointing from $K_1$ to $K_2$ if $F\in\Fhi$ with  $F = \partial K_1\cap\partial K_2$; the orientation is arbitrary but fixed in what follows.}\\
    \nb & \parbox[t]{11cm}{the unit outward pointing normal to $\Omega$ at $x$ if $x\in\Fhb$.}
  \end{cases}
\end{equation*}
With the fixed unit normal vector and $K_1, K_2$ numbered accordingly, we define
\begin{equation*}
 U^- = \begin{cases}
    \restr{U}{K_1} & \text{if } F\in\Fhi,\\
    \restr{U}{K}   & \text{if } F\in\Fhb,
  \end{cases}
  \qquad
 U^+ = \begin{cases}
    \restr{U}{K_2} & \text{if } F\in\Fhi,\\
    g              & \text{if } F\in\Fhb,
  \end{cases}
\end{equation*}
for boundary data $g$ to be defined for a given example. We then define
the jump and average operators on a facet $F$ as
\begin{equation*}
  \jump{U} = U^- - U^+,
  \qquad\text{and}\qquad
  \mean{U} = \frac{1}{2}\left(U^-  + U^+\right).
\end{equation*}

\subsection{Spatial semi-discretisation}
The derivation of a DG formulation for \eqref{eqn.hyperbolic-equation}
can be found in many textbooks, e.g.,~\cite{HW08,PE12,Kop09,Gir20}. The idea is to
multiply \eqref{eqn.hyperbolic-equation} with an arbitrary test
function $V_h\in W_h$, integrate over $\Omega$, divide the integral
into element contributions, integrate by parts and choose a numerical
flux on element boundaries (since $U_h$ and $V_h$ are discontinuous).
This gives the spatially semi-discrete form
on every element $K\in\Th$
\begin{equation}\label{eqn.spatially-semi-discrete}
  \lprod{\partial_t U_h}{V_h}{K} - \lprod{F(U_h)}{\nabla V_h}{K} + \lprod{F_{n}(U_h)}{V_h}{\partial K} = \lprod{G(U_h)}{V_h}{K},
\end{equation}
where $\lprod{\cdot}{\cdot}{\omega}$ denotes the $L^2$-inner-product on $\omega \subset \Omega$.
For the numerical flux, we choose the Lax-Friedrich flux, which is given by
\begin{equation*}
  F_{n}(U_h) = \mean{F(U_h)}\nb + \frac{\Lambda}{2}\jump{U_h},
\end{equation*}
where $n$ is the outward-pointing unit normal vector on a given
element boundary and $\Lambda$ corresponds to the largest Eigenvalue
of $F'(U)$. We take
\begin{equation}\label{eqn.lambda_max}
  \Lambda = \max\{ \vert U_h^+\cdot \nb\vert + \vert v_r \eb_z\cdot \nb\vert + c_m, \vert U_h^-\cdot \nb\vert + \vert v_r \eb_z\cdot \nb\vert + c_m \},
\end{equation}
where $c_m$ is the speed of sound in moist air
\begin{equation*}
  c_m = \sqrt{\frac{\gamma_m p}{\rho}},
  \quad\text{with}\quad
  \gamma_m = \frac{q_d c_{vd} + q_v c_{vv} + (q_c + q_r) c_{l} + q_d R_d + q_v R_v}{q_d c_{vd} + q_v c_{vv} + (q_c + q_r) c_{l}}.
\end{equation*}
The ratio $\gamma_m$ is the isentropic expansion factor for moist air \cite{DAB14}, extended to include the specific heat and mass ratio for liquid water. In the case of $q_v=q_c=q_r=0$, this reduces to the isentropic expansion factor for dry air, and \eqref{eqn.lambda_max} reduces to the Lax-Friedrichs Flux for the dry Euler equations. Consequently, this can be viewed as the natural extension to the case with moisture and rain.
Note that \eqref{eqn.lambda_max} might not give the exact eigenvalue, however it is an
upper bound for the fastest wave and therefore sufficient, since
$\Lambda$ only needs to be large enough,~\cite{HW08}.
On domain boundaries, we will either
consider periodic or slip-wall (reflective) boundary conditions. Solid
walls are characterised by $\ub \cdot \nb=0$. To implement this, we set
$\mb_h^+ = (\mb^- - 2 (\mb_h^-\cdot \nb)\nb)$ on slip-wall boundaries and
$(U_h^+)_i = (U_h^-)_i$ for the remaining components.

We then arrive at the spatially semi-discrete form by summing over all elements
\begin{equation}\label{eqn.semi-discrete}
  \lprod{\partial_t U_h}{V_h}{\O} + F_h(U_h)(V_h) = G_h(U_h)(V_h),
\end{equation}
with
\begin{align*}
  F_h(U_h)(V_h) &\coloneqq \sum_{K\in\Th} - \lprod{F(U_h)}{\nabla V_h}{K} + \sum_{F\in\Fh} \lprod{F_{n}(U_h)}{\jump{V_h}}{F},\\
 G_h(U_h)(V_h) &\coloneqq \sum_{K\in\Th} \lprod{G(U_h)}{V_h}{K},
\end{align*}
where we define the jump $\jump{V_h} = V_h$ on domain boundary facets $F \in \Fhb$.

\subsection{Time-integration}
\label{sec.discr:subsec.time-discr}
Let $[0, T_\text{end}]$ be the time-interval of interest and consider
the constant time-step $\dt = T_\text{end} / N$ for some given
$N\in\mathbb{N}$. We will use an explicit time-stepping scheme to
advance the equations in time, since we know this results in a very efficient
scheme, in the sense that it is highly parallelizable, thus able to utilize HPC hardware. As usual with explicit schemes, the time-step needs to be
chosen sufficiently small, dependant on $h$ and $k$, to obtain stability. 
We note that IMEX schemes~\cite{GKC13,AGC17} are able to improve on the constant in the time-step restriction, by treating the terms leading to the fastest waves implicitly.

To this end, let $\{\varphi_i\}_{i=1}^{N_W}$ be the basis of $W_h$ and
$\ubu$ be the coefficient vector of $U_h$ in this basis. Furthermore,
let $\Mbu, \Fbu, \Gbu: \mathbb{R}^{N_W}\rightarrow \mathbb{R}^{N_W}$
denote the mass matrix flux and forcing operators, i.e.,
\begin{equation}
  \Mbu_{ij} = \lprod{\varphi_i}{\varphi_j}{\Omega}, \quad \Fbu(\ubu)_i
    = F_h(U_h)(\varphi_i), \quad\text{and}\quad \Gbu(\ubu)_i
    = G_h(U_h)(\varphi_i).
\end{equation}
Then the explicit Euler scheme for \eqref{eqn.semi-discrete} reads as: For
all $n=1,\dots, N$, find $\ubu$ such that
\begin{equation}\label{eqn.fully-discrete.EE}
  \Mbu \ubu^n = \Mbu \ubu^{n-1} - \dt \Fbu(\ubu^{n-1}) + \dt \Gbu(\ubu^{n-1}).
\end{equation}
Note that solving the above system can be implemented very efficiently
since (beside assembling several source vectors) it only involves
inverting the mass-matrix. As there is no coupling between different
elements, the latter is block diagonal (or even diagonal with the
correct choice of basis) and thus solving the system can be
implemented in a parallel manner and/or even matrix-free to better utilise
modern high-performance computing architectures. This is also illustrated in \Cref{sec.numex:subsec.gravitywaves:subsubsec.results} below.

\begin{rmk}
In practice, we will use the four stage, third order strong stability
preserving Runge-Kutta scheme (SPPRK(4,3)),~\cite{Kra91,GKS11}. To ease
the readability, we use the explicit Euler scheme in the presentation
of the method.
\end{rmk}

\subsection{Dependant variable reconstruction}
\label{sec.discr:subsec.reconstruction}

To solve system \eqref{eqn.fully-discrete.EE}, we need to evaluate the
discrete flux and source terms $F_h$ and $G_h$. Unfortunately, they
depend on the variables $\rho_v, \rho_c, T$ which are not present in the
(conservative) variables $U_h$, but derived from the non-linear
relationships~\eqref{eqn.implicit_condensation}.

After every time-step, we construct $\rho_{v,h}^n,\rho_{c,h}^n,T_h^n$
from $U_h^n$, by solving system \eqref{eqn.implicit_condensation} in
every quadrature point\footnote{NGSolve uses a product of
Gauss-Legendre rules, which in the case of simplicial elements are
mapped onto the simplex using a Duffy transformation.} by Newtons
method, taking the state at the last time-step as the initial guess.
The resulting function is then only well defined in the quadrature
points. The polynomial solution in $\PP^k(\Th)$ is then obtained via
an $L^2$-projection. As a result, we implicitly enforce the saturation
requirement which allows for both condensation of water vapour and
(instantaneous) evaporation of cloud water in undersaturated regions.
With these secondary, dependent variables reconstructed, we can then
assemble $F_h$ and $G_h$ for the next time-step or stage of the
Runge-Kutta scheme. Our approach is similar to that of the
one-step-scheme presented in~\cite{DAB14}, where
system~\eqref{eqn.implicit_condensation} was solved globally by an
iterative and decoupled Newton scheme. In contrast, we solve
\eqref{eqn.implicit_condensation} pointwise and fully coupled using a
standard Newton approach, and then reconstruct the piecewise
polynomial solution to evaluate the numerical flux. This has the
advantage of converging quadratically, and being easily computable in
parallel.

\subsection{Stabilisation through artificial diffusion}
Phase changes due to rain dynamics can lead to spurious oscillations
which eventually could have an effect on the stability of our method.
To dampen these oscillations and assure stability, we apply a proper
scaled artificial diffusion term to our system.
In the discontinuous Galerkin
setting, this is done through the symmetric interior penalty bilinear
form
\begin{equation}\label{eqn.sip-diffusion}
  A_h(U_h, V_h) = \begin{aligned}[t]
    & \sum_{K\in\Th}\lprod{\alpha\nabla U_h}{\nabla V_h}{K} \\
    &+  \sum_{F\in\Fh}\left[- \lprod{\alpha\mean{\nabla U_h} n}{\jump{V_h}}{F} - \lprod{\alpha\mean{\nabla V_h} n}{\jump{U_h}}{F} + \nicefrac{\sigma h^2}{k} \lprod{\alpha\jump{U_h}}{\jump{V_h}}{F}\right],
  \end{aligned}
\end{equation}
with the penalty parameter $\sigma >0$ chosen sufficiently large and a
vector-valued diffusion coefficient $\alpha>0$.

Convection dominated convection-diffusion problems can become unstable
in the case of a mesh Péclet number (Pe) larger than one. The mesh
Péclet number (for $\ub$-transported quantities) is given by
\begin{equation*}
  Pe = \frac{\Vert\ub\Vert_{\infty, K} h}{2\alpha}.
\end{equation*}
For pure convection problems (i.e. $\alpha = 0$), this can be seen as
infinitely large.
Motivated by this, we choose $\alpha$ in
\eqref{eqn.sip-diffusion} element-wise and component-wise for the
components of $U_h$ as
\begin{equation*}
  \alpha_K = \gamma 0.5 h^{1 - d/2} \Vert \ub \Vert_{0,K} \quad\text{and}\quad \alpha_K = \gamma 0.5 h^{1 - d/2}\Vert \ub - v_r \eb_z\Vert_{0, K},
\end{equation*}
where the latter is only chosen for the rain density $\rho_r$, since
this is the velocity-field by which it is transported. Note that we
use the element wise $L^2$-norm, omit the scaling by the mesh size and
include a scaling parameter $\gamma$. This is because we cannot easily
compute $\Vert \ub \Vert_{\infty, K}$ and because
$\Vert \ub \Vert_{0, K} \simeq c h^{d / 2} \Vert \ub \Vert_{\infty, K}$,
with a constant $c>0$ independent of $h$. For our numerical experiments,
we choose $\gamma$ as small as possible, without getting an unstable
method, and project the piecewise constant $\alpha_K$ into a
piecewise-linear, continuous function, such that $\alpha_K$ is well
defined on every facet. This is different to~\cite{TMQ22,GGD12}, where
the added artificial viscosity is scaled with a single constant
coefficient.

Finally, we note that for the case of the moist equations without
rain, i.e., without any explicit phase changes, no artificial
diffusion was found to be necessary for stability of the method. We
therefore only add diffusion for the examples with rain.

\begin{rmk}[Loss of higher-order convergence]
By adding the artificial diffusion term \eqref{eqn.sip-diffusion} to
our system, we are solving the original problem including a
perturbation of order $\mathcal{O}(h)$. As a result, we cannot expect
high order convergence of the DG method, even if higher-order elements
are used. For a fixed artificial viscosity (independent of $h$), we 
expect convergence to the exact solution of the perturbed system with the given fixed 
viscosity, but not towards the solution of the Euler equations (no viscosity)
under consideration here.
\end{rmk}

\subsection{Explicit Sponge Layer}
\label{sec.discr:subsec.sponge}

In some problems, the top of the domain should represent an open domain.
That is, waves should exit the domain and not reflect off the top boundary.
Consequently, a slip-boundary condition, representing a solid wall, is
unphysical. However, mathematical modelling of a non-reflective boundary is
a non-trivial task. In order to prevent waves from reflecting off the top
solid wall boundary, we implement a Rayleigh sponge layer to damp gravity
waves, as is common in the atmospheric flow literature
\cite{TMQ22, SBB14, GGD12, KG12, SSG16}.
The idea is essentially to relax the solution towards a known far-field 
condition (usually the hydrostatic background state).

In our code, we relax the moist mass-fraction perturbations, the velocity and
energy perturbation towards the known velocity and zero in the case of scalar
perturbations. In all the cases below, the velocity is also zero, so the
damping of the mass-fraction perturbations, velocity and energy perturbation
is equivalent to damping the density perturbations, momentum and energy 
density perturbation. That is, at the end of each time step, we compute
\begin{equation}\label{eqn.sponge}
  U_h^n = (1 - \delta_R(z))\widetilde{U}_h^{n, \text{pred}} + \delta_R(z)\widetilde{U}_h^\text{far},
\end{equation}
where $\widetilde{U}_h^{n, \text{pred}} = (\rho_m', \rho_r', \rho\ub, E')$ 
is the solution resulting from the time-stepping scheme, and 
$\widetilde{U}_h^\text{far}=(0,0,\rho\overline{\ub},0)$ is the far-field state 
for these variables. This is in contrast to \cite{SBB14, KG12}, where the 
damping is applied to the entire state vector (velocity and scalars). Since
the mass-faction damping can lead to mass loss in the sponge layer, we add a
correction step and locally add the moist mass lost through damping to the
dry density. Consequently, the total mass is still preserved globally by our
scheme. We note that mass fixing in numerical weather forecasting is done in
the literature \cite{DF13, MDAF19} with both local corrections and by
spreading the mass correction in the entire domain.

In \eqref{eqn.sponge}, we use the specific blending function
\begin{equation*}
  \delta_R(z) =
    \begin{cases}
      \frac{\alpha}{2}\left(1 - \cos\big(\pi \frac{z - z_D}{z_T - z_D}\big)\right) &\text{for }z \geq z_D,\\
      0 &\text{for }z<z_D.
    \end{cases}
\end{equation*}
Here $z_D$ is the height of the bottom of the sponge layer, $z_T$ is
the height of the top boundary and $0< \alpha \leq 1$ is a parameter to
tune the intensity of the sponge layer.

We implement \eqref{eqn.sponge} through an $L^2$-projection. The sponge
layer can also be implemented by adding a corresponding $L^2$-term to the
right-hand side, however, the choice of $\alpha$ is more difficult then.

\section{Numerical Examples}\label{sec.numex}

All numerical examples are implemented using the finite element library
Netgen/NGSolve~\cite{Sch97, Sch14} and \url{www.ngsolve.org}. We run
all examples either in a distributed memory (MPI) parallel or shared
memory parallel fashion (similar to openMP). The python scripts implementing our specific
examples and the results presented here are freely available on github
and archived on zenodo~\cite{HLSvW23_zenodo}.

We consider the following computation examples:
\begin{itemize}[labelwidth=6em, labelsep=1em, leftmargin=7em, rightmargin=1em, itemindent=!]
  \item[{\hyperref[sec.numex:subsec.gravitywaves]{Example 1:}}] Inertia gravity waves in a saturated atmosphere without rain dynamics to test our method with respect to optimal, high-order convergence, and parallel scalability.
  \item[{\hyperref[sec.numex:subsec.brianfristch]{Example 2:}}] The widely used benchmark by \citeauthor{BF02}~\citeyear{BF02} of a rising thermal in a saturated atmosphere, again without rain dynamics. We investigate the method with respect to large gradients due to the rising thermal.
  \item[{\hyperref[sec.numex:subsec.gravitywavesnoclouds]{Example 3:}}] We revisit the inertia gravity wave problem, but with an initial condition of a saturated vapour density but no clouds, in order to investigate cloud formation in our method.
  \item[{\hyperref[sec.numex:subsec.hydrostatic_mountain]{Example 4:}}] Conservation of a hydrostatic state in a domain with a large steep mountain. Unstructured meshes are of particular interest in this example due the steep terrain under consideration. This example also includes the possibility of rain forming, if the solution deviates from the hydrostatic state.
  \item[{\hyperref[sec.numex:subsec.rising_thermal_rain_2d]{Example 5:}}] A two-dimensional rising thermal leading to cloud formation and rain dynamics as formulated by \citeauthor{GC91}~\citeyear{GC91}.
  \item[{\hyperref[sec.numex:subsec.rising_thermal_rain_3d]{Example 6:}}] The three-dimensional extension of Example 5 as described in \cite{GC93}.
  \item[{\hyperref[sec.numex:subsec.squall_line]{Example 7:}}] Two dimensional squall line storm simulation based on \cite{GGD12,TMQ22}.
\end{itemize}

\subsection{Example 1: Inertia gravity waves in a saturated atmosphere}
\label{sec.numex:subsec.gravitywaves}
As our first example, we consider the case of a moist version of the
non-hydrostatics gravity waves adapted from~\cite{SK94} and presented
in~\cite{BGS20}. This example results in a smooth solution, such that
we can use it to test high-order convergence of our method.

\subsubsection{Set-up}
\label{sec.numex:subsec.gravitywaves:subsubsec.setup}
The domain is $\O=(0, 300\,\unit{km})\times(0, 10\,\unit{km})$. Periodic
boundary conditions are applied on the left and right boundaries, while
we have solid walls (slip) boundary conditions on the bottom and top
boundaries. The time interval under consideration is $[0, 3600\,\unit{s}]$.
The hydrostatic base state is defined via the wet equivalent potential
temperature
\begin{equation}\label{eqn.gravitywaves.hydrotheta_e}
  \overline{\theta}_e = \Theta_0 \exp(N^2 z / g),
\end{equation}
with $\Theta_0=300\,\unit{\kelvin}$ and
$N^2=10^{-1}\,\unit{\per\second\squared}$. The total water fraction is
$q_w=0.02$ and the pressure boundary condition is $p(0)=\pref$. The
perturbation applied to the hydrostatic base state is
\begin{equation*}
  \theta_e' = \frac{\Delta\Theta}{1 + a^{-2}(x - L/2)^2}\sin\left(\frac{\pi z}{H}\right),
\end{equation*}
where $H$ is the domain height, $L$ the domain length,
$a=5\times 10^{3}\,\unit{m}$ and $\Delta\Theta=0.01\,\unit{\kelvin}$.
The perturbation is applied under the requirements that $q_w$
and the pressure remain unchanged, and the air is saturated everywhere.
The initial velocity is prescribed as $u=(20, 0)^T
\,\unit{\meter\per\second}$. Additional details of the computation of the
hydrostatic base state and initial perturbation are provided in
\Cref{appendix.gravitywaves}.

\subsubsection{Results}
\label{sec.numex:subsec.gravitywaves:subsubsec.results}
\begin{figure}
  \centering
  \includegraphics{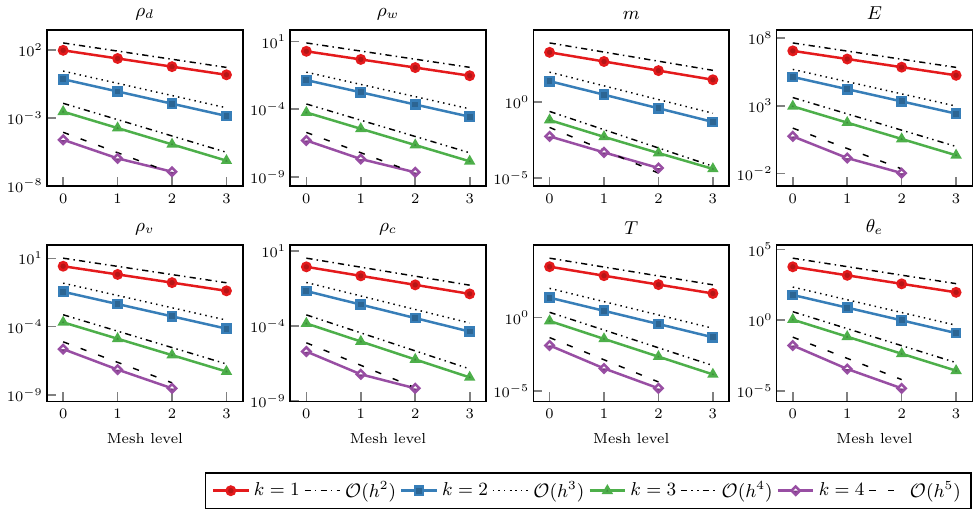}
  \caption{Example 1: Convergence results of the primary and secondary variables, and wet equivalent potential temperature on a series of structured quadrilateral meshes for different polynomial orders.}
  \label{fig.ex1_convergence}
\end{figure}

We consider this problem on a structured quadrilateral mesh, starting
with $h=1000\,\unit{m}$ and a series of three mesh refinements. On each
of these meshes, we consider the polynomial order $k=1,2,3,4$. The
time-step is chosen to be within the time-step restriction of the explicit
time-stepping scheme, that is $\dt=1\,\unit{s}$ for $h=1000\,\unit{m}$ with 
$k=1$, and $\dt=0.025\,\unit{s}$ for $h=125\,\unit{m}$ with $k=4$.
To compute convergence rates, we consider the
results on the finest mesh ($h=125\,\unit{m}$) with order $k=4$ as the
reference solution. We then use the $\ell^2(L^2)$-type
(space-time) norm
\begin{equation*}
  \Vert x_h - x_{h, \text{ref}}\Vert^2 = \sum_{i=1}^m \Vert x_h(t_i) - x_{h, \text{ref}}(t_i) \Vert^2_{L^2(\Omega)},
\end{equation*}
with a total of $m=30$ equidistant points in time. The convergence
results are shown in \Cref{fig.ex1_convergence}. We observe the
expected optimal order of convergence of $k+1$. There is some less
than optimal convergence between the finest two meshes considered with
$k=4$. We attribute this to the numerical reference solution not being
sufficiently accurate.

\begin{figure}
  \centering
  \begin{minipage}[b]{.49\textwidth}
    \centering
    \includegraphics{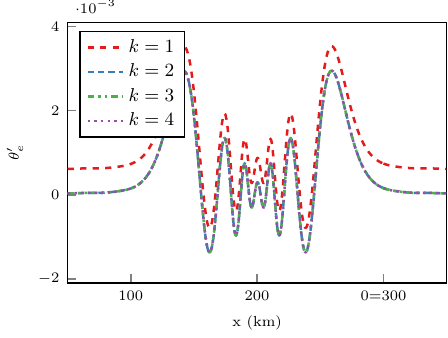}
  \end{minipage}
  \begin{minipage}[b]{.49\textwidth}
    \centering
    \includegraphics[height=4.45cm]{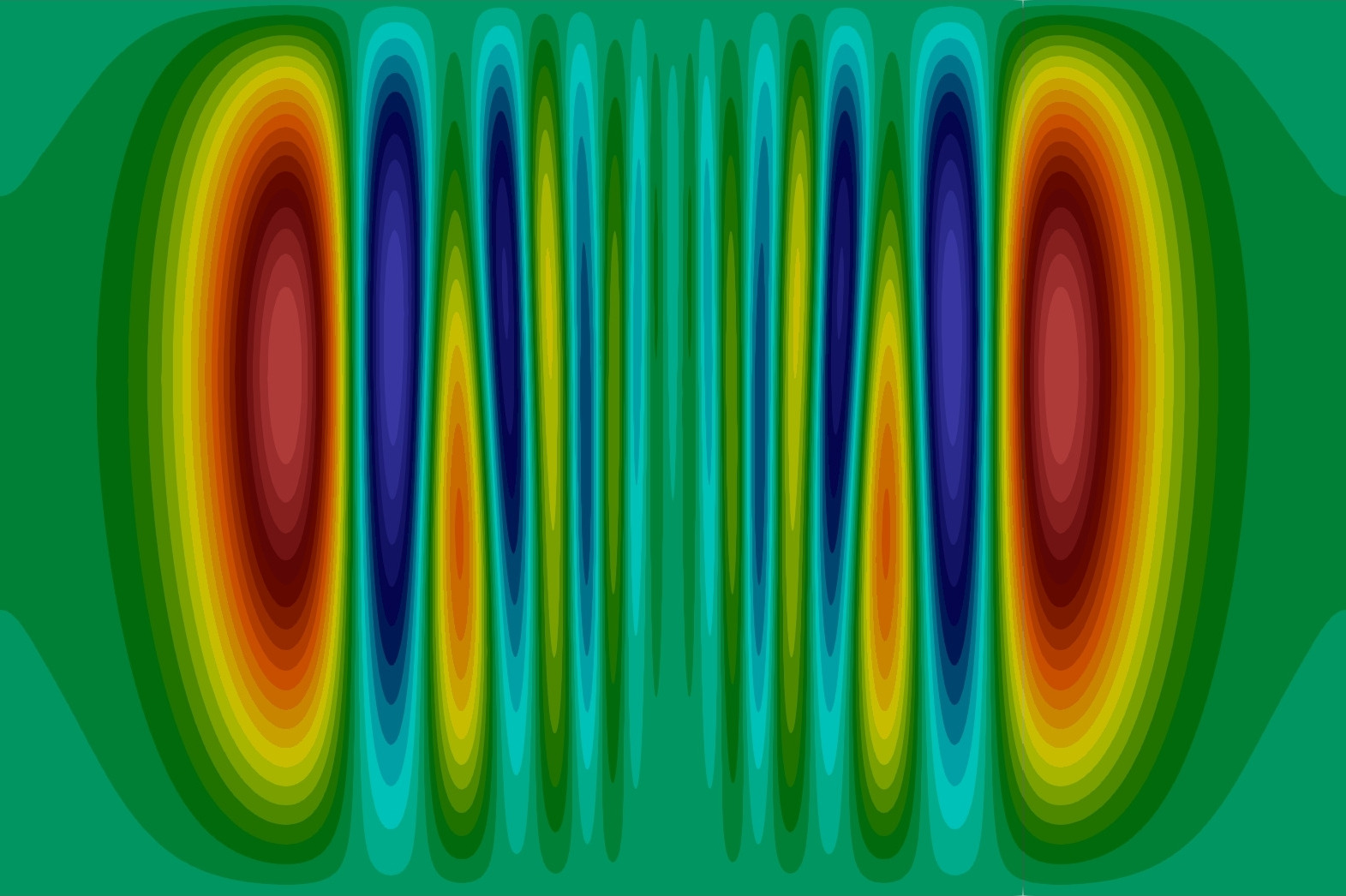}
    \includegraphics[height=4.45cm, trim={0 30cm 0 30cm}, clip=true]{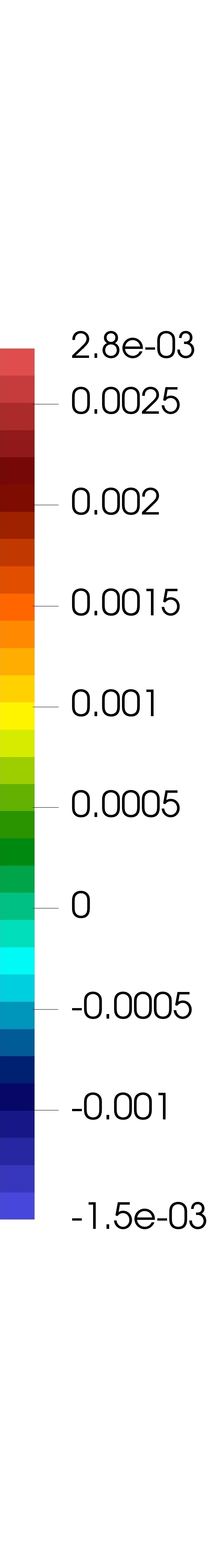}
    \vspace*{25pt}
  \end{minipage}
  \caption{Example 1: Inertia gravity waves in a saturated atmosphere. Perturbation of the wet equivalent potential temperature from the hydrostatic base state. Left: Profile along the line $z=5\,\unit{km}$ at $t=2520\,\unit{s}$ for different polynomial orders on a structured quadrilateral mesh with $h=125\,\unit{m}$, centred around $x=200\,\unit{km}$. Right: Solution at $t=3600\,\unit{s}$ computed using $\QQ^3$ elements with $h=500\,\unit{m}$, centred around $x=222\,\unit{km}$.}
  \label{fig.ex1_theta_e_pert}
\end{figure}

The perturbation of the wet equivalent potential temperature
along the line $z=5\,\unit{km}$ at $t=2520\,\unit{s}$ can be seen on the left of \Cref{fig.ex1_theta_e_pert}. Here we see the expected symmetry around the $x\approx 200\,\unit{km}$, similar to \cite{GR08}, even though the latter only considered dry dynamics. Due to our initial perturbation being located at $x=150\,\unit{km}$ as in \cite{BGS20} rather than at $x=100\,\unit{km}$ in \cite{GR08}, the axis of symmetry is also shifted in our results. Furthermore, we chose the time $t=2520\,\unit{s}$, rather than $t=2500\,\unit{s}$ as in \cite{GR08} as not all our time-steps hit this point in time.  We observe that the lowest order results with $k=1$ are visibly different from the higher order results, while for $k=2,3,4$, the results are indistinguishable. On the right of \Cref{fig.ex1_theta_e_pert}, we see the perturbation of the wet equivalent potential temperature at $t=3600\,\unit{s}$ in the bulk of the domain computed
using $\QQ^3$ elements on the mesh with $h=500\,\unit{m}$. This matches the results presented in \cite{BGS20} very well. However, we emphasize that this is not a variable in our system of equations, but a quantity obtained by post-processing.

\begin{figure}
  \centering
  \includegraphics{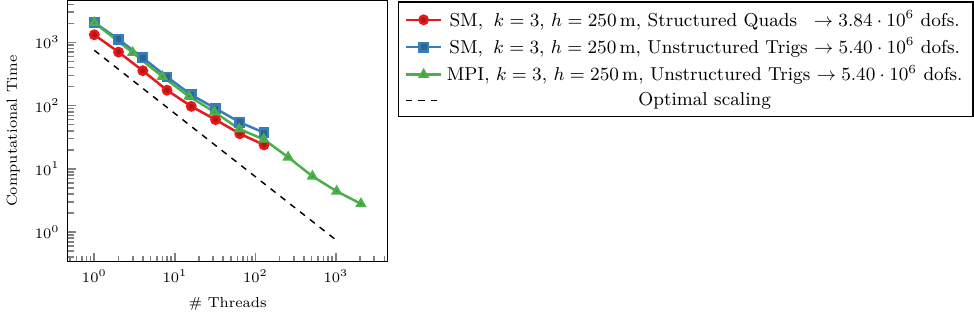}
  \caption{Example 1: Strong scaling results for the time-stepping loop using NGSolve's shared memory (SM) and distributed memory (MPI) parallel implementations. The number of degrees of freedom (dofs) are counted for the finite element space for the eight primary variables.}
  \label{fig.ex1_scaling}
\end{figure}

Finally, we consider the parallel scalability of our method. As discussed in \Cref{sec.discr:subsec.time-discr}, every time step requires vector assemblies, and the only large systems that have to be solved are mass-matrix problems. Due to the choice of basis, the matrices are diagonal and solving the system is implemented matrix-free, which is known to be a scalable approach. Our scheme's second major computationally expensive part is the dependant variable reconstruction, discussed in \Cref{sec.discr:subsec.reconstruction}. However, since this is done completely locally, we expect this to also scale well in parallel. A strong scaling plot for one hundred time steps is shown in \Cref{fig.ex1_scaling} using either NGSolve's own shared memory parallel loops or an MPI-distributed memory parallel version. As we can see, both versions scale well. We have nearly perfect scaling up to 16 parallel threads and good scaling even up to 2048 MPI ranks across 16 nodes.

\subsection{Example 2: Bryan-Fritsch moist benchmark}
\label{sec.numex:subsec.brianfristch}
Having established the optimal order error convergence of our method,
we consider a benchmark problem proposed in~\cite{BF02} as our second
example. This problem includes water vapour and clouds but still no rain.

\subsubsection{Set-up}
The domain is given by $\O=(0, 20\,\unit{km})\times(0, 10\,\unit{km})$. Periodic
boundary conditions are applied on the left and right boundaries and
rigid wall (slip) boundary conditions are applied on the top and bottom
boundaries, i.e., $\ub\cdot\nb = 0$. The time interval under
consideration is $[0, 1000\,\unit{s}]$. The hydrostatic base state is
defined by assuming a saturated atmosphere, $q_w = 0.02$, a constant
wet equivalent potential temperature of $\overline{\theta}_e =
320\,\unit{K}$ and the boundary condition for the pressure of $\pbar =
\pref = 10^5\,\unit{Pa}$. The initial perturbation is computed by
perturbing the density potential temperature in a circular region,
such that there is positive buoyancy in this region while the pressure
remains unchanged. Based on this, we compute the temperature
perturbation, from which we then reconstruct the vapour and cloud
densities using the saturation vapour pressure. This leads to the
initial condition $\rho_d' = 0$, $\rho_m' = 0$, $\rho u=(0, 0)$, $E'$,
$\rho_v'$, $\rho_c'$ and $T'$. Details on the computation of the
hydrostatic state and perturbation are given in
\Cref{appendix.bryanfritsch}

\subsubsection{Results}
\begin{figure}
  \centering
  \includegraphics[width=10cm]{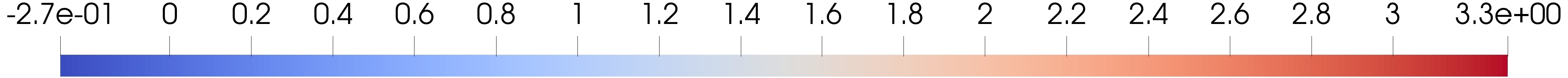}\\[2pt]
  \begin{subfigure}[b]{0.24\textwidth}\centering
    \includegraphics[width=4cm, trim={20.23cm 10.1cm 20.23cm 5.05cm}, clip=true]{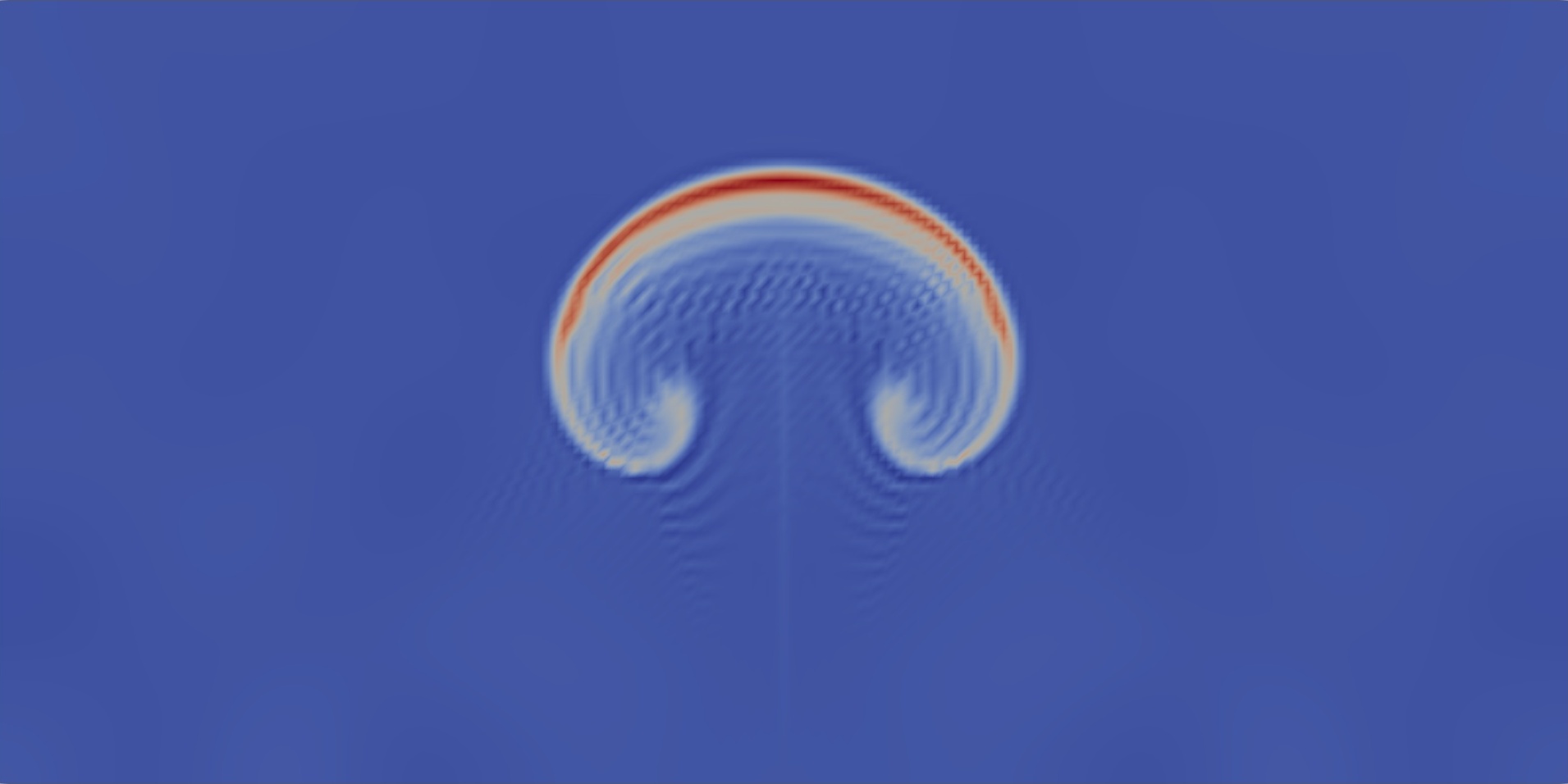}
    \caption{Structured quadrilateral mesh with $h\!=\!100\,\unit{m}$, $k\!=\!1$, $\dt\!=\!0.08\,\unit{s}$}
  \end{subfigure}
  \begin{subfigure}[b]{0.24\textwidth}\centering
    \includegraphics[width=4cm, trim={20.23cm 10.1cm 20.23cm 5.05cm}, clip=true]{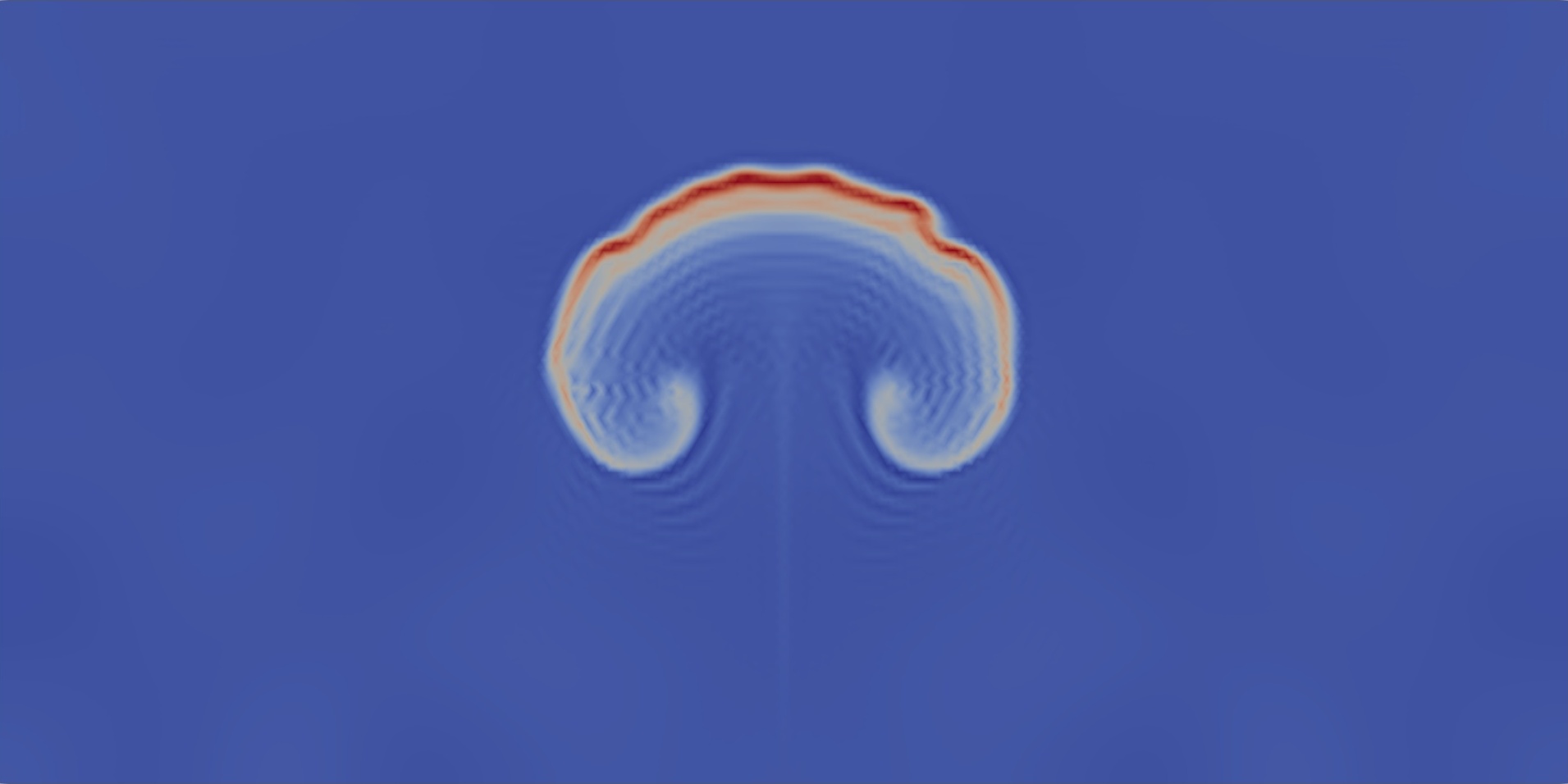}
    \caption{Unstructured simplicial mesh with $h\!=\!100\,\unit{m}$, $k\!=\!1$, $\dt\!=\!0.08\,\unit{s}$}
  \end{subfigure}
  \begin{subfigure}[b]{0.24\textwidth}\centering
    \includegraphics[width=4cm, trim={20.23cm 10.1cm 20.23cm 5.05cm}, clip=true]{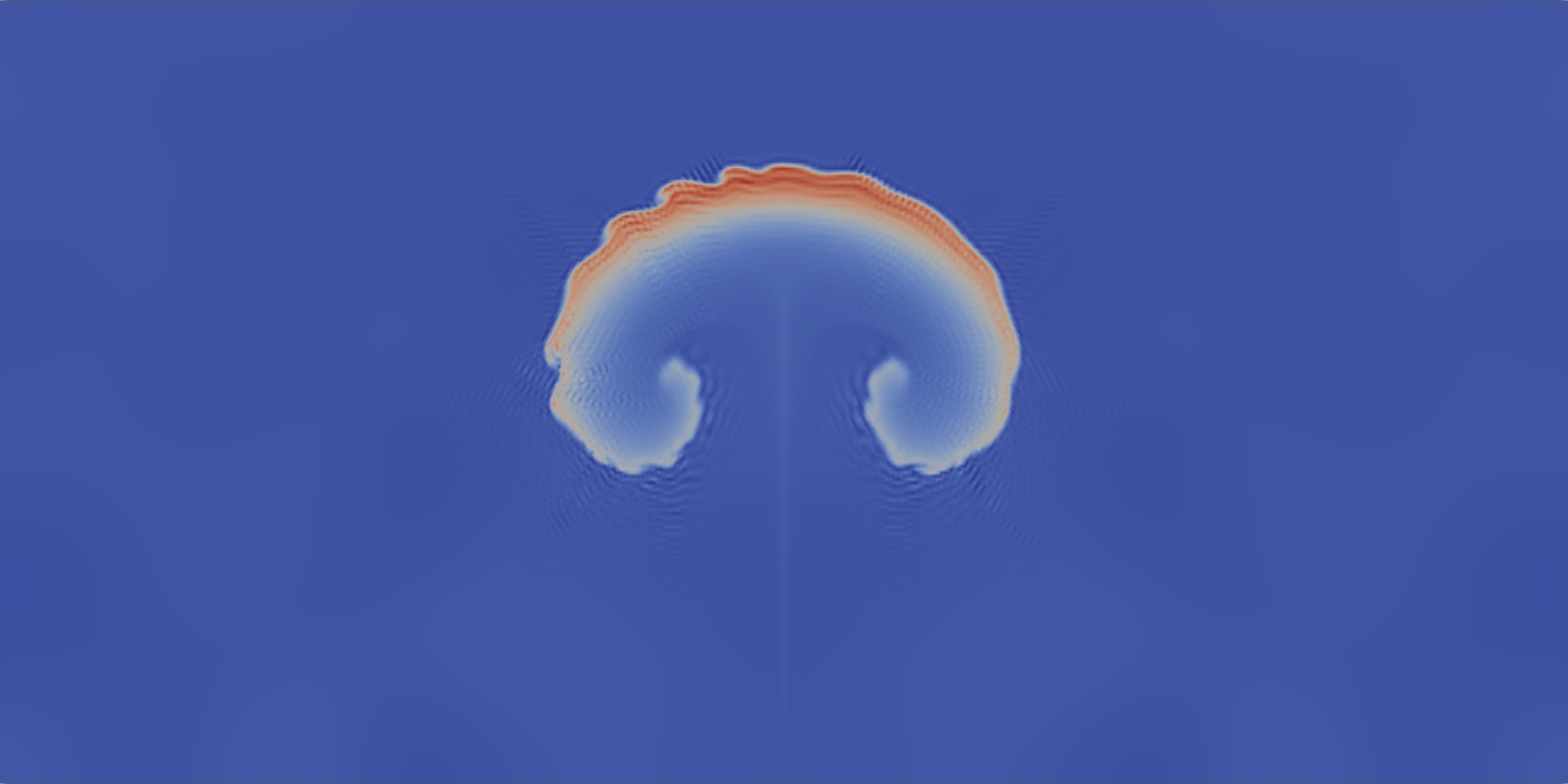}
    \caption{Unstructured simplicial mesh with $h\!=\!100\,\unit{m}$, $k\!=\!2$, $\dt\!=\!0.04\,\unit{s}$}
  \end{subfigure}
  \begin{subfigure}[b]{0.24\textwidth}\centering
    \includegraphics[width=4cm, trim={20.23cm 10.1cm 20.23cm 5.05cm}, clip=true]{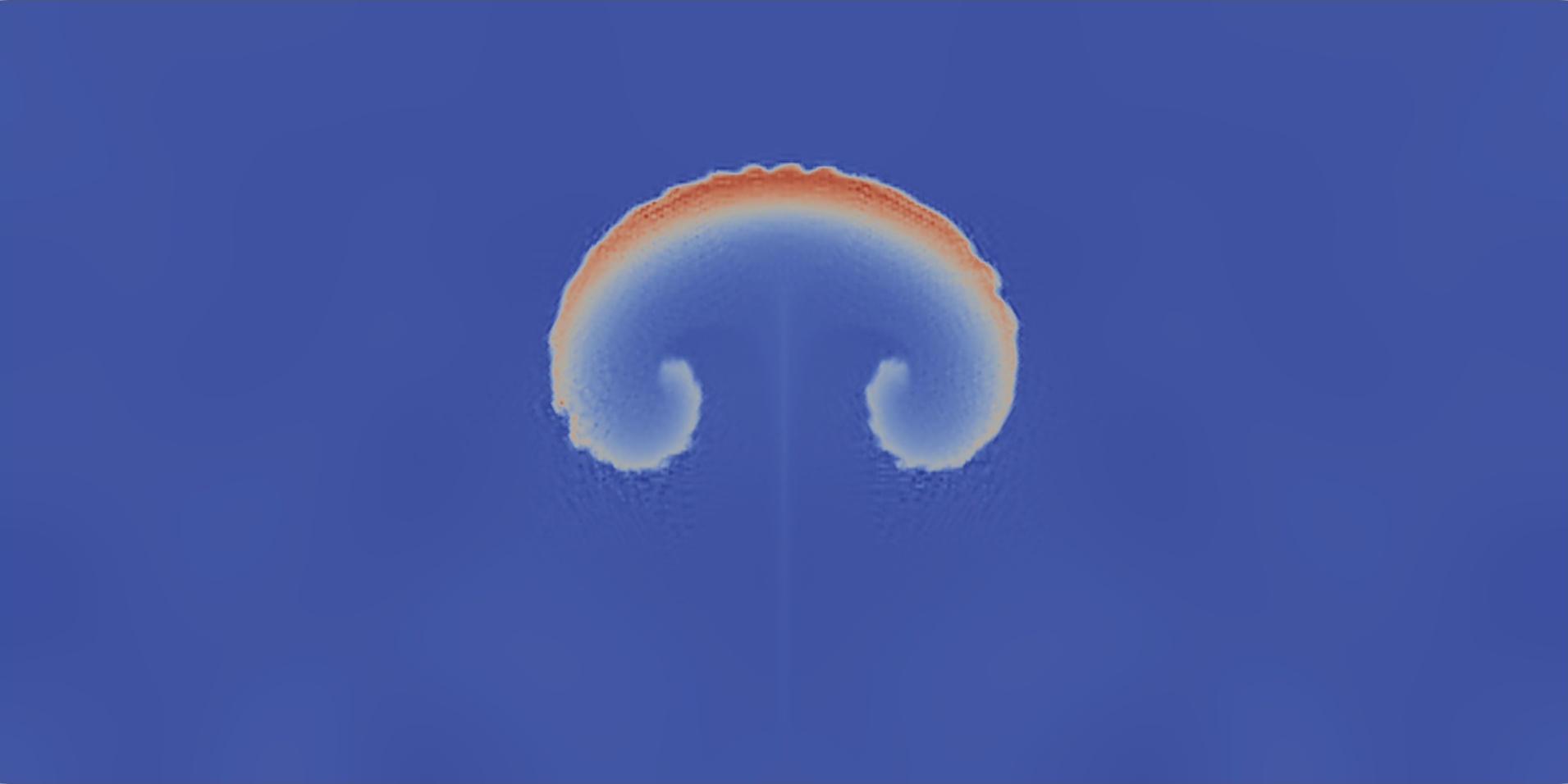}
    \caption{Unstructured simplicial mesh with $h\!=\!100\,\unit{m}$, $k\!=\!3$, $\dt\!=\!0.025\,\unit{s}$}
  \end{subfigure}
\begin{subfigure}[b]{0.24\textwidth}\centering
    \includegraphics[width=4cm, trim={20.23cm 10.1cm 20.23cm 5.05cm}, clip=true]{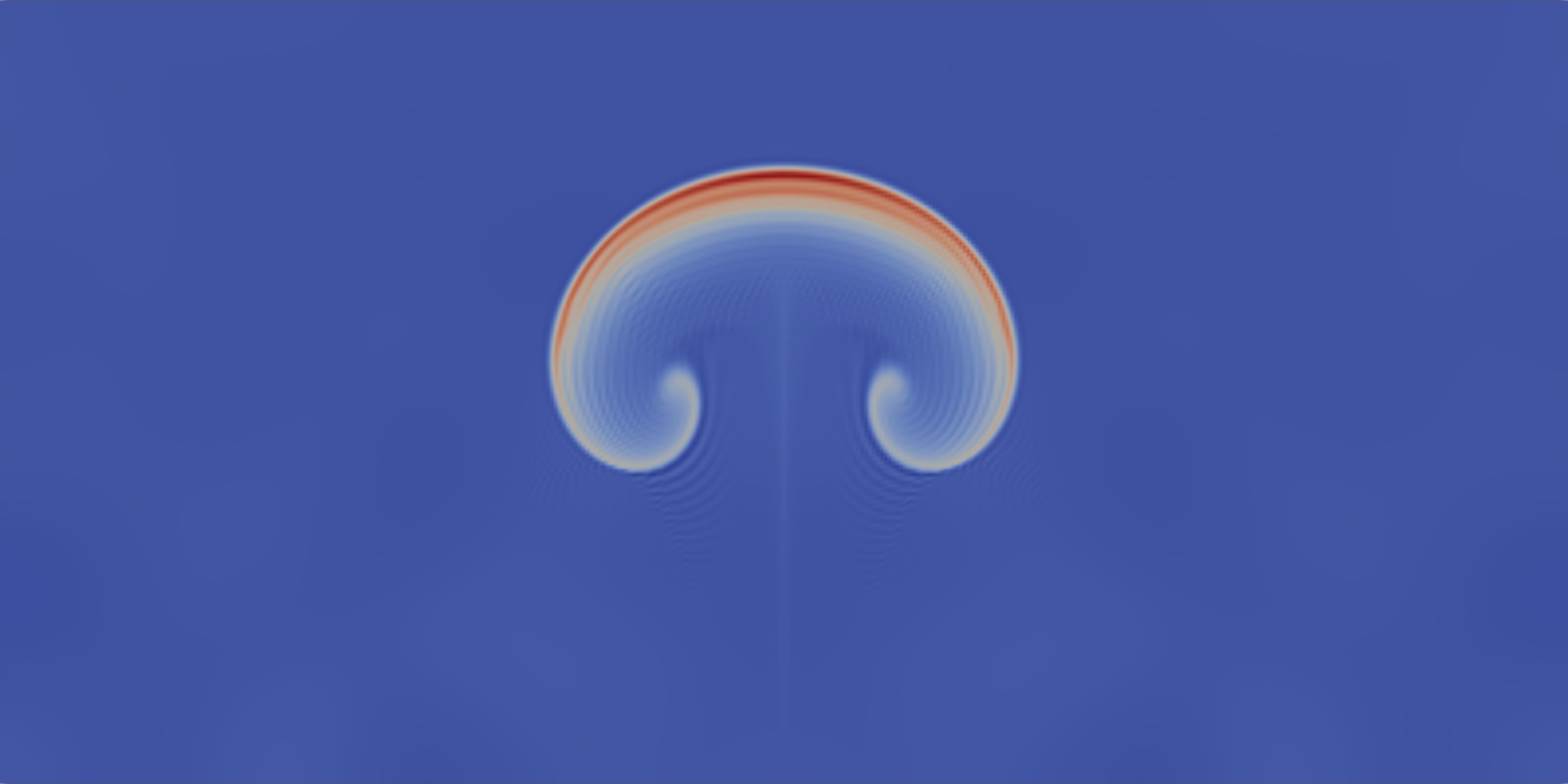}
  \caption{Structured quadrilateral mesh with $h\!=\!50\,\unit{m}$, $k\!=\!1$, $\dt\!=\!0.04\,\unit{s}$}
  \end{subfigure}
  \begin{subfigure}[b]{0.24\textwidth}\centering
    \includegraphics[width=4cm, trim={20.23cm 10.1cm 20.23cm 5.05cm}, clip=true]{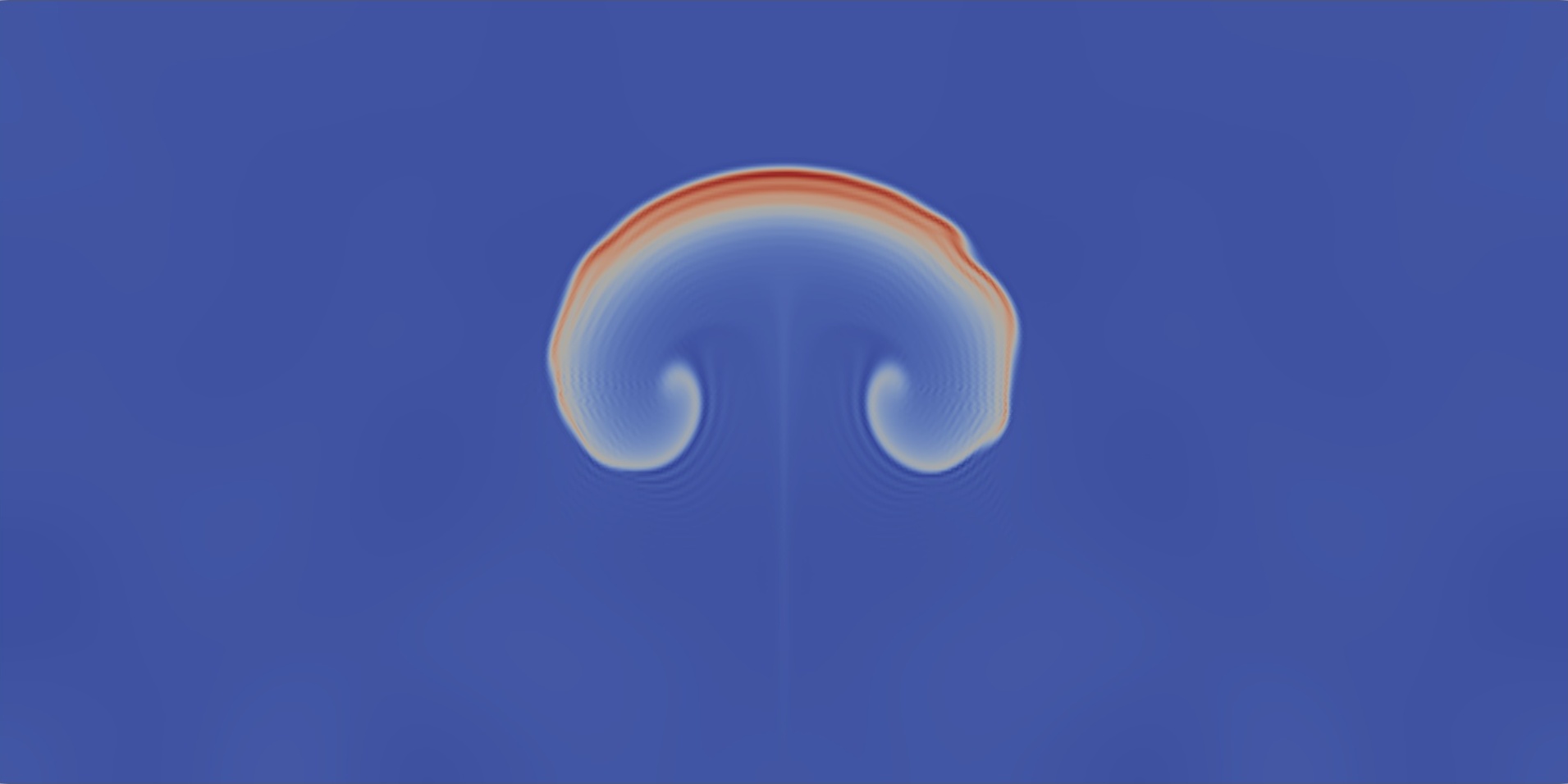}
    \caption{Unstructured simplicial mesh with $h\!=\!50\,\unit{m}$, $k\!=\!1$, $\dt\!=\!0.04\,\unit{s}$}
  \end{subfigure}
  \begin{subfigure}[b]{0.24\textwidth}\centering
    \includegraphics[width=4cm, trim={20.23cm 10.1cm 20.23cm 5.05cm}, clip=true]{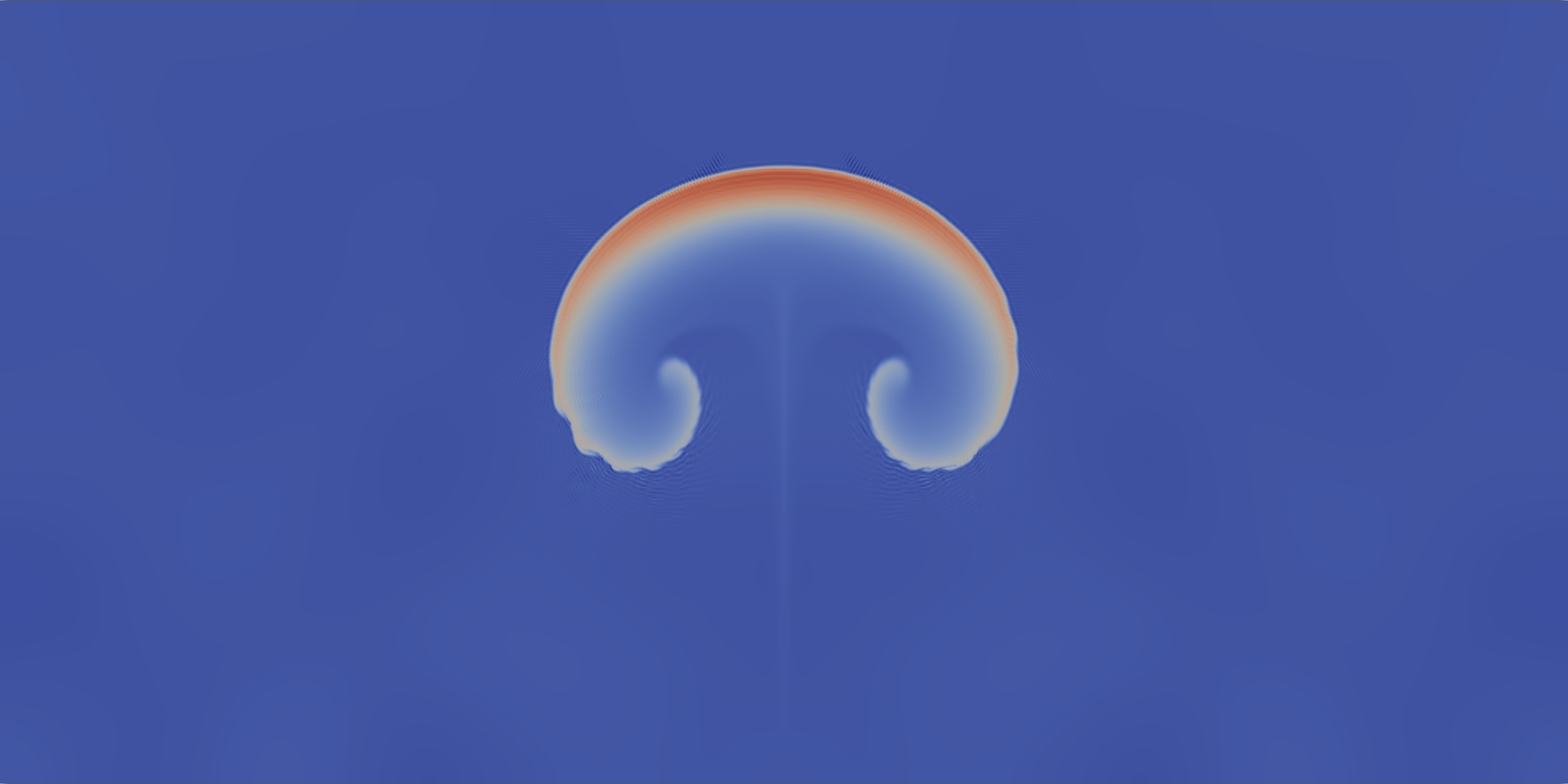}
    \caption{Unstructured simplicial mesh with $h\!=\!50\,\unit{m}$, $k\!=\!2$, $\dt\!=\!0.02\,\unit{s}$}
  \end{subfigure}
  \begin{subfigure}[b]{0.24\textwidth}\centering
    \includegraphics[width=4cm, trim={20.23cm 10.1cm 20.23cm 5.05cm}, clip=true]{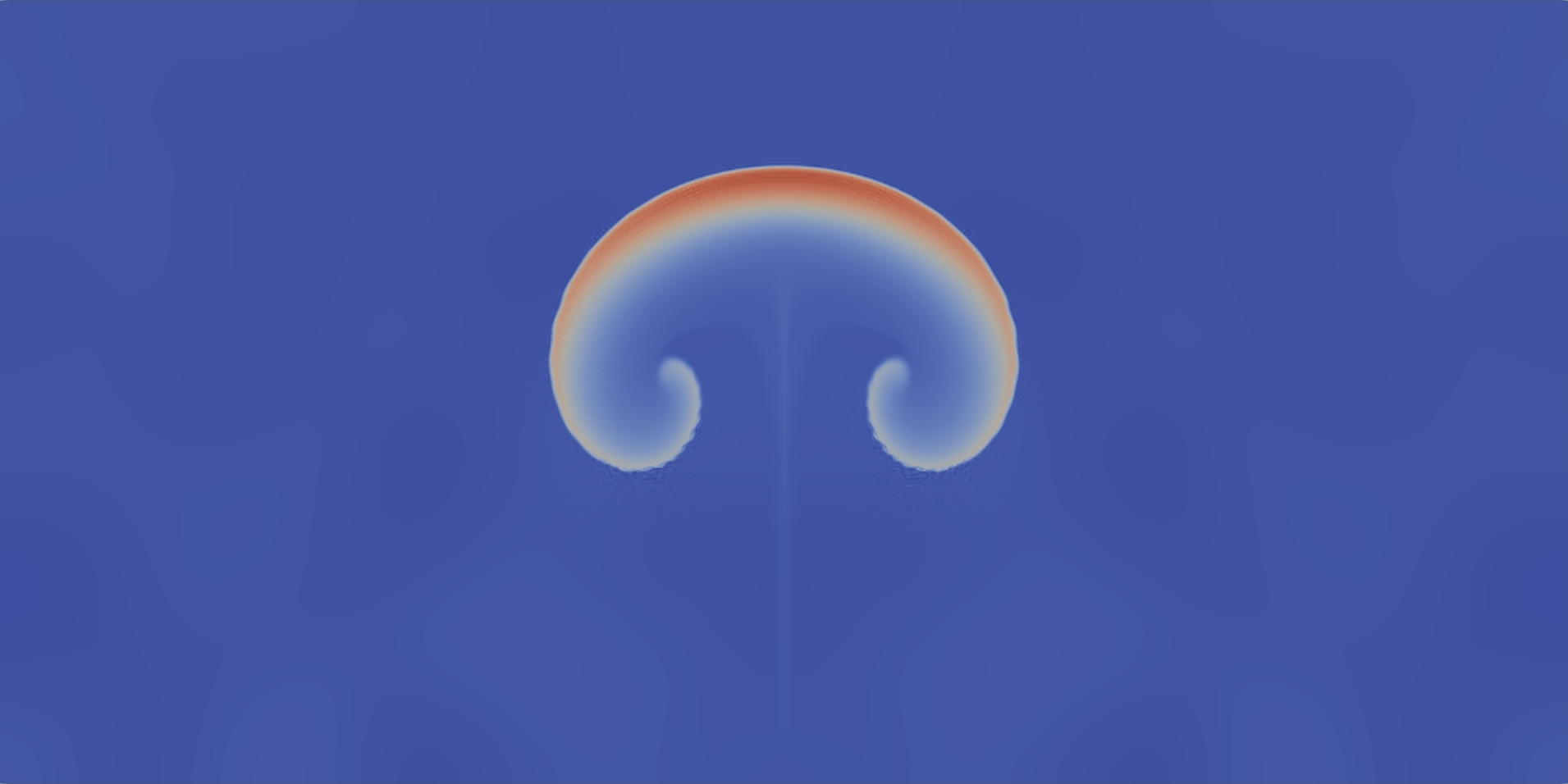}
    \caption{Unstructured simplicial mesh with $h\!=\!50\,\unit{m}$, $k\!=\!3$, $\dt\!=\!0.0125\,\unit{s}$}
  \end{subfigure}
  \caption{Example 2: Density potential temperature perturbation at $t=1000\,\unit{s}$ in the region $(6\,\unit{km}, 14\,\unit{km})\times(3\,\unit{km},8.5\,\unit{\km})$ on a series of different triangular and quadrilateral meshes using different order finite elements.}
  \label{fig.ex2.theta_rho_pert}
\end{figure}

The quantity of interest in this benchmark is the perturbation of the
density potential temperature $\theta_\rho$, the definition of which
is given below in \eqref{eqn.appendix.theta_rho}. The resulting perturbation
$\theta_\rho'$ over a series of meshes and polynomial orders can be
seen in \Cref{fig.ex2.theta_rho_pert}. The time-step is chosen to be
within the time-step restriction of the explicit time-stepping scheme.
The exact values are given in \Cref{fig.ex2.theta_rho_pert}.
We emphasize that the density potential temperature quantity is not a
variable we solve for in our system, but a quantity which has to be
post-processed from the available data.

The final states are comparable to that presented in
\cite{BF02}. There are some oscillations inside the temperature bubble
visible throughout, although we note
that they become smaller under mesh refinement and increasing
polynomial order. Unsurprisingly, the symmetry of the solution is
maintained on the structured quadrilateral mesh, and there are some
instabilities visible at the boundary of the temperature bubble for
higher-order elements. However, these also appear to become smaller
under mesh refinement and increasing polynomial order.

Finally, we note that for all polynomial orders considered, the
temperature bubble is consistent. This contrasts with the results in
\cite{BGS20}, where additional plumes appeared for their higher-order
case $k=1$. The authors of the latter paper attributed this to a physical
instability which is damped in the lowest order case $k=0$, due to higher
numerical diffusing.

\subsection{Example 3: Inertia gravity waves in a saturated atmosphere without initial clouds}
\label{sec.numex:subsec.gravitywavesnoclouds}
In the previous examples, the vapour density is saturated
throughout. To investigate cloud formation, we consider the inertia
gravity waves set-up, but define the initial condition such that vapour
density is fully saturated, but without any cloud water present.

\subsubsection{Set-up}
The spatial and temporal set-up is identical to that in
\Cref{sec.numex:subsec.gravitywaves:subsubsec.setup}. The hydrostatic
base state is computed using the same $\overline{\theta}_e$ given in
\eqref{eqn.gravitywaves.hydrotheta_e}, but now with $\rhobar_v =
\rhobar_{vs}$ and $\rhobar_c = 0$. The same perturbation is added to the
wet equivalent temperature and the initial state is again given by
assuming that the pressure is not changed by the perturbation and the
presence of a saturated atmosphere without clouds. Additional details
are provided in \Cref{appendix:subsec.gravitywavesnoclouds}.

\subsubsection{Results}
\begin{figure}
  \centering

  \begin{minipage}[b]{.49\textwidth}
    \centering
    \includegraphics{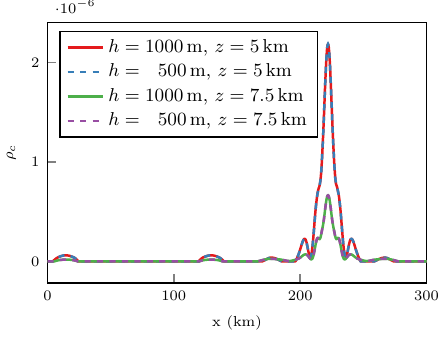}
  \end{minipage}
  \begin{minipage}[b]{.49\textwidth}
    \centering
    \includegraphics[height=4.5cm]{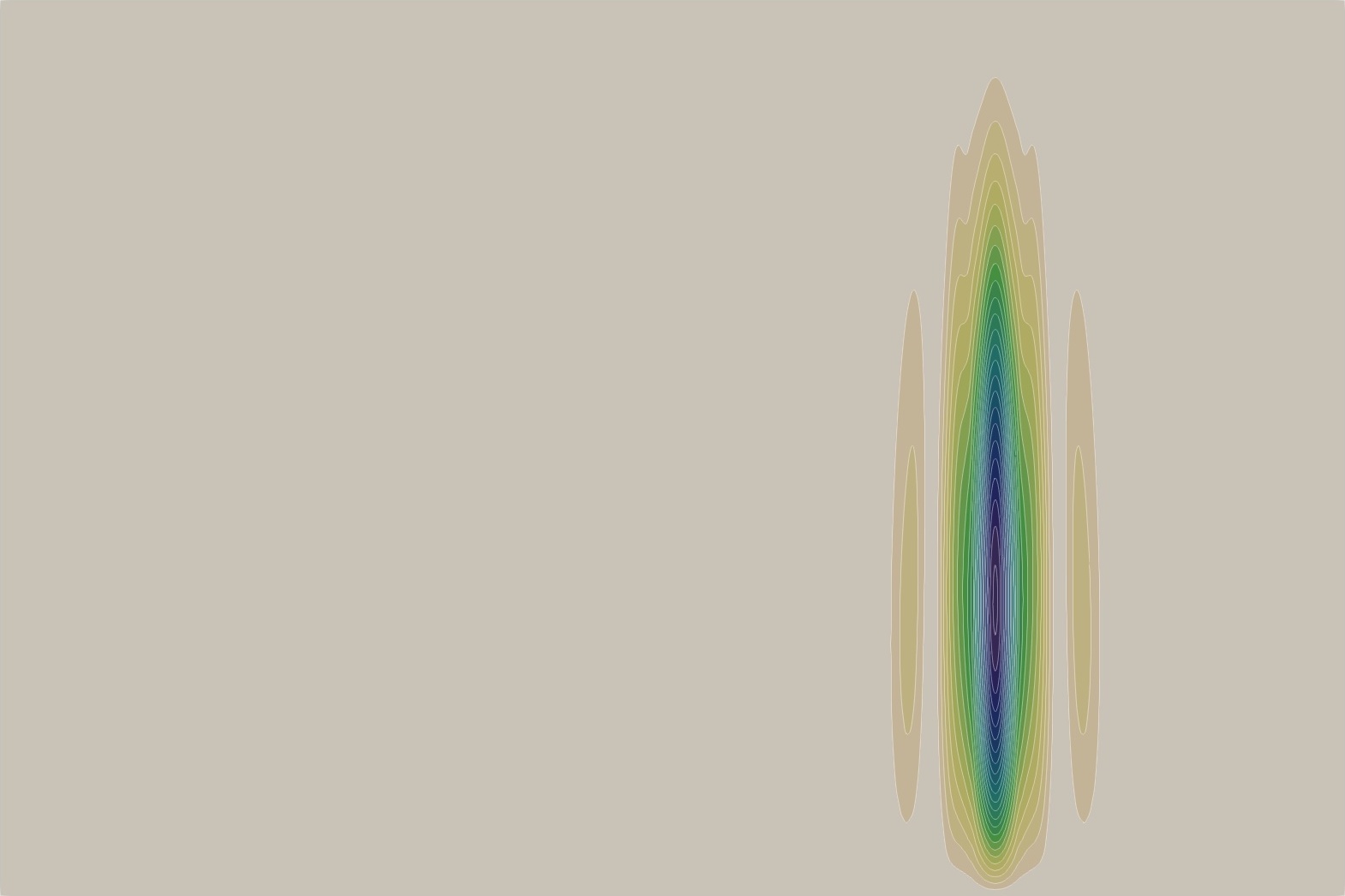}
    \includegraphics[height=4.5cm, trim={0 30cm 0 30cm}, clip=true]{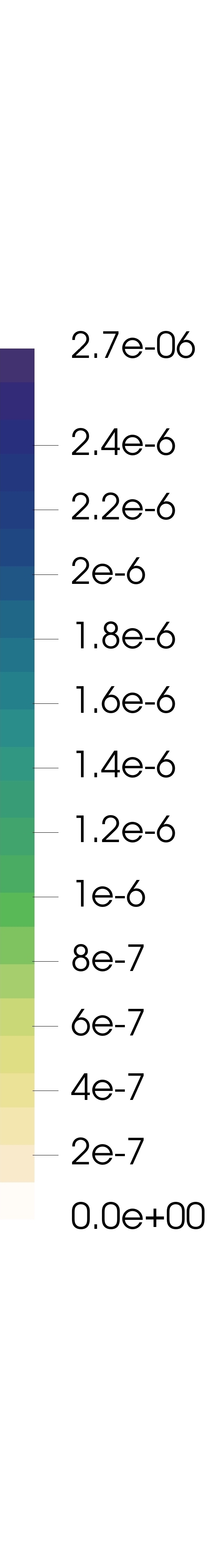}
    \vspace*{25pt}
  \end{minipage}
  \caption{Example 3: Cloud density (\,\unit{\kg\per\m^2}) and contours at $t=3600\,\unit{s}$ resulting from inertia gravity waves in a saturated atmosphere without initial clouds. Computed using $\QQ^3$ elements on a structured quadrilateral mesh. Left: Profiles along the lines $z=5\,\unit{km}$ and $7.5\,\unit{km}$. Right: Density in the whole domain, computed with $h=500\,\unit{m},\dt=0.15\,\unit{s}$.}
  \label{fig.ex3_clouds_results}
\end{figure}

We consider elements of order three on the two coarsest meshes used in
the previous example, since the errors in \Cref{sec.numex:subsec.gravitywaves} 
were already very small for this choice. The mesh sizes and time
steps are $(h,\dt)=(1000\,\unit{m}, 0.3\,\unit{s})$ and
$(h,\dt)=(500\,\unit{m},0.15\,\unit{s})$, resulting in $2.4\times10^5$ and
$9.6\times10^5$ degrees of freedom for the finite element space of the
primal variables, respectively\footnote{The compute wall times on a single
node with two AMD EPYC 7713 64-Core Processors with
hyper-threading enabled using 256 shared memory parallel threads was
$279\,\unit{s}$ ($h=1000, \dt=0.3$) and $1390\,\unit{s}$ ($h=500, \dt=0.15$) , respectively.}.
Profiles of the cloud density along two horizontal lines, 
and the density in the entire volume are shown in \Cref{fig.ex3_clouds_results}
for $t=3600\,\unit{s}$.
Looking at the results in \Cref{fig.ex3_clouds_results}, we see that 
the cloud formation is consistent between the two meshes, and in fact
there is no visible difference between the resulting cloud profiles.

\subsection{Example 4: Atmosphere at rest with a steep mountain}
\label{sec.numex:subsec.hydrostatic_mountain}

Approaches based on coordinate following discretisations can struggle in
cases of severe slopes, even with the atmosphere at
rest~\cite{SBB14,Zaen12}. As an example where unstructured meshes are
particularly attractive, we consider a domain with a single Gaußian
mountain, based on the dry example described in~\cite{Zaen12}, but
consider the full moist system \eqref{eqn.fullsysterm.perturbation}.

\subsubsection{Set-up}
The domain has a height of $40\,\unit{km}$ and a width of $35\,\unit{km}$.
The terrain profile is given by
\begin{equation*}
  z_m(x) = h_0 \exp\bigg(- \frac{(x - x_c)^2}{a^2}\bigg),
\end{equation*}
where $x_c$ is the centre of the domain and we consider the height
$h=7000\,\unit{m}$ and half-width $a=2000\,\unit{m}$, as in~\cite{Zaen12}.
This results in a maximal slope of about $3$. On the vertical
boundaries, we consider periodic boundary conditions and solid wall
boundaries at the top and bottom. 

For this, we need to include the sponge layer discussed in
\Cref{sec.discr:subsec.sponge}. Specifically, we choose 
$z_D=15\,\unit{km}$ and $\alpha = 0.1$.

We define the hydrostatic base state through the temperature profile
\begin{equation}\label{eqn.hydrostatic_mountain:Tprofile}
  \Tbar{}(z) = T_\text{str} + (T_\text{sl} - T_\text{str})\exp\left(-\frac{z}{H_\text{scal}}\right),
\end{equation}
where $T_\text{sl} = 288.15\,\unit{K}$, $T_\text{str} = 213.15\,\unit{K}$
and $H_\text{scal}=10000\,\unit{m}$. The water densities are then defined
by requiring vapour saturation throughout and setting the cloud and rain
densities to zero. Further details are provided in
\Cref{appendix.mountain} This setting is challenging, since deviations
from the hydrostatic state lead to dynamics in all variables. In
particular, deviations will cause cloud and rain formation,
c.f., \Cref{sec.eqn.subsec.micro-physics}.

\subsubsection{Results}

\begin{figure}
  \centering
  \begin{minipage}[t]{1.8cm}
    \centering
    \vspace*{-130pt}
    \includegraphics[height=4cm]{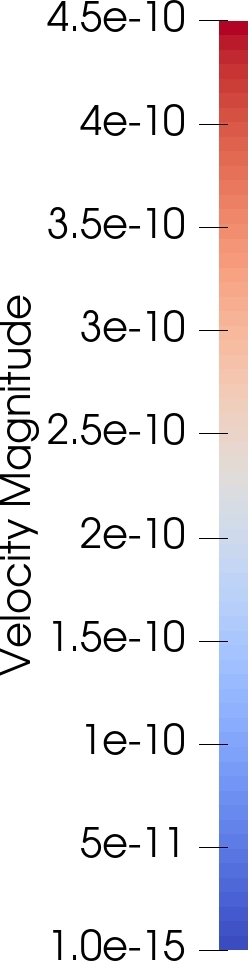}
  \end{minipage}
  \begin{minipage}[t]{9.5cm}
    \centering
    \includegraphics[height=5cm]{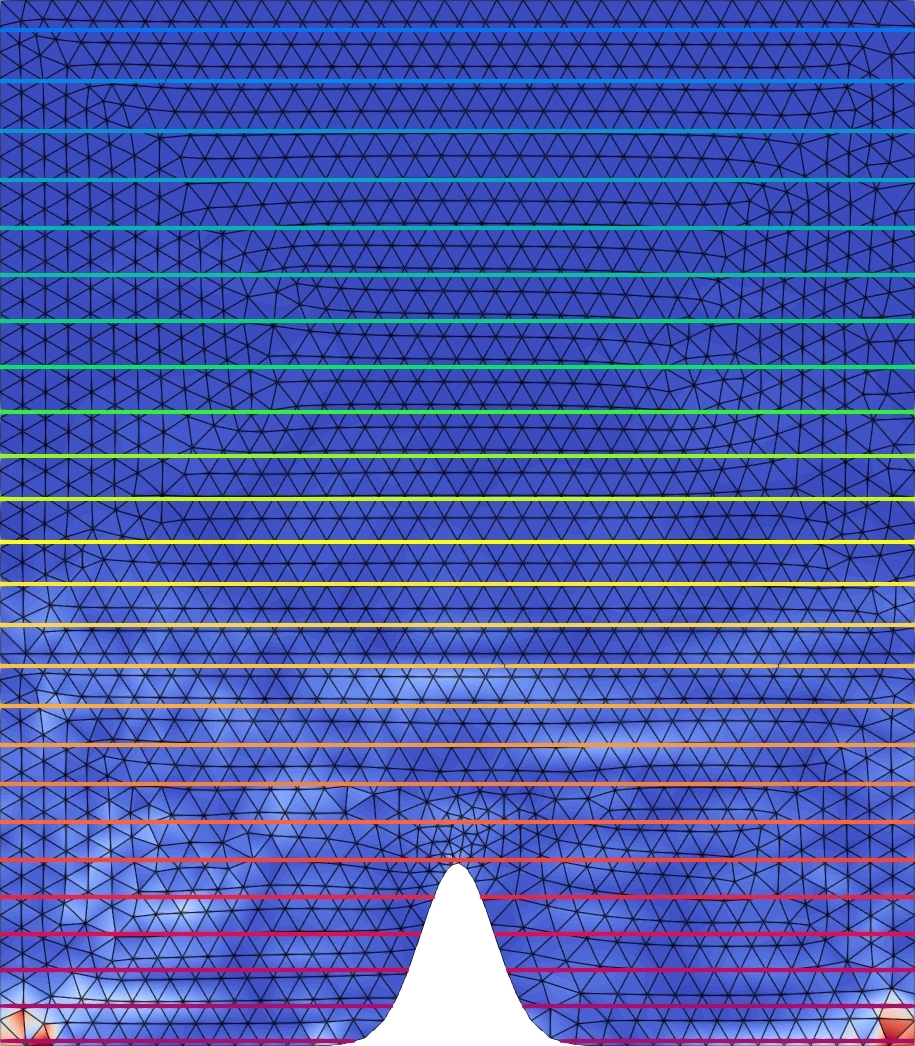}
    \hspace*{.5cm}
    \includegraphics[height=5cm]{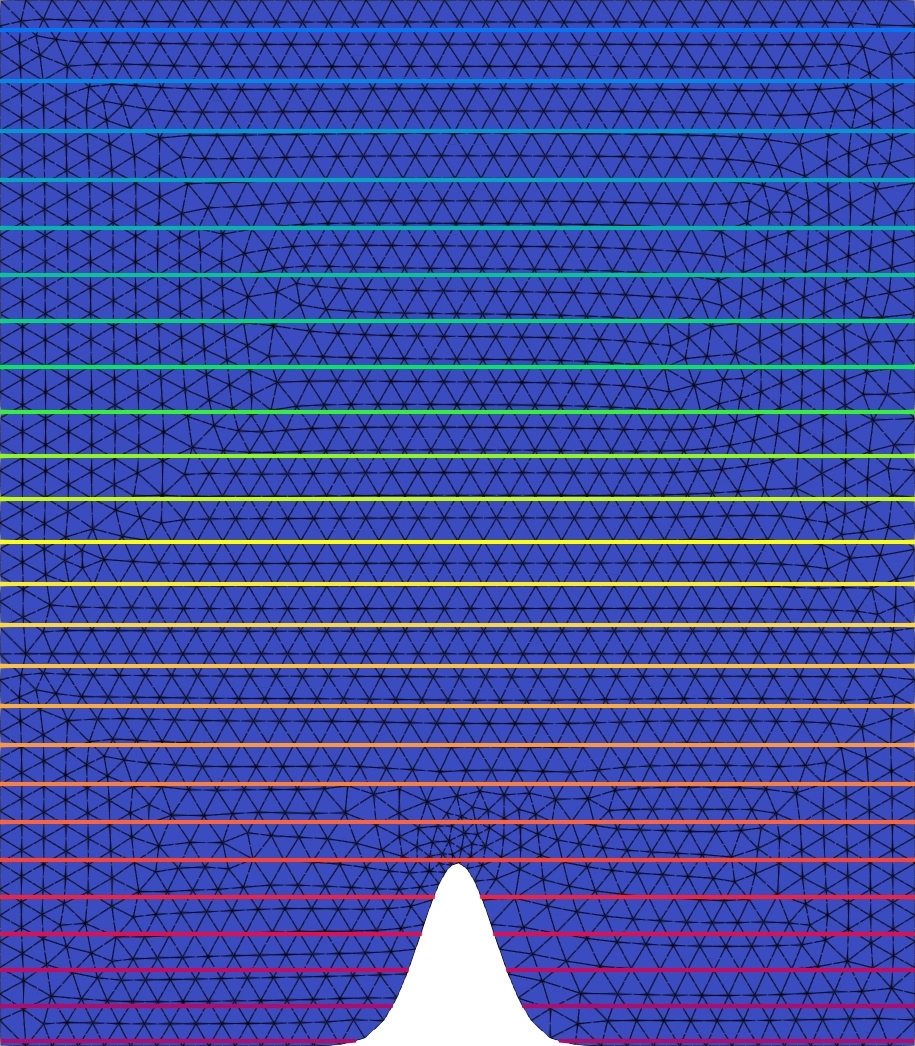}
  \end{minipage}
  \begin{minipage}[t]{1.8cm}
    \centering
    \vspace*{-130pt}
    \includegraphics[height=4cm]{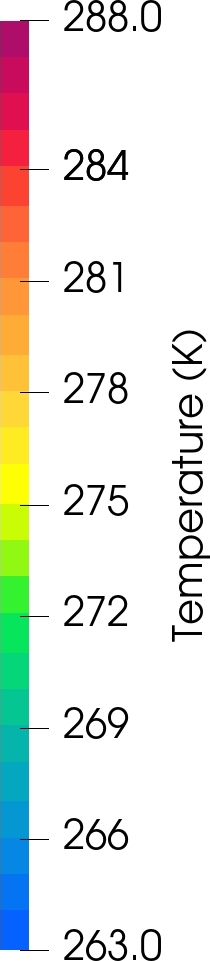}
  \end{minipage}
  \caption{Example 4: Results for the atmosphere-at-rest experiment with a steep mountain after $6\,\unit{h}$ on unstructured triangular meshes. Solid lines indicate temperature contour lines ($1\,\unit{K}$ intervals) and the shading indicates spurious velocity. Left: $h=1000\,\unit{\m}$, $k=1$, $\dt=0.4\,\unit{s}$, Right: $h=1000\,\unit{m}$, $k=2$, $\dt=0.2\,\unit{s}$.}
  \label{fig.ex4_results.steep_mountain}
\end{figure}

We consider an unstructured triangular mesh, with linear elements and 
mesh size $h=1000\,\unit{m}$, of the domain together with elements of
order $k=1,2$. The time step is chosen according to the stability limit,
resulting in $\dt=0.4\,\unit{s}$ and $0.2\,\unit{s}$, respectively. The 
resulting spurious velocity, temperature contours and the mesh can be seen in
\Cref{fig.ex4_results.steep_mountain}. Here we see that the solution has
remained stable, with the largest spurious velocities visible at the
coarsest elements and the bottom boundary. Furthermore,  there are no
oscillations visible in the temperature contours, even close to the
mountain profile. In fact, it appears that for the case $k=2$, no
spurious velocities are present. A closer inspection of the data shows
that the velocity is non-zero up to approximately $10^{-12}$, which is
not visible on the colour scale. Furthermore, we note that while the
simulation allowed for the formation of rain, which occurs as soon as
the cloud density is non-zero, see \eqref{eqn.source-terms}, the rain
density remained stable with a perturbation of $7\times 10^{-13}$
and $3\times10^{-15}$ for $k=1$ and $k=2$, respectively. Finally, while
our sponge layer allows for mass exchange between the moist and dry densities
in the sponge layer, the dry density perturbation observed is of order
$10^{-11}$ and $10^{-13}$ for $k=1$ and $k=2$, respectively, and with the
largest perturbations visible in the area below the sponge layer.

\subsection{Example 5: Rising thermal with rain in two spatial dimensions}
\label{sec.numex:subsec.rising_thermal_rain_2d}
This test problem is based on the one described in~\cite{GC91} and again
involves a rising thermal in two dimensions, but this time in an
undersaturated atmosphere and allowing for the development of rain.

\subsubsection{Set-up}
The domain is $\O=(0, 3.6\,\unit{km})\times(0, 2.4\,\unit{km})$, periodic
boundary conditions are set on the left and right boundaries and solid
wall (slip) boundary conditions at the top and bottom boundaries.
Note that while consequently dry air, water vapour and cloud water
cannot be transported out of the domain, rain water can leave the
domain since $v_r \nb_z$ is not necessarily zero at the bottom boundary.
The time interval under consideration is $[0,
600\,\unit{s}]$.

The hydrostatic base state is defined by specifying a relative humidity
and the dry potential temperature with
\begin{equation*}
  \overline{\HC}  = 0.2,\quad\text{and}\quad \overline{\theta}_d = \Theta e^{Sz},
\end{equation*}
where $\Theta$ is the dry potential temperature corresponding to
$T_\text{surf} = 283\,\unit{K}$ and $p = 8.5 \times 10^4\,\unit{Pa}$. This
pressure is also used for the pressure boundary condition at $z=0$. The
stratification is given by $S = 1.3 \times 10^{-5}\,\unit{\per\meter}$.
There are no clouds and no rain. Details on the computation of the
hydrostatic base state are given in
\Cref{appendix:subsec.risingthermalrain}.

The initial perturbation of the initial condition is then given by a
circular bubble, where the air is saturated but no clouds are present,
without changing the pressure and dry potential temperature. Further
details on the initial condition are again provided in
\Cref{appendix:subsec.risingthermalrain}.

\subsubsection{Results}
We consider an unstructured simplicial mesh of the domain with
$h=25\,\unit{m}$ and orders $k=1,2$, as well as a second mesh with
$h=12.5\,\unit{m}$ and $k=2$. The time step is chosen as large as the
time-step restriction allows, resulting in $\dt=0.02\,\unit{s}, 0.01\,\unit{s}$
and $0.005\,\unit{s}$, respectively. For stability, we include artificial
diffusion here and choose the parameter to be $\gamma=0.06$, which we
found to be the lower limit to preserve stability in this case.

\begin{figure}[p]
  \centering
  \includegraphics[height=3.3cm]{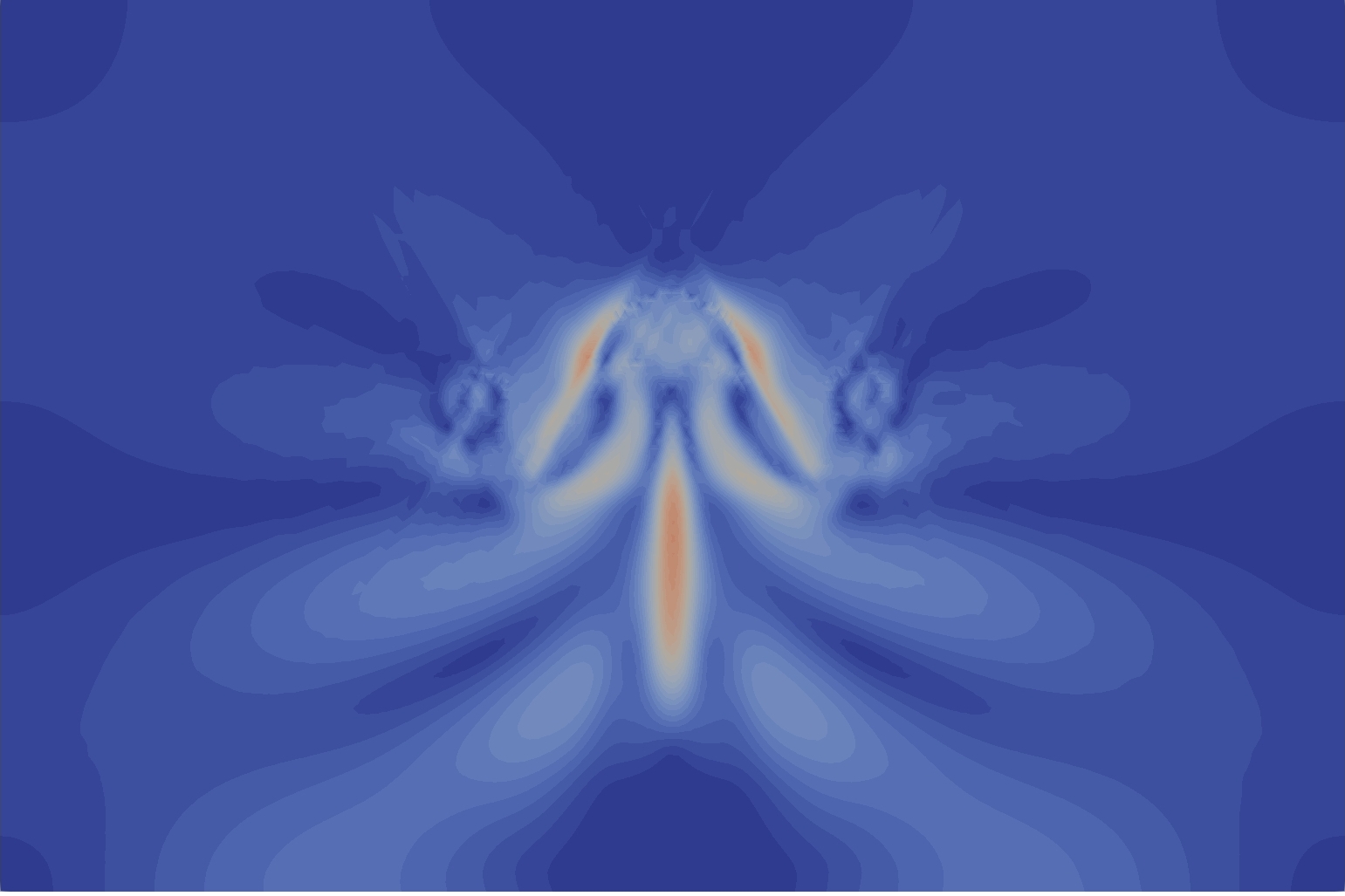}
  \includegraphics[height=3.3cm]{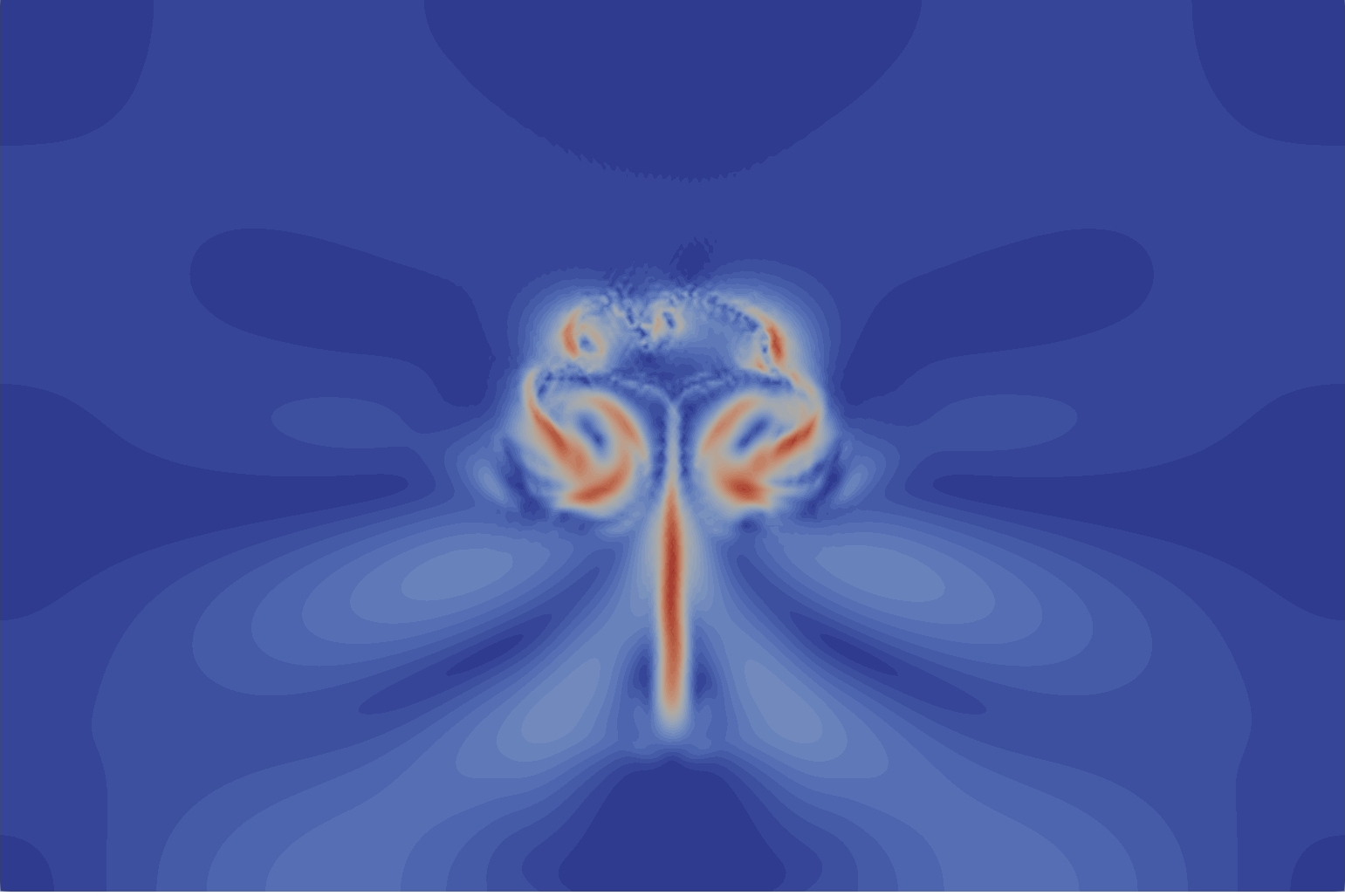}
  \includegraphics[height=3.3cm]{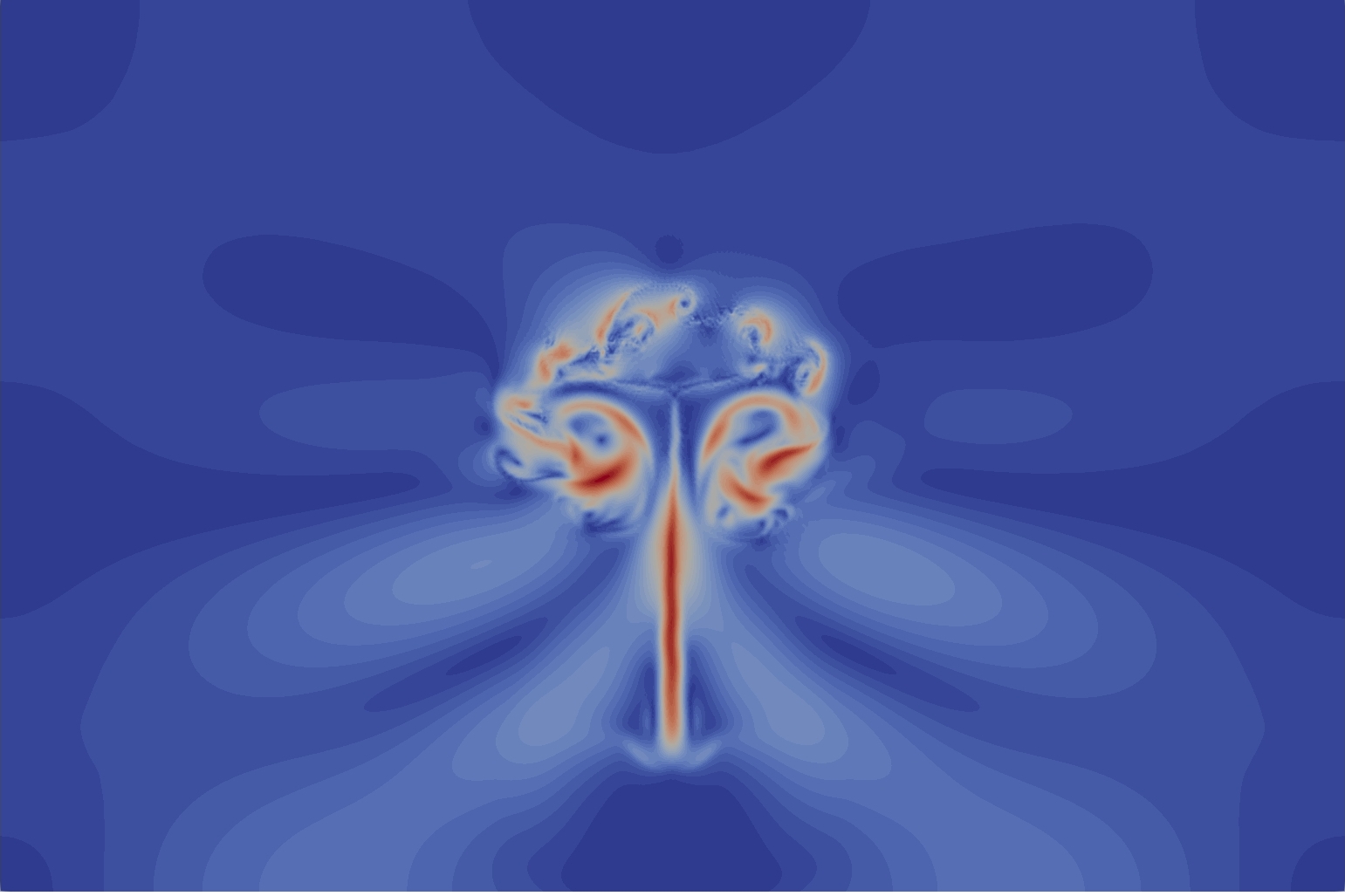}
  \includegraphics[height=3.3cm]{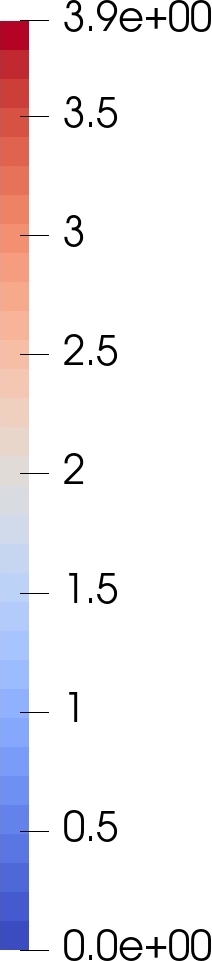}
  \caption{Example 5: Dry air velocity (\,\unit{\m\per\s}) solution at $t=600\,\unit{s}$ for a two-dimensional rising thermal leading to precipitation. Discretisation parameters from left to right: $(h,k,\dt)=(25\,\unit{m},1,0.02\,\unit{s}), (25\,\unit{m},2, 0.01\,\unit{s}), (12.5\,\unit{m},2, 0.005\,\unit{s})$.
  Computed on unstructured simplicial meshes.}
  \label{fig.ex5_results.rising_thermal_rain_2d.velocity}
\end{figure}

\begin{figure}
  \centering
  \begin{minipage}[t]{15cm}
    \includegraphics[height=3.3cm]{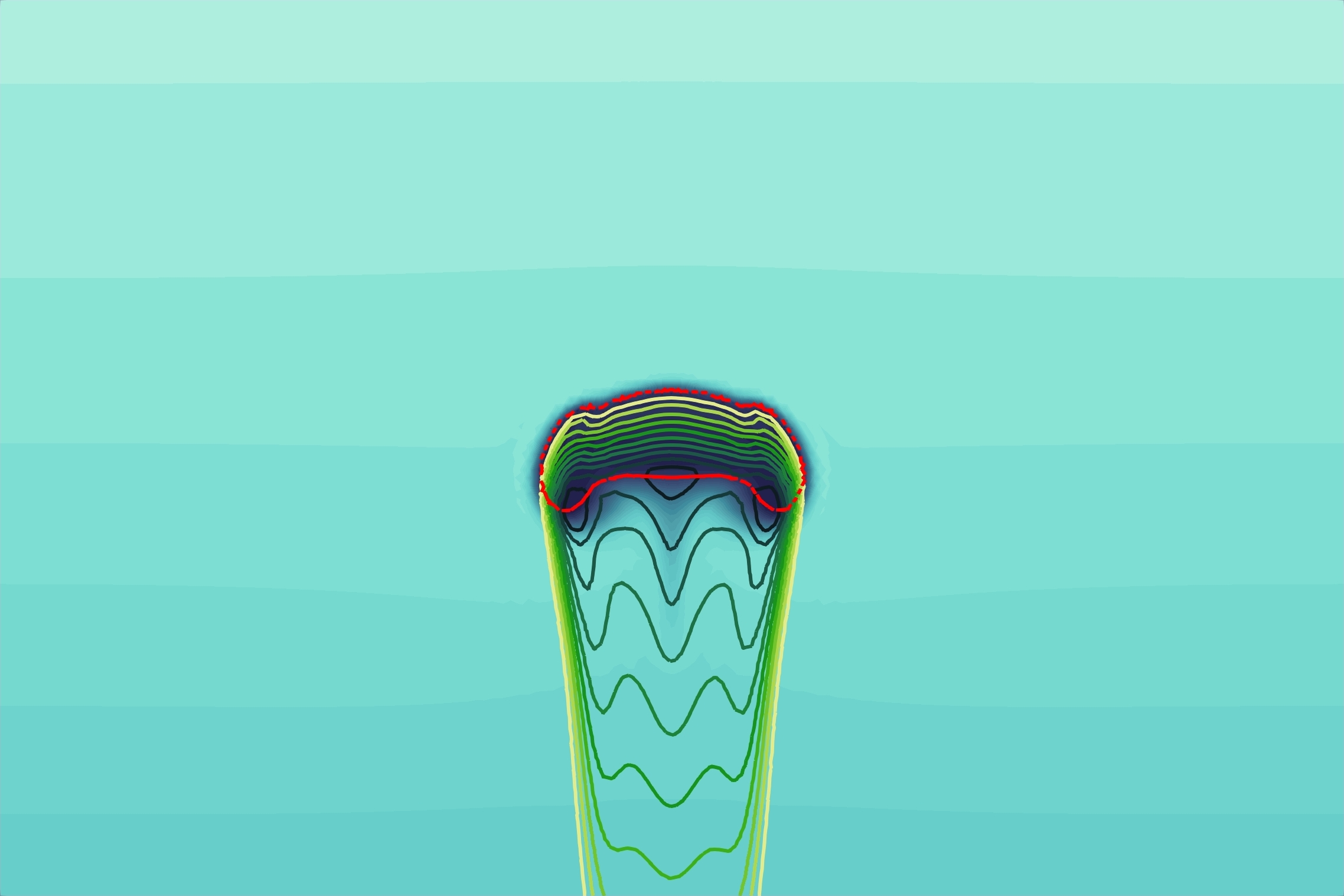}
    \includegraphics[height=3.3cm]{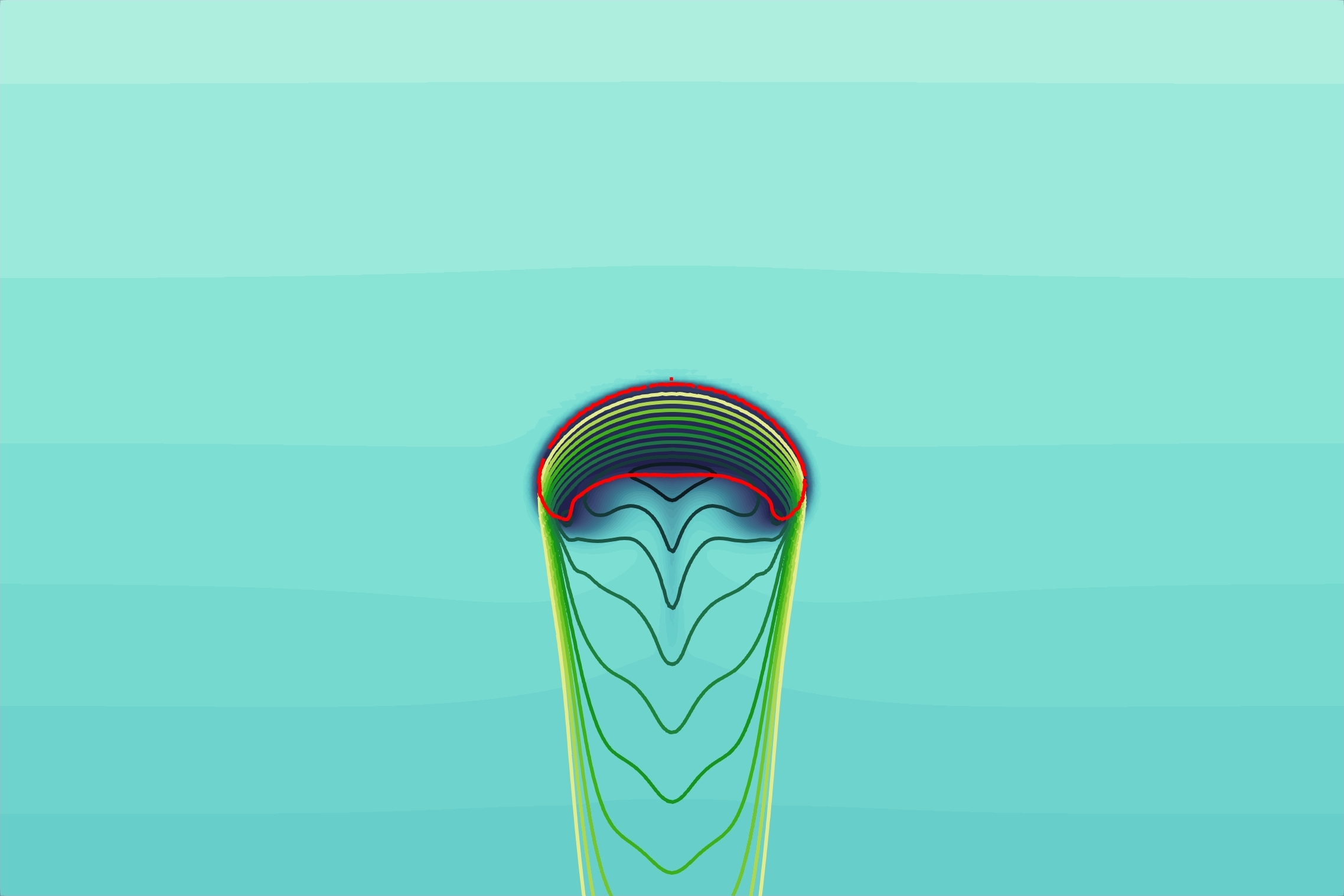}
    \includegraphics[height=3.3cm]{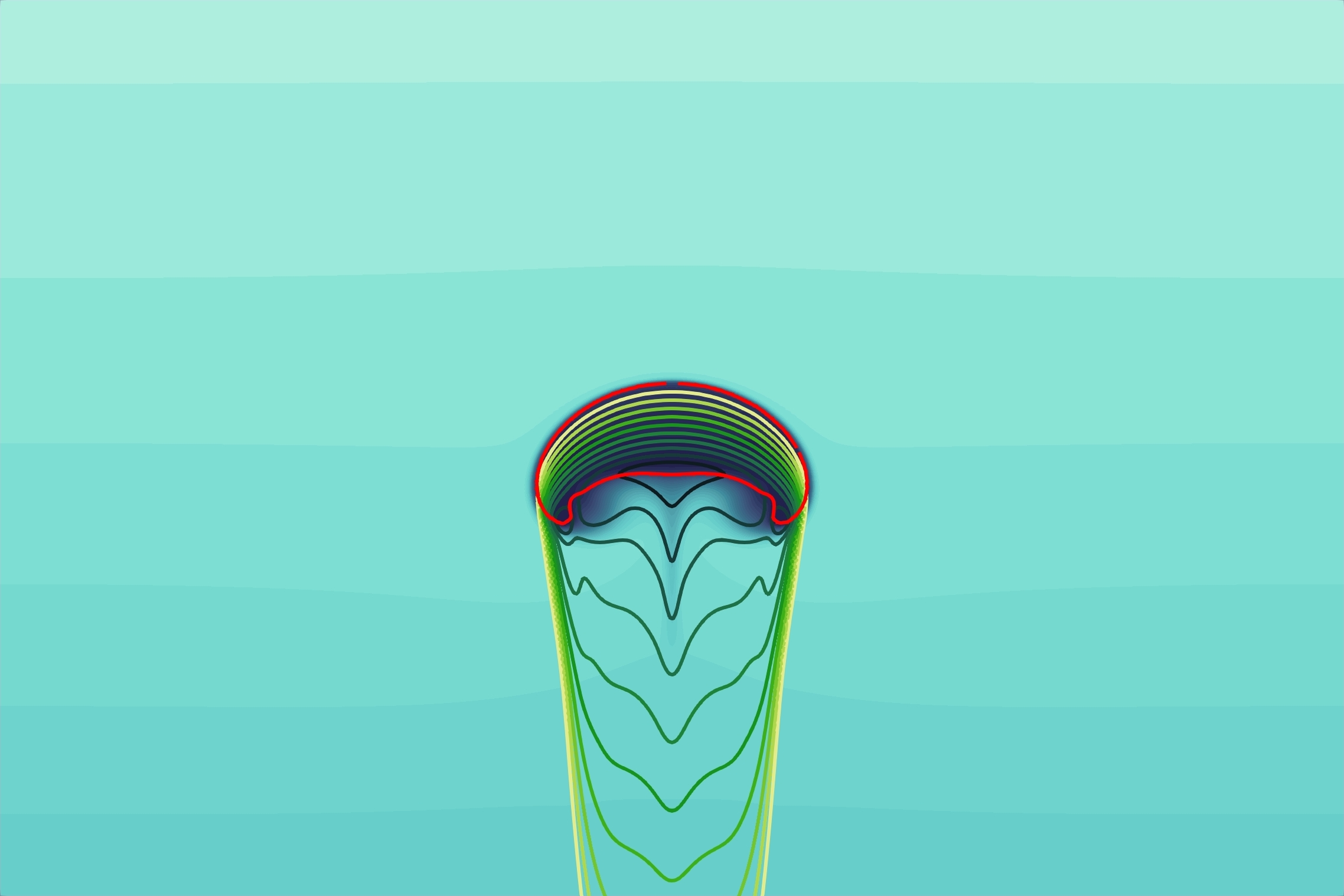}\\[4pt]
    \includegraphics[height=3.3cm]{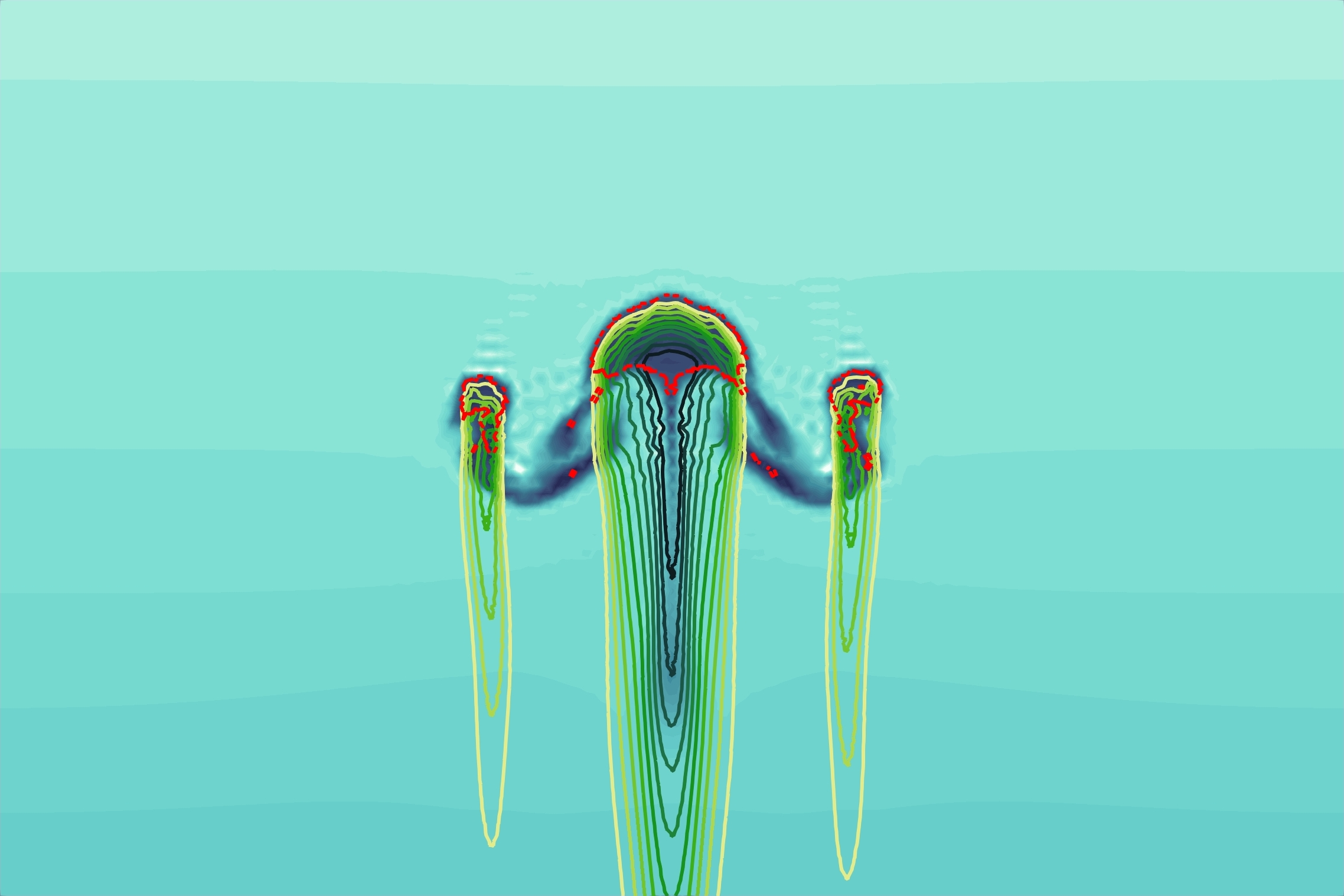}
    \includegraphics[height=3.3cm]{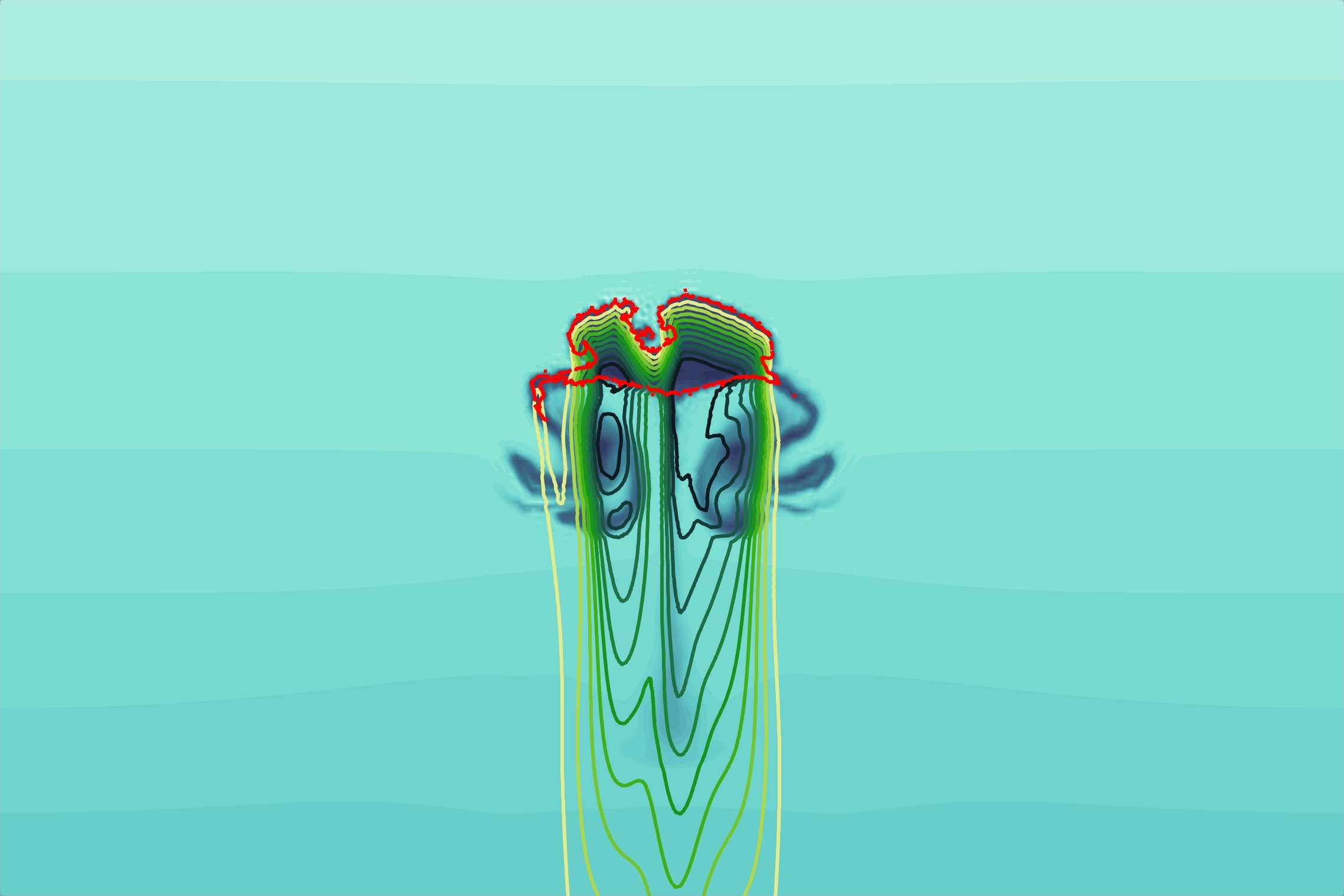}
    \includegraphics[height=3.3cm]{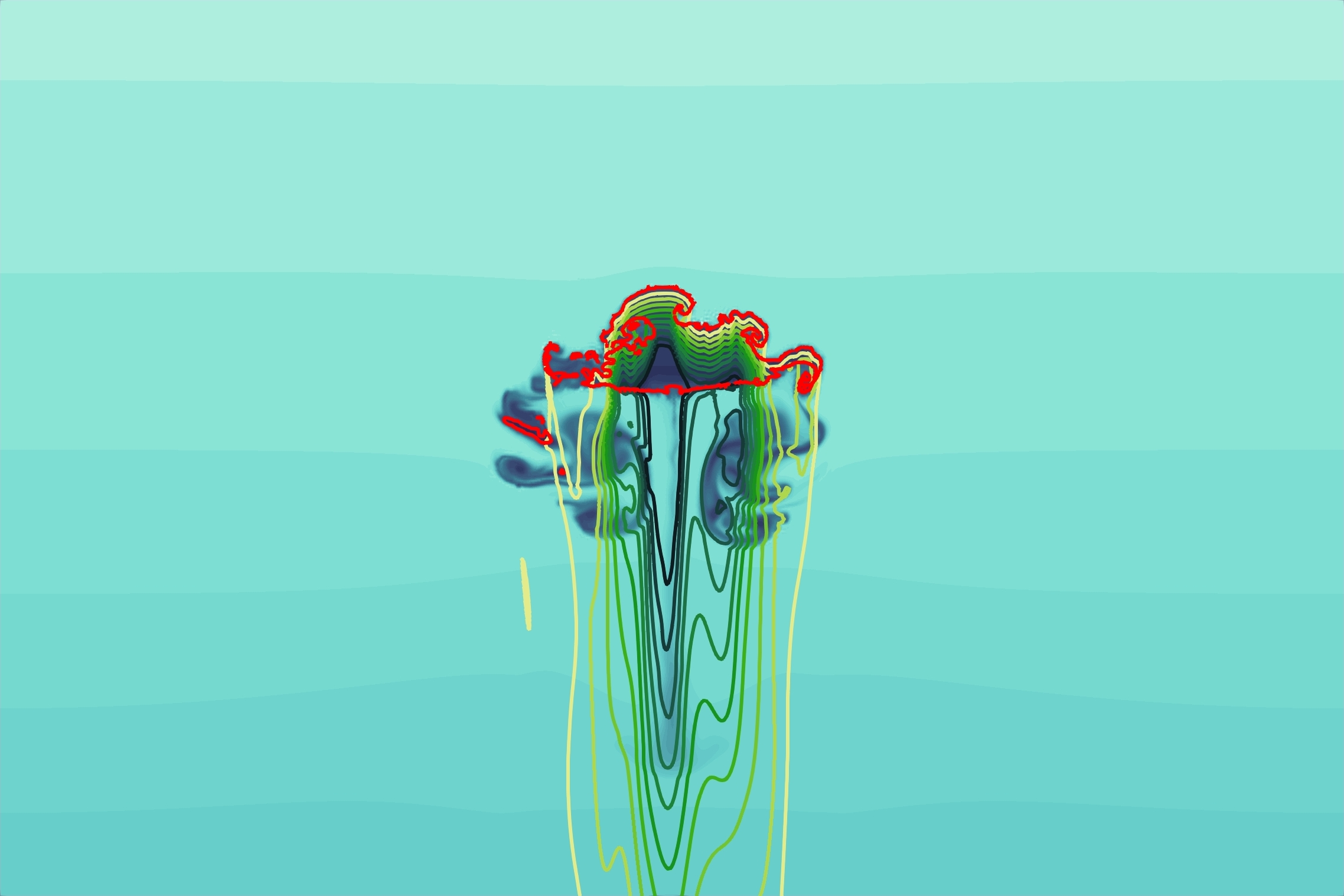}
  \end{minipage}
  \begin{minipage}[t]{0.96cm}
    \vspace*{-40pt}
    \includegraphics[height=3.3cm]{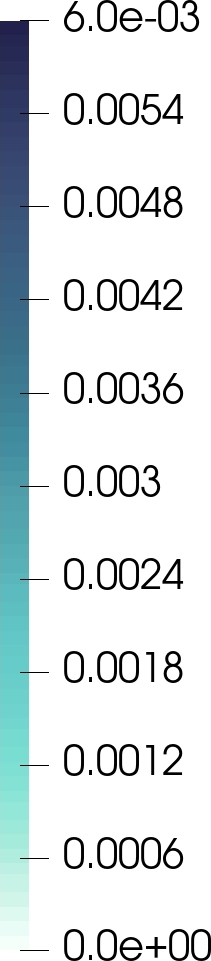}
  \end{minipage}
  \caption{Example 5: Water vapour density (\,\unit{\kg\per\m^2}) solution (background), cloud density (\,\unit{\kg\per\m^2}) edge contour at $0.001$ (red) and rain water density (\,\unit{\kg\per\m^2}) contours at $\{10^{-6}, 2\times10^{-6}, 3\times10^{-6}, 4\times10^{-6}, 5\times10^{-6}, 6\times10^{-6}, 7\times10^{-6}, 8\times10^{-6}, 9\times10^{-6}, 10^{-5} \}$ (greens), for a two-dimensional rising thermal leading to precipitation. Top row: $t=300\,\unit{s}$, bottom row: $t=600\,\unit{s}$. Discretisation parameters as in \Cref{fig.ex5_results.rising_thermal_rain_2d.velocity}.}
  \label{fig.ex5_results.rising_thermal_rain_2d.vapour_cloud_rain}
\end{figure}

\begin{figure}
  \centering
  \includegraphics{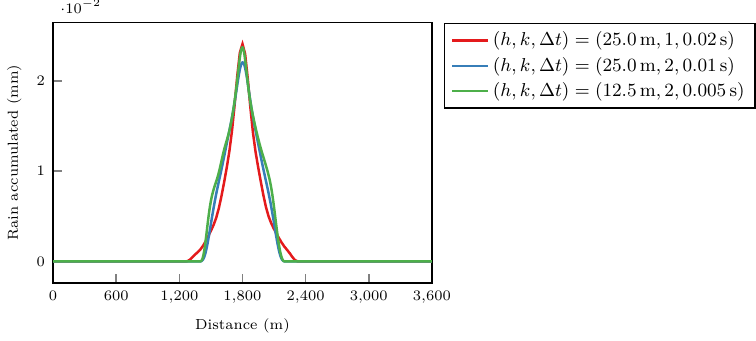}
  \caption{Example 5: Total rain fallout at the bottom of the domain for a two-dimensional rising thermal leading to precipitation after $t=600\,\unit{s}$.}
  \label{fig.results.rising_thermal_rain_2d_rain_fall}
\end{figure}

The resulting velocity field at $t=600\,\unit{s}$ can be seen in
\Cref{fig.ex5_results.rising_thermal_rain_2d.velocity}, and the water
vapour density, cloud boundary contour and rain density contours at $t=300\,\unit{s},600\,\unit{s}$ in
\Cref{fig.ex5_results.rising_thermal_rain_2d.vapour_cloud_rain}. We also
show the total rain transported out of the bottom of the domain in
\Cref{fig.results.rising_thermal_rain_2d_rain_fall}.

Looking at the velocity solution, we see that the results are
consistent, and higher order results in faster velocities and more details
in the velocity field.
Looking at the water vapour, cloud and rain density solutions at 
$t=300\,\unit{s}$ in
\Cref{fig.ex5_results.rising_thermal_rain_2d.vapour_cloud_rain}, we have 
similar results for all three discretisations. In particular,
the results for $k=2$ are consistent with the literature~\cite{GC91}. At
$t=600\,\unit{s}$, the densities have more pronounced differences;
notably, in the case $k=1$, the cloud has split into three separate
clouds, which is mirrored in the rain contours and in the more
spread out total rainfall seen in
\Cref{fig.results.rising_thermal_rain_2d_rain_fall}. We note that the 
single rain column is consistent with results in the literature~\cite{GS96}.
With regard to the water vapour, we see that the $k=2$ results show 
significantly more small-scale features than the $k=1$ solution.
Finally, we note that we do not preserve the symmetry of the initial 
condition due to the use of unstructured meshes. This lack of symmetry can be
observed consistently in 
\Cref{fig.ex5_results.rising_thermal_rain_2d.velocity},
\Cref{fig.ex5_results.rising_thermal_rain_2d.vapour_cloud_rain}
and \Cref{fig.results.rising_thermal_rain_2d_rain_fall}.

\subsection{Example 6: Rising thermal with rain in three spatial dimensions}
\label{sec.numex:subsec.rising_thermal_rain_3d}
We take the extension of the previous example as presented in~\cite{GC93}.

\subsubsection{Set-up}
The domain is $\O=(0, 3.6\,\unit{km})\times(0, 3.6\,\unit{km})\times
(0, 2.4\,\unit{km})$, periodic boundary conditions are set on the
horizontal boundaries and solid wall (slip) boundary conditions at the
top and bottom boundaries. The hydrostatic state and the initial
conditions are defined analogously to those in
\Cref{sec.numex:subsec.rising_thermal_rain_2d}, however, the bubble is
now three-dimensional.

\subsubsection{Results}
We consider an unstructured tetrahedral mesh of the domain with
$h=100\,\unit{m}$. We consider the polynomial orders $k=1,2$, resulting in
approximately $4\times 10^6$ and $10^7$ unknowns for the DG space for
the primary variables. The time-step is chosen according to the
time-step restriction, resulting in $\dt=0.04\,\unit{s}$ and 
$\dt=0.02\,\unit{s}$, respectively. The artificial diffusion parameter is again chosen as
small as possible, resulting in $\gamma=0.6$.

\begin{figure}
  \centering
  \begin{minipage}[t]{1.8cm}
    \centering
    \vspace*{-130pt}
    \includegraphics[height=4cm]{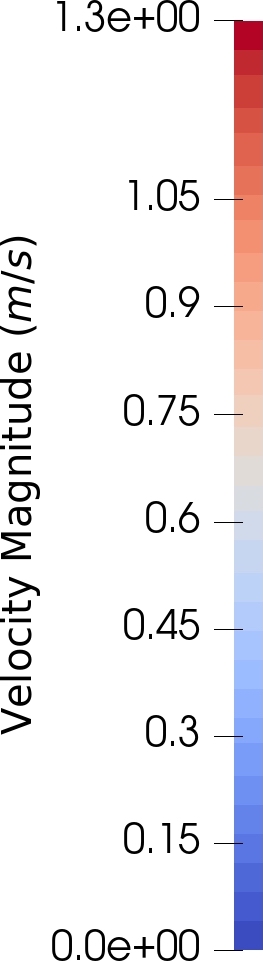}
  \end{minipage}
  \begin{minipage}[t]{10cm}
    \centering
    \includegraphics[height=5cm]{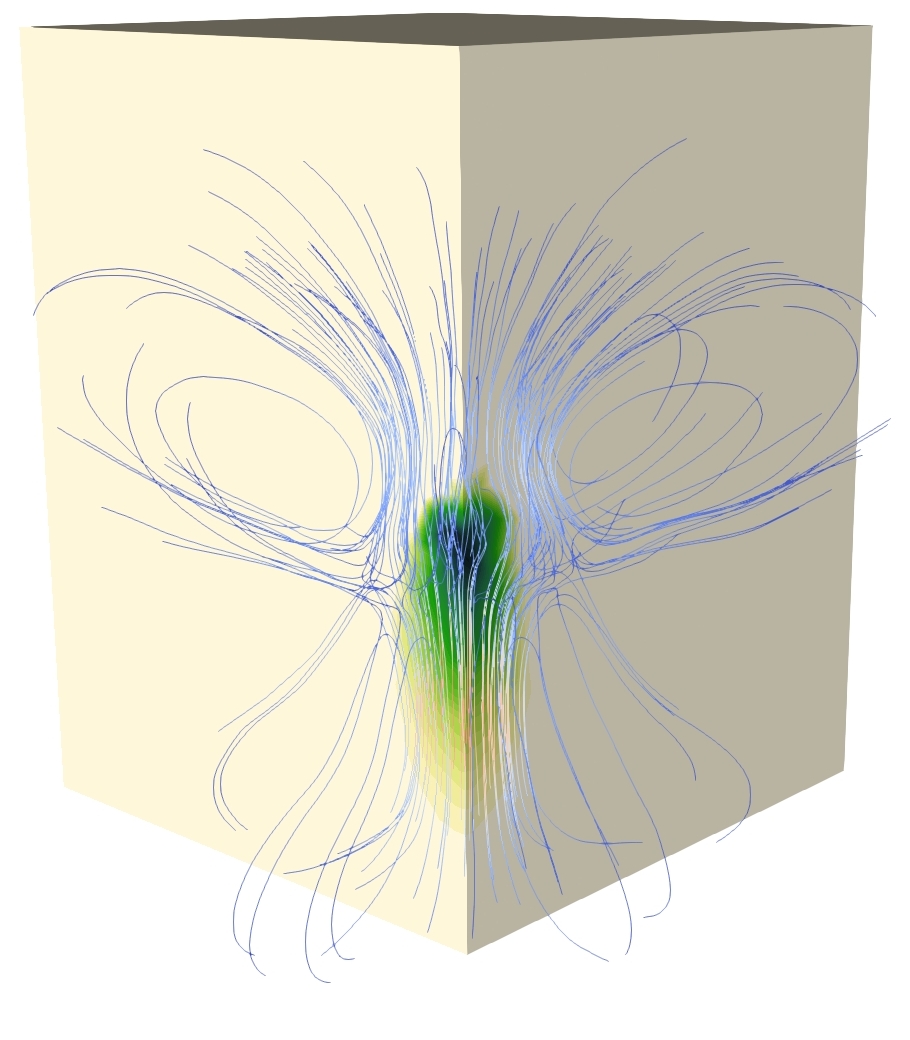}
    \hspace*{0.5cm}
    \includegraphics[height=5cm]{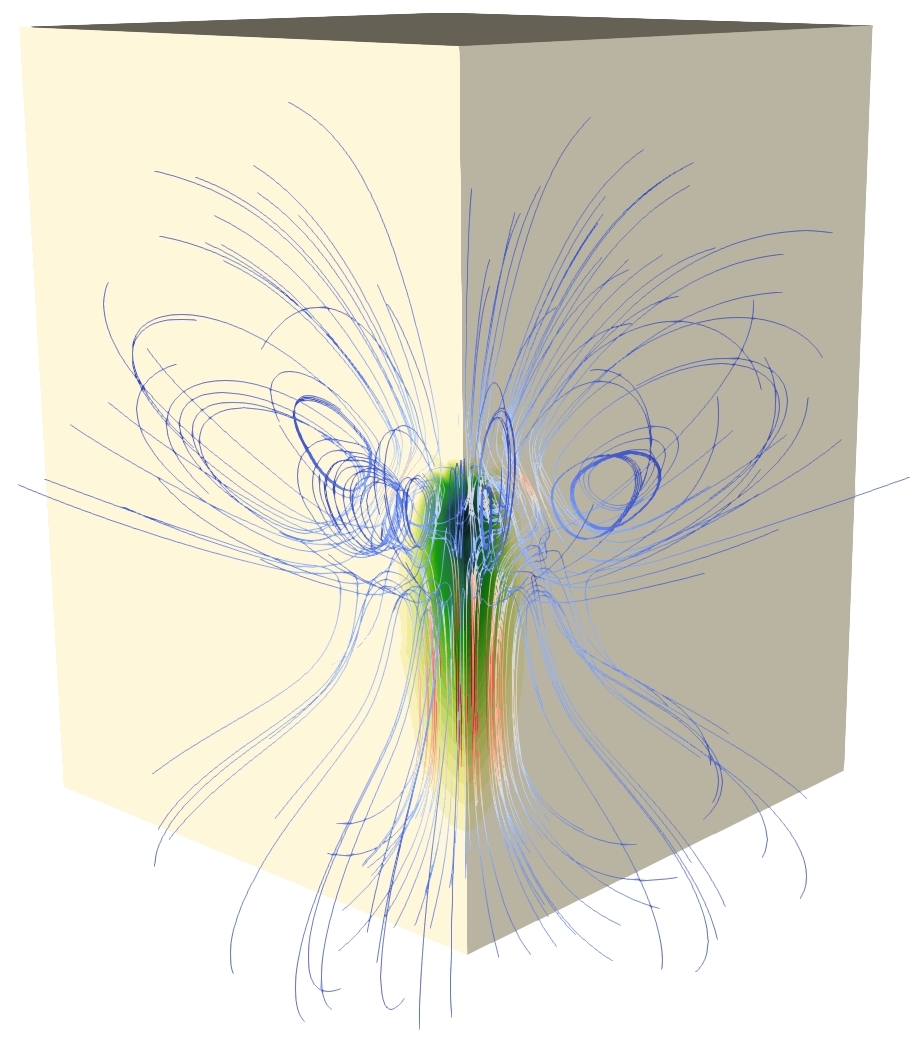}
  \end{minipage}
  \begin{minipage}[t]{1.8cm}
    \centering
    \vspace*{-130pt}
    \includegraphics[height=4cm]{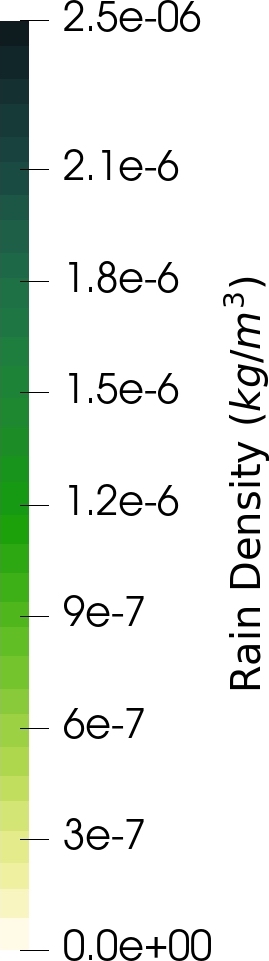}
  \end{minipage}
  \caption{Example 6: Velocity and rain water density results for a three-dimensional rising thermal leading to precipitation on an unstructured tetrahedral mesh at $t=360\,\unit{s}$. Left: $h=100\,\unit{m}$, $k=1$, $\dt=0.04$, Right: $h=100\,\unit{m}$, $k=2$, $\dt=0.02$}
  \label{fig.ex6_results.rising_thermal_rain_3d_vel}
\end{figure}

\begin{figure}
  \centering
  \includegraphics[height=3.3cm]{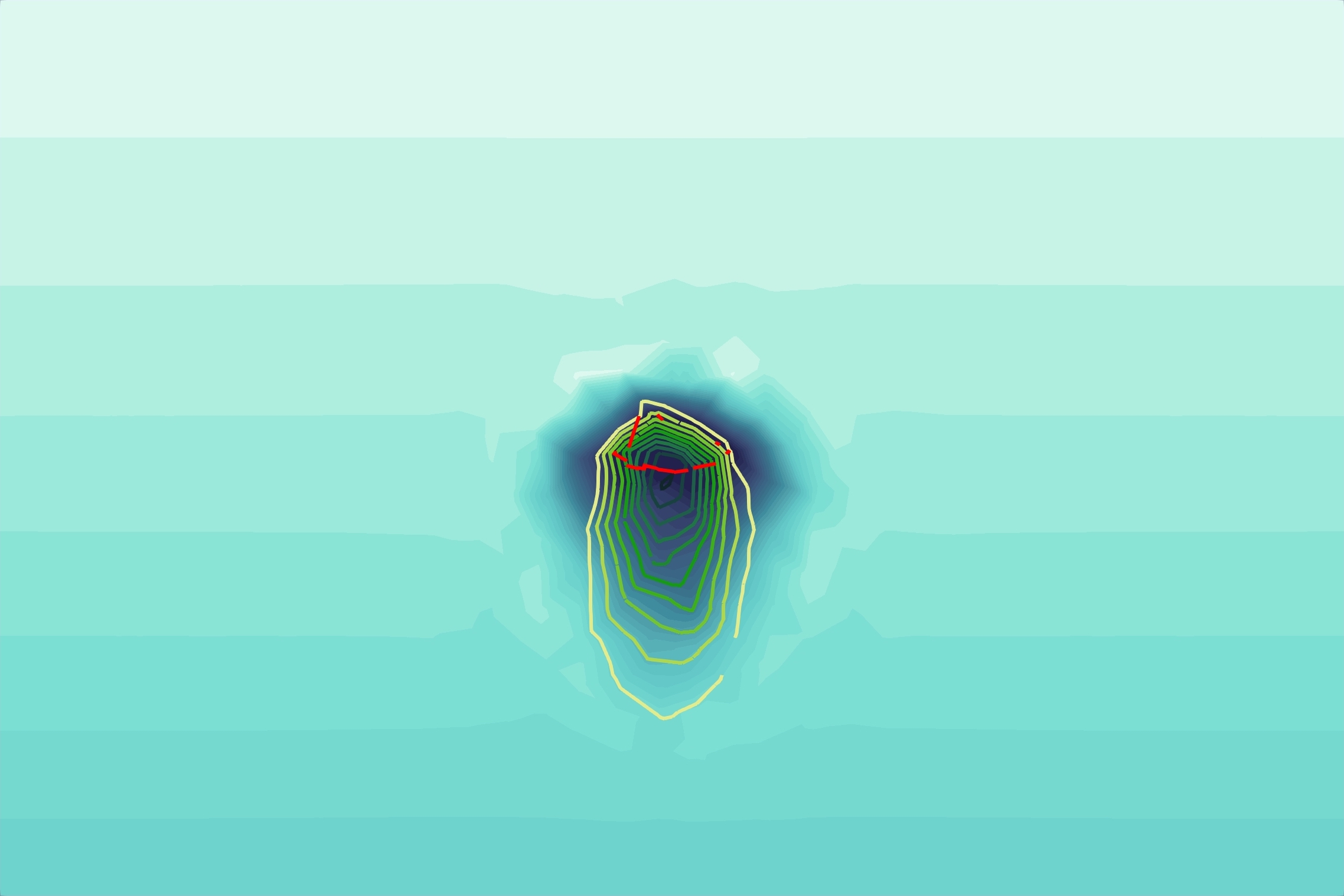}
  \includegraphics[height=3.3cm]{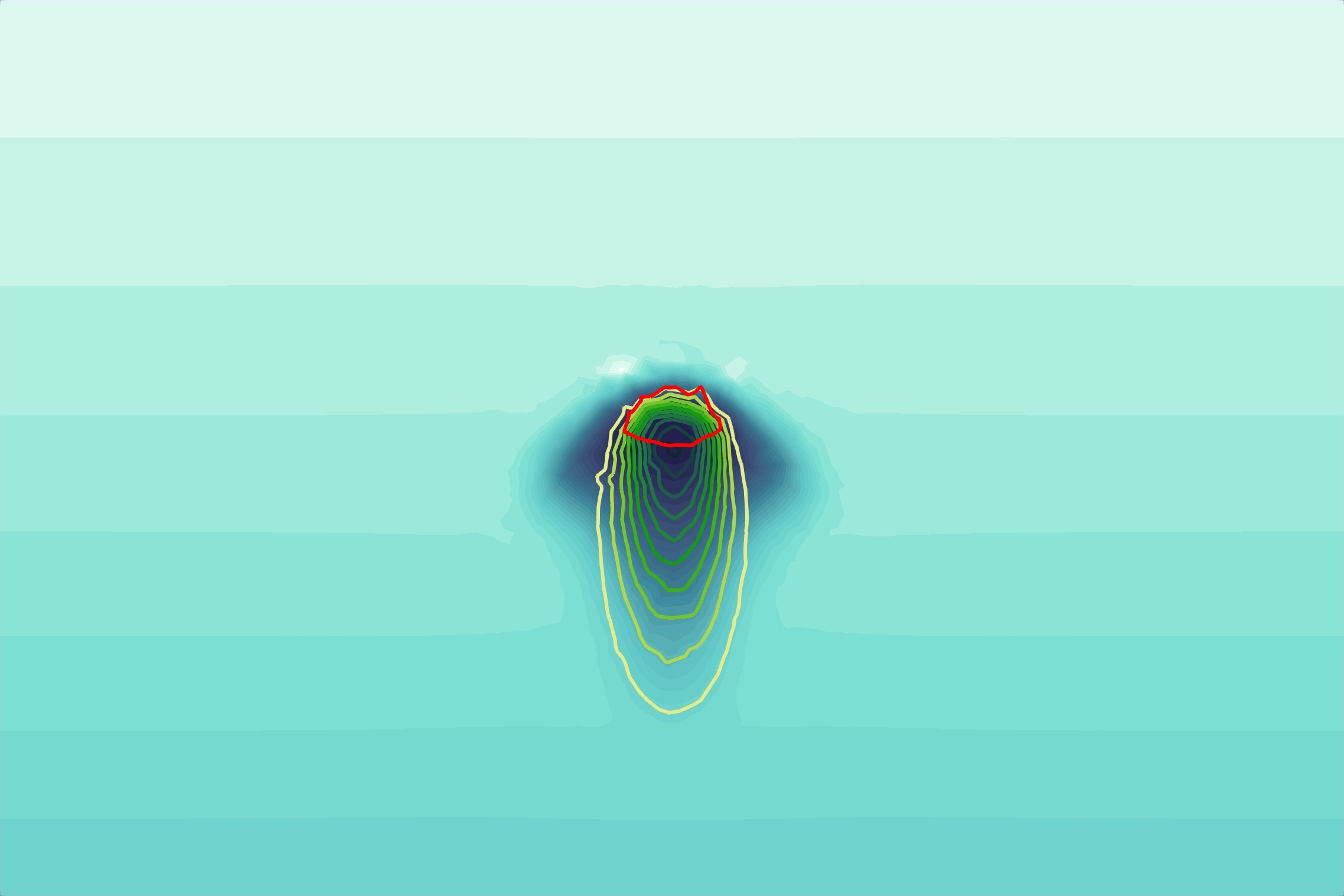}
  \includegraphics[height=3.3cm]{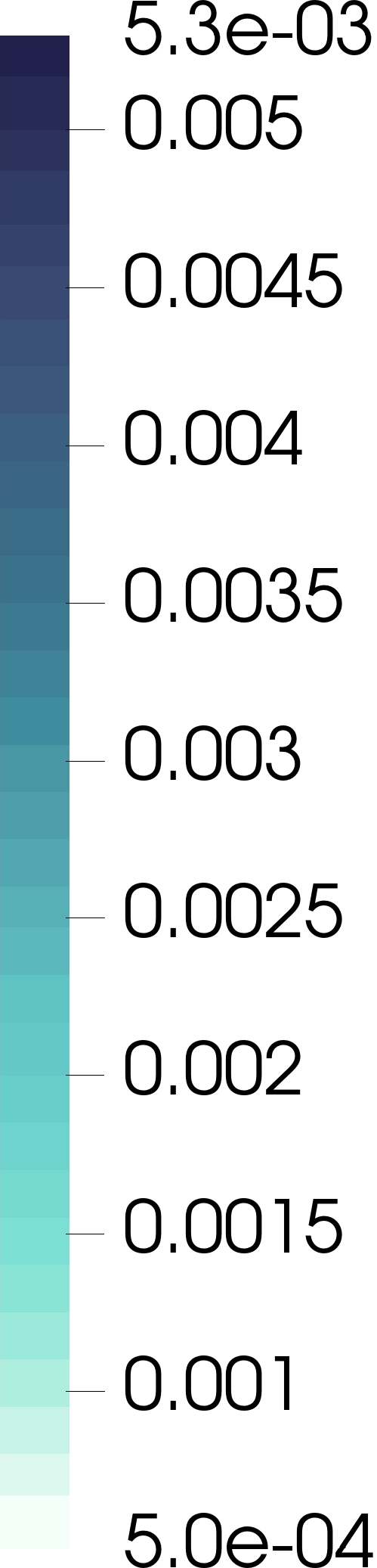}
  \caption{Example 6: Water vapour density (background), cloud density contour at $0.0001$ (red) and rain density contours at $\{2.4\times10^{-7}, 4.8\times10^{-7}, 7.2\times10^{-7}, 9.6\times10^{-7}, 1.2\times10^{-6}, 1.44\times10^{-6}, 1.68\times10^{-6}, 1.92\times10^{-6}, 2.16\times10^{-6}, 2.4\times10^{-6}\}$ (greens) for a three-dimensional rising thermal leading to precipitation on an unstructured tetrahedral mesh in the plane $y=1.8\,\unit{km}$ at $t=360\,\unit{s}$
   }
  \label{fig.ex6_results.rising_thermal_rain_3d_combined}
\end{figure}

The rain density for a quarter of the domain, together with velocity
streamlines in the remaining domain, can be seen in
\Cref{fig.ex6_results.rising_thermal_rain_3d_vel}, and a vertical slice through the centre of the domain ($y=1.8\,\unit{km}$) with the water vapour density, cloud outline contour and rain contours can be seen in \Cref{fig.ex6_results.rising_thermal_rain_3d_combined}. Here, we again see faster
and more dynamics in the velocity solution for the case $k=2$ and that
the rain density is concentrated more towards the centre of the domain,
consistent with our two-dimensional results. The higher velocity is also
consistent with the observation that the rain falls from a higher
height in the case $k=2$. Furthermore, we observe that the higher order simulation preserves the symmetry of significantly better than the case for $k=1$.

\subsection{Example 7: Squall Line}
\label{sec.numex:subsec.squall_line}

As a final example, we consider an idealised test-case presented by \cite{GGD12}. The initial condition for this consists of a synthetic vertical profile, based on a typical environment typical for mid-latitude squall lines~\cite{GGD12, WKR88,RKW88}. While we remain in the original two-dimensional setting from \cite{GGD12}, we take the precise values from \cite{TMQ22}. In particular, the set-up in \cite{TMQ22} considers a smaller domain in the horizontal direction, and a shear wind in the opposite direction.

\subsubsection{Set-up}

The domain is $\Omega=(0, 150\,\unit{km})\times(0, 24\,\unit{\km})$. Periodic boundary conditions are applied on the lateral boundaries and a free-slip on the top and bottom boundary conditions. To avoid the reflection of waves from the non-physical top boundary condition, we apply an explicit sponge layer in the top $6\,\unit{\km}$ of the domain. The time interval under consideration is $[0, 9000\,\unit{\s}]$.

The hydrostatic background state is computed based on a given vertical profile for the potential temperature $\theta$ and the water vapour mass fraction $q_v$. The specific values for this are taken from \cite[Table A1]{TMQ22}. Further details are given in \Cref{appendix:subsec.squall_line}.

The initial condition is then given by a temperature perturbation bubble,
which we apply to the temperature under the assumption that the pressure is unchanged as above. Details are again given in \Cref{appendix:subsec.squall_line}. The densities remain unchanged. However, the velocity is initialized with a horizontal shear flow with $\ub_x = 12\,\unit{\m\per\s}$ at $z=0$ and decreases linearly to zero at $z=2.5\,\unit{\km}$. 

\subsubsection{Results}
\begin{figure}
  \centering
  \includegraphics[width=8cm]{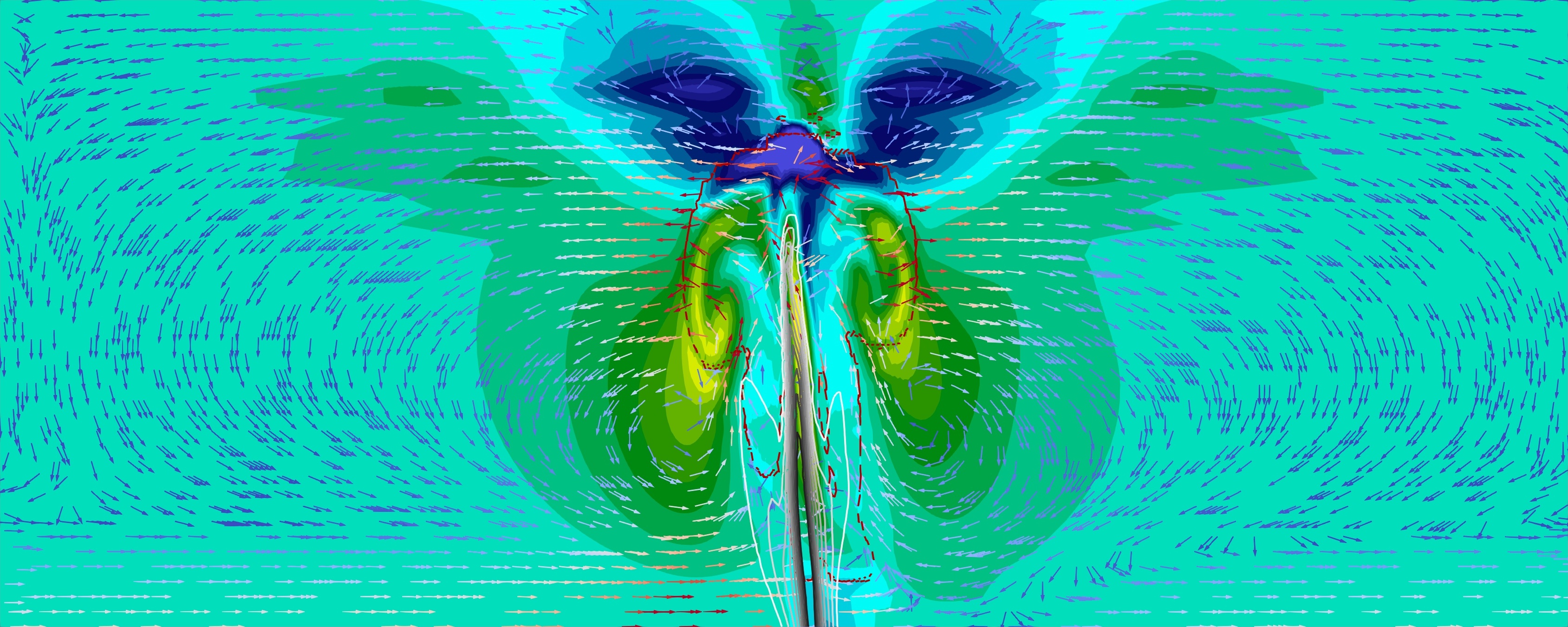}
  \includegraphics[width=8cm]{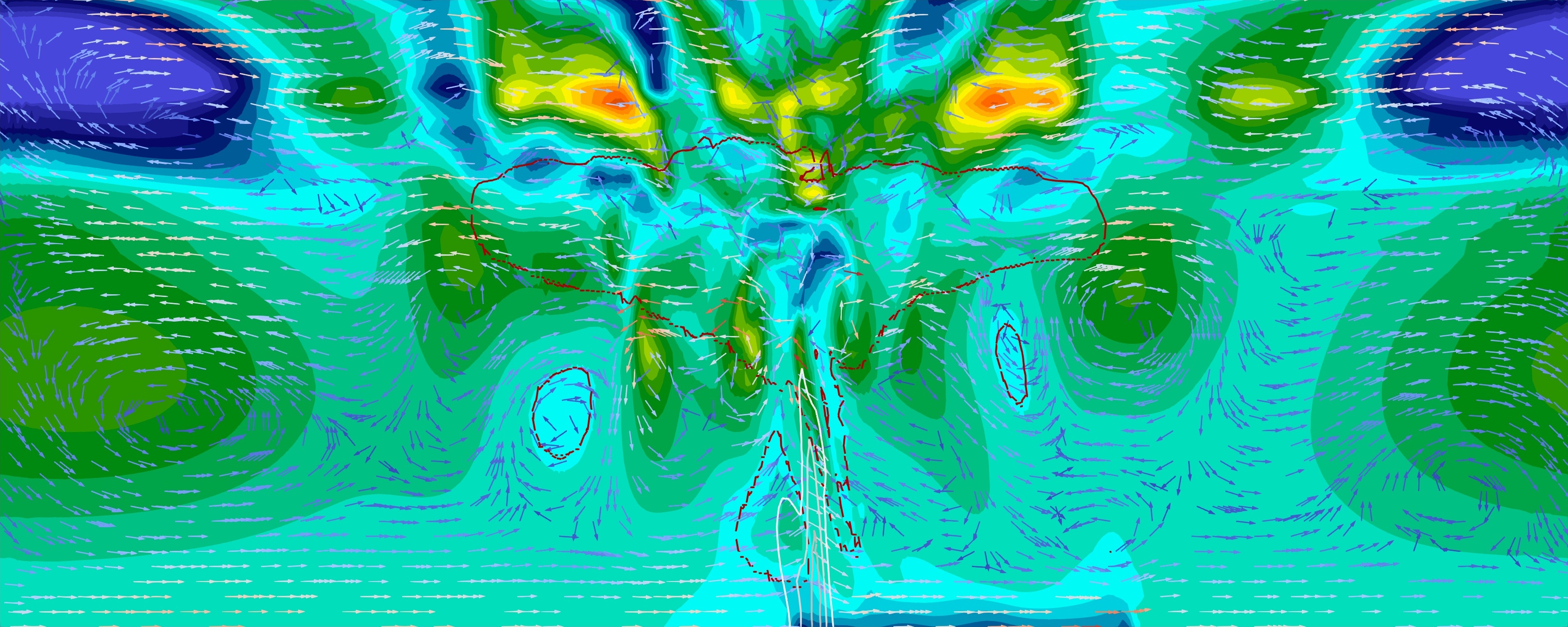}\\[2pt]
  \includegraphics[width=8cm]{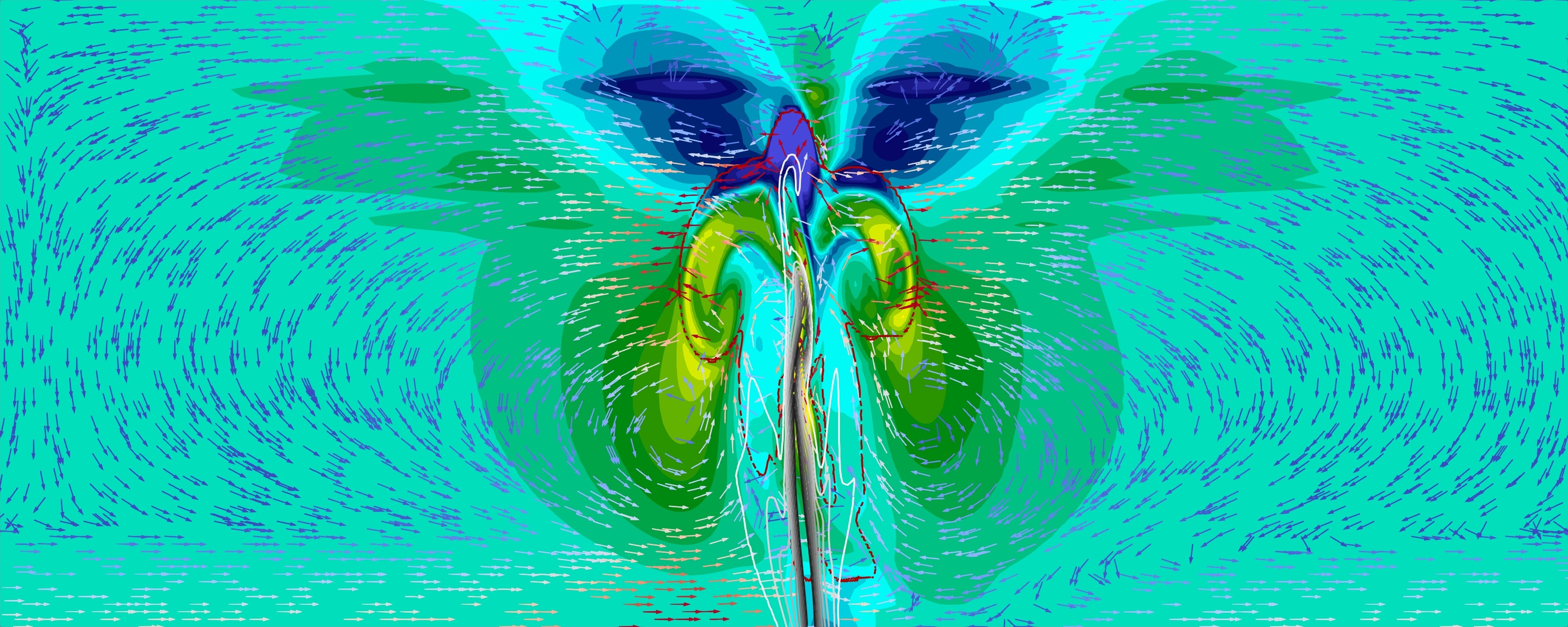}
  \includegraphics[width=8cm]{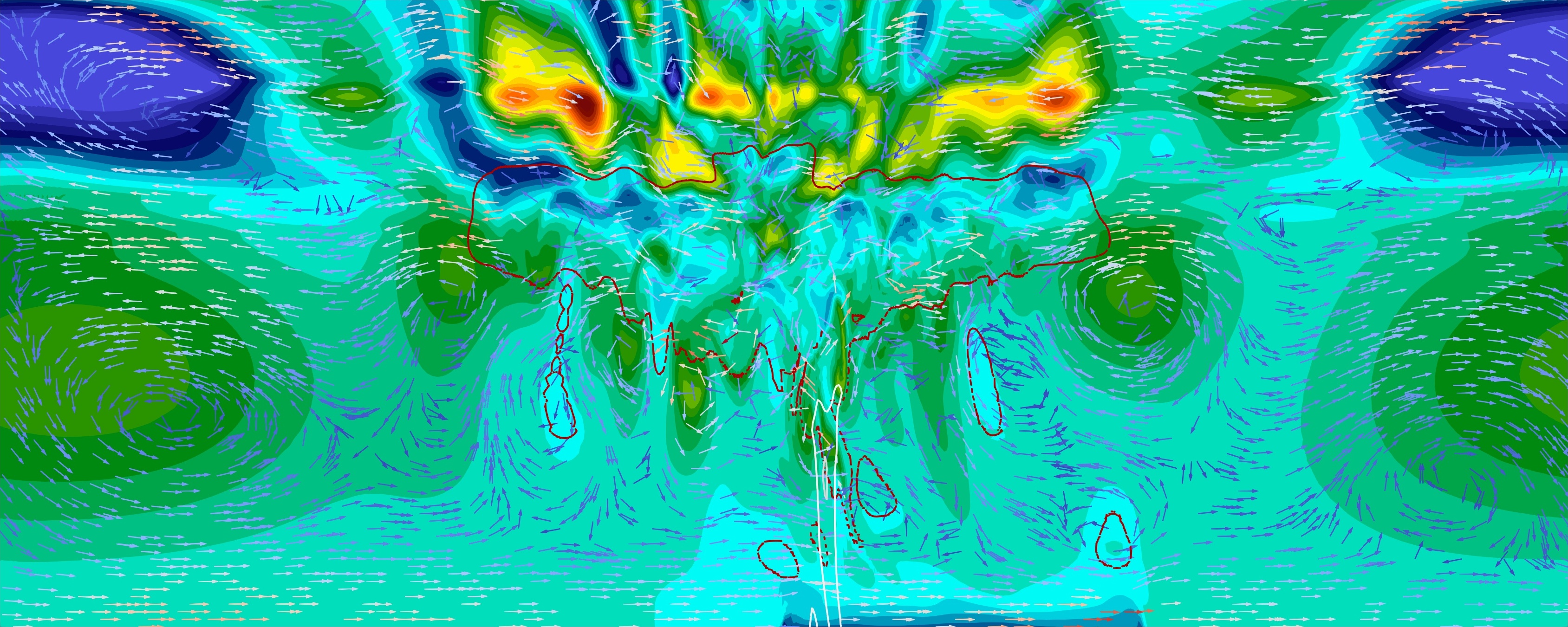}\\[2pt]
  \includegraphics[width=8cm]{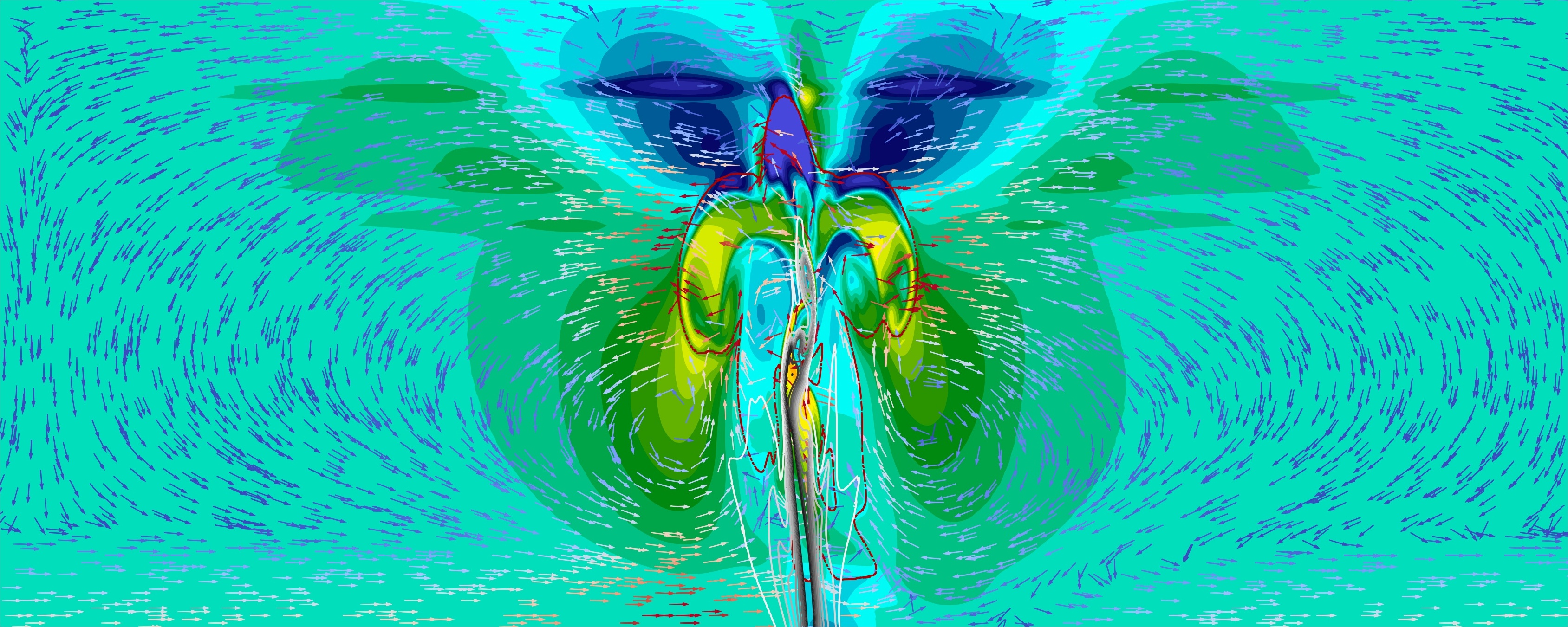}
  \includegraphics[width=8cm]{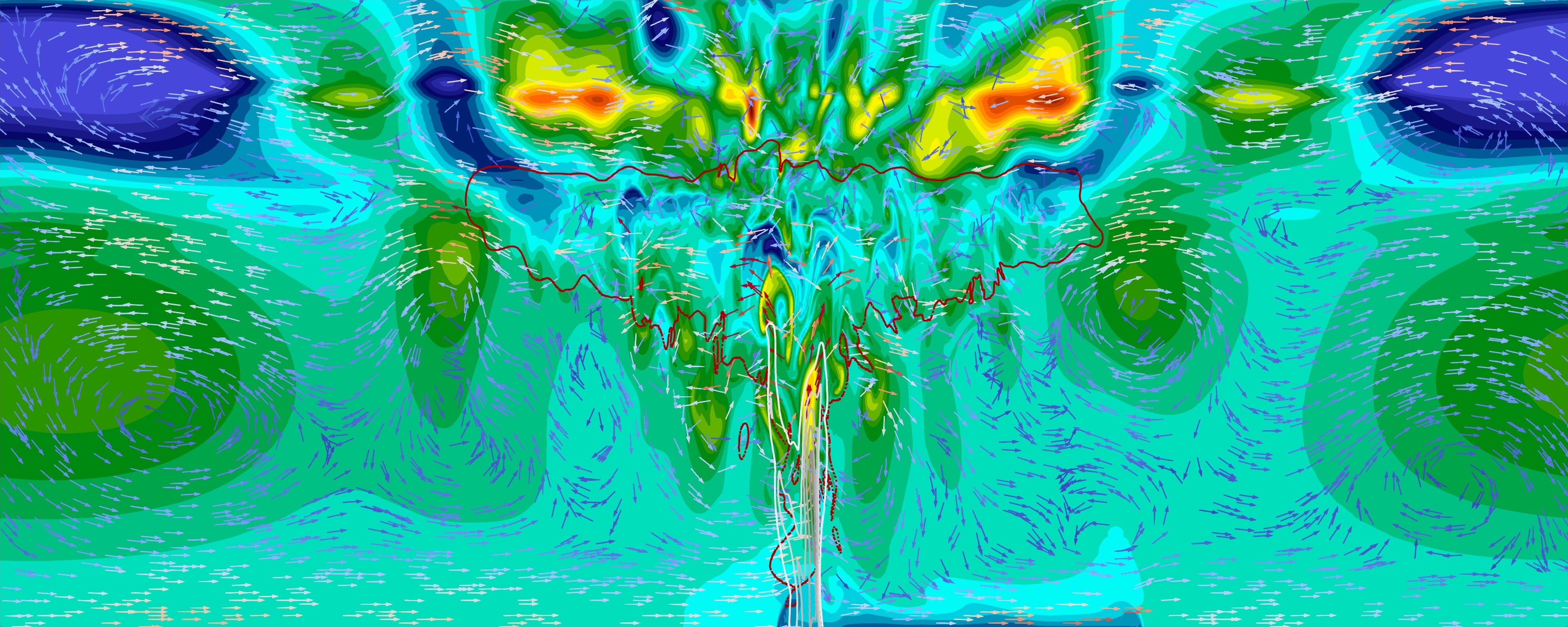}\\[6pt]
  \includegraphics[width=16.1cm]{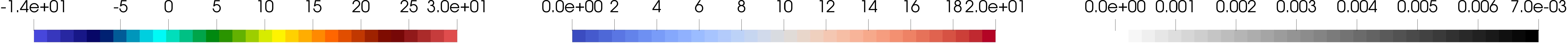}
  \caption{Example 7: Potential temperature perturbation (background), cloud density contour at $10^{-5}$ (red), rain density contours (gray), and velocity field (blue to red arrows) for a two dimensional squall-line test below $z=18\,\unit{km}$. Left: $t=1500\,\unit{s}$, Right: $t=1500\,\unit{s}$. Computed on unstructured simplicial meshes with $k=1$. Top: $h=500\,\unit{m}$, $\dt=0.2\,\unit{s}$, centre: $h=250\,\unit{m}$, $\dt=0.1\,\unit{s}$, bottom: $h=125\,\unit{m}$, $\dt=0.05\,\unit{s}$.}
  \label{fig.ex7:different_meshes}
\end{figure}

\begin{figure}
  \centering
  \includegraphics[width=8cm]{img/ex7_h125k1dt0.05_t1500.jpeg}
  \includegraphics[width=8cm]{img/ex7_h125k1dt0.05_t3000.jpeg}\\[2pt]
  \includegraphics[width=8cm]{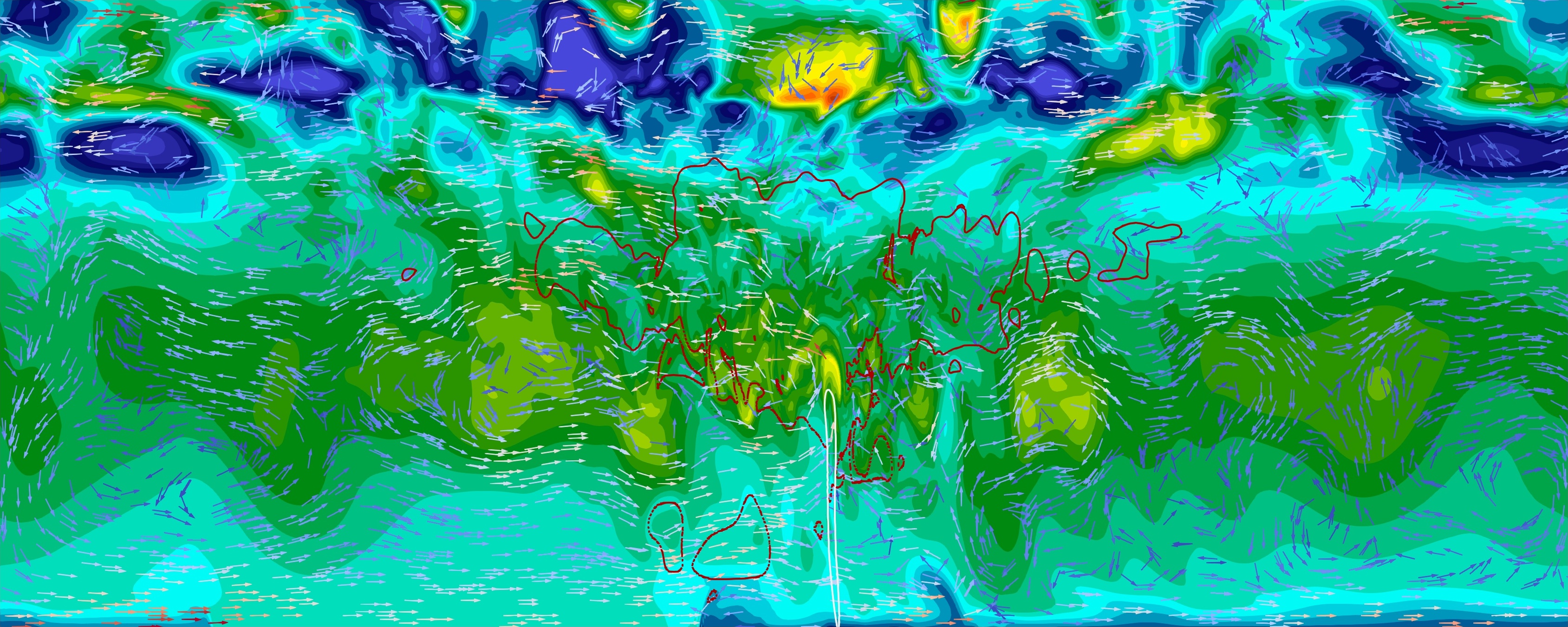}
  \includegraphics[width=8cm]{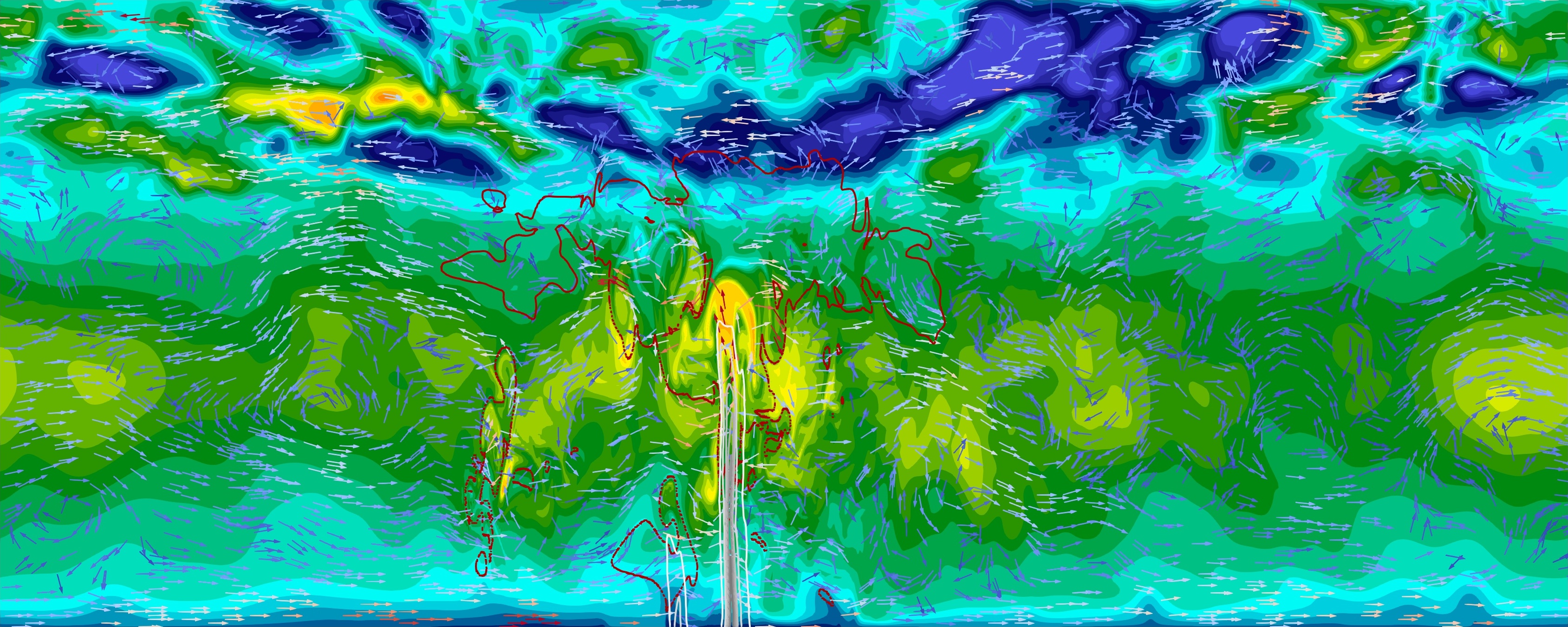}\\[6pt]
  \includegraphics[width=16.1cm]{img/ex7_colorbar}
  \caption{Example 7: Potential temperature perturbation (background), cloud density contour at $10^{-5}$ (red), rain density contours (gray), and velocity field (blue to red arrows) for a two dimensional squall-line test below $z=18\,\unit{km}$. Computed on unstructured simplicial meshes with $h=125\,\unit{m}$, $k=1$, and $\dt=0.05\,\unit{s}$. From left to right and top to bottom: $t=1500, 3000, 6000, 9000\,\unit{s}$.}
  \label{fig.ex7:temporal_evolution}
\end{figure}

\begin{figure}
  \centering
  \includegraphics{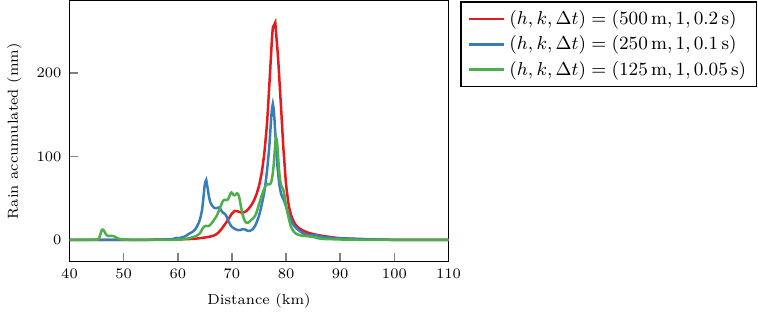}
  \caption{Example 7: Total rain fallout at the bottom of the domain for a two dimensional squall-line test after $t=9000\,\unit{s}$.}
  \label{fig.ex7:rain_fallout}
\end{figure}

We consider unstructured simplicial meshes of the domain with $h=500, 250$, and $125\,\unit{m}$, respectively. On these meshes, we consider $k=1$ and time steps as large as the time-step restriction allows. We present the resulting perturbation of the potential temperature and velocity field, together with the cloud outline and rain contours at $t=1500\,\unit{s}$ and $t=3000\,\unit{s}$ in \Cref{fig.ex7:different_meshes}. We further show the storm's evolution in the same quantities resulting from our finest mesh at  $t=1500, 3000, 6000,$ and $9000\,\unit{s}$ in \Cref{fig.ex7:temporal_evolution}. Finally, we plot the total rain fallout at the bottom of the domain in \Cref{fig.ex7:rain_fallout}.

The main cloud forms at around $t=300\,\unit{s}$, and the precipitating water reaches the ground at around $t=900\,\unit{s}$ (not shown). We note that both these times are significantly earlier than those reported in \cite{MMV13}. However, the latter uses a different equation for the temperature/energy, slightly different source terms governing phase changes and a different formula for the saturation vapour pressure. Nevertheless, this occurs in the convective tower as expected, and a pool of colder and denser air forms downwind of the tower, see \Cref{fig.ex7:different_meshes}, consistent with the literature~\cite{GGD12,MMV13,TMQ22}. The main cloud then spreads out into the expected anvil shape, and the profile of the potential temperature perturbation is consistent with the literature. We attribute the visual differences to the fact that our energy equation takes the effects due to moisture into account, and there are differences in the source terms governing phase changes.

The total accumulated rain is larger than that reported in the literature. However, it reduces with mesh refinement, consistent with faster storm development and less total accumulated rain reported in \cite{GGD12,TMQ22}. The accumulated rain on the lowest mesh is indeed consistent with the results in the literature. We also observe multiple peaks in the continuous Galerkin discretisation in \cite{TMQ22}. However, in that work, the secondary peaks were observed in the downwind direction rather than the upwind direction. However, more accumulated rain in the upwind direction is consistent with the literature, e.g.,~\cite{MMV13}, where the two-dimensional case was also considered. This contrasts \cite{TMQ22}, where a fully three-dimensional simulation was run with a single element with polynomials of order four in the $y$-direction.

\section{Conclusions}
\label{sec.conclusion}

We presented a discontinuous Galerkin method for the equations governing moist atmospheric flows in their conservative form. In particular, we retained thermodynamic details in these equations, such as the different specific heats and heat capacities, which are often neglected in the literature. The source terms governing phase changes were based on Kessler's microphysics closure as taken from the COSMO model. To avoid the difficult modelling of a source term governing condensation, we considered a single \emph{moist} density for both water vapour and cloud water in our hyperbolic balance-law equation set. To recover the individual vapour and cloud densities from the moist density, we had to solve an additional non-linear problem. This reconstruction was necessary in order to evaluate the (numerical) flux and in particular the source terms governing water phase changes. Since the problem of reconstructing the densities is algebraic, we presented a novel approach by solving this non-linear problem in each quadrature point, which is highly parallelisable. For time-stepping, we used an explicit scheme such that only mass-matrix problems had to be solved to advance the solution in time. As the Legendre polynomial basis is orthogonal, the mass matrix is diagonal and the mass-matrix problems could be solved matrix-free. The combination of the local density and temperature reconstruction with the explicit time-stepping scheme led to a highly parallelisable scheme. To stabilise the method, we added an artificial viscosity term with local, dynamic, and asymptotically vanishing viscosity parameter. This is in contrast with the artificial viscosity widely used in the literature on DG for atmospheric flows which uses a constant and non-vanishing scaling, which in turn changes the problem at hand.

We validated our scheme on a number of numerical examples taken from the literature. We illustrated the optimal high-order convergence of the method for polynomial orders from one to four in the case of moist dynamics without rain. We also illustrated that the method works on both structured tensor-product meshes and unstructured simplicial meshes. This flexibility with respect to the mesh and the higher-order convergence are particular advantages of discontinuous Galerkin approaches. The mesh flexibility also allows for good parallel load balancing in every direction. We also noted that even on unstructured meshes, higher-order schemes performed well in retaining structures (symmetry) of the initial data and geometry. Overall, we obtained results in good agreement with the literature.

Our examples, which included rain, were also stable, again with a higher-order choice for the polynomial order, resulting in finer details. Our results for the storm resulting from a squall-line set-up were consistent with the available literature. However, it remains an open and interesting question to study the exact effects of retaining the thermodynamic details in the energy equation, which are usually neglected.

\section*{Data availability statement}
The code used to realise the results presented in this paper is
available online in the github repository
\url{https://github.com/hvonwah/cloud-models-code}, and is archived on
zenodo \cite{HLSvW23_zenodo}.

\section*{Acknowledgments}
SH, JS and HvW acknowledge funding by the Austrian Science Fund
(FWF) through the research program ``Taming complexity in partial
differential systems'' (F65). The computational results presented have
been achieved in part using the Vienna Scientific Cluster (VSC). We
would also like to thank Rupert Klein for many fruitful discussions.

\printbibliography

\appendix
\section{Thermodynamic parameters}\label{appendix:constants}

The values of the thermodynamic variables used are given in \Cref{tab.parameters}.

\begin{table}
  \centering
  \caption{Thermodynamic equation of state parameters for moist air at
    reference temperature. Source:~\citeauthor{CBvdH11}~\citeyear{CBvdH11}.}
  \label{tab.parameters}
  \begin{tabular}{llll}
    \toprule
    Parameter & {Value} & Units & Description\\
    \midrule
    $c_l$      & \num{4218}    & \unit{\joule\per\kg\per\kelvin} & Specific heat of liquid water (at $\Tref$).\\
    $c_{pd}$   & \num{1005}    & \unit{\joule\per\kg\per\kelvin} & Specific heat of dry air at constant pressure.\\
    $c_{pv}$   & \num{1850}    & \unit{\joule\per\kg\per\kelvin} & Specific heat of water vapour at constant pressure.\\
    $c_{vd}$   & \num{718}     & \unit{\joule\per\kg\per\kelvin} & Specific heat of dry air at constant volume.\\
    $c_{vv}$   & \num{1390}    & \unit{\joule\per\kg\per\kelvin} & Specific heat of water vapour at constant volume.\\
    $\eref$    & \num{610.7}   & \unit{\pascal} & Saturation vapour pressure with respect to water (at $\Tref$).\\
    $\Lref$    & \num{2.835e6} & \unit{\joule\per\kg} & Latent heat of vaporization (at $\Tref$).\\
    $R_d$      & \num{287.05}  & \unit{\joule\per\kg\per\kelvin} & Gas constant for dry air.\\
    $R_v$      & \num{461.51}  & \unit{\joule\per\kg\per\kelvin} & Gas constant for water vapour.\\
    $\Tref$    & \num{273.15}  & \unit{\kelvin} & Reference temperature.\\
    $\pref$    & \num{1.0e5}   & \unit{\pascal} & Reference pressure.\\
    $\epsilon$ & \num{0.622}   & & $\epsilon\equiv R_d / R_v$.\\
    \bottomrule
  \end{tabular}
\end{table}

\section{Computation of initial conditions}\label{appendix.setups}

\subsection{Inertia gravity waves in a saturated atmosphere}
\label{appendix.gravitywaves}

The hydrostatic state is determined by solving the (1-dimensional)
problem: Find $\pbar, \rhobar_{d}, \rhobar_{vs}$ and $\Tbar$ such that
\begin{align*}
  \partial_z \pbar &= - (1 + q_w) \rhobar_{d} g,\\
  \pbar &= (\rhobar_{d}R_d + \rhobar_{vs} R_v)\Tbar,\\
    \rhobar_{vs} R_v \Tbar &=
  \eref \left(\frac{\Tbar}{\Tref}\right)^\frac{c_{pv} - c_l}{R_v}
    \exp \left(\frac{\Lref - (c_{pv} - c_l)\Tref }{R_v}\left(\frac{1}{\Tref}
      - \frac{1}{\Tbar}\right) \right),\\
  \theta_e &= \Tbar \left(\frac{\rhobar_{d} R_d \Tbar}{\pref}\right)^{-R_d/(c_{pd} + c_l q_w)}
    \exp\left(\frac{(\Lref + (c_{pv} - c_l)(\Tbar - \Tref))\rhobar_{vs}}{\rhobar_d(c_{pd} + c_{l}q_w)\Tbar}\right).
\end{align*}
and the pressure boundary condition is $\pbar(0)=\pref$. The definition of
$\theta_e$ is taken from~\cite{Ema94}.
The initial perturbation is computed by solving the non-linear 
problem: Find $(\rho_d, \rho_{vs}, T)$ such that
\begin{subequations}\label{eqn.appendix.gravitywaves_pertubation}
  \begin{align}
  \overline{\theta}_e + \theta_e' &= T \left(\frac{\rho_{d} R_d T}{\pref}\right)^{-R_d/(c_{pd} + c_l q_w)}
    \exp\left(\frac{(\Lref + (c_{pv} - c_l)(T - \Tref))\rho_{vs}}{\rho_d(c_{pd} + c_{l}q_w)T}\right),\\
    \pbar &=  (\rho_d R_d + \rho_{vs} R_v) T,\\
     \rho_{vs} R_v T &=
      \eref \left(\frac{T}{\Tref}\right)^\frac{c_{pv} - c_l}{R_v}
        \exp \left(\frac{\Lref - (c_{pv} - c_l)\Tref }{R_v}\left(\frac{1}{\Tref}
          - \frac{1}{T}\right) \right),
  \end{align}
\end{subequations}
and set $\rho_d'=\rho_d - \rhobar_d$, $\rho_v'=\rho_{vs} -
\rhobar_{vs}$, $\rho_m' = \rho_d' - \rho_{vs}'$ and $T'=T-\Tbar$. The
total initial energy (perturbation) can the be computed from
\eqref{eqn.mass-total-energy}.

\subsection{Bryan-Fritsch moist benchmark}
\label{appendix.bryanfritsch}

The hydrostatic base state is computed as in
\Cref{appendix.gravitywaves}.
The initial perturbations are then defined by perturbing the
\emph{density potential temperature}
\begin{equation}\label{eqn.appendix.theta_rho}
  \theta_\rho = T \left(\frac{\pref}{p} \right)^{R_d / c_{pd}}(1 + q_v / \epsilon) = \theta_d(1 + q_v / \epsilon),
\end{equation}
where $\theta_d(p, T)$ is the \emph{dry potential temperature}. This
perturbation is defined as
\begin{equation*}
  \theta' = 2 \cos^2\left(\frac{\pi L}{2}\right),
  \qquad\text{with}\qquad
  L = \min\left\{\left\Vert\left(\frac{x - x_c}{x_r},
    \frac{z - z_c}{z_r}\right)\right\Vert_2, 1 \right\},
\end{equation*}
where $x_c = 10\,\unit{\km}$, $z_c=2\,\unit{\km}$ and $x_r, z_r = 2\,\unit{\km}$. The
perturbation to the potential temperature is then introduced by assuming that
the resulting buoyancy is the same as the buoyancy in a similar the dry test
case with $\theta_d=300\,\unit{\kelvin}$. This leads to the equation
\begin{equation*}
  \theta_d(p, T)\left[1 + \frac{q_{vs}(p, T)}{\epsilon}\right]
    = \theta_{\rho0}(1 + q_w)\left(\frac{\theta_d'}{300} + 1\right),
\end{equation*}
which we solve for $T$ in every point, using the hydrostatic
pressure $p=\pbar$ and \eqref{eqn.saturation-ratio} for the saturation
vapour mixing ratio. With this temperature, we can then reconstruct the
vapour and cloud densities using the saturation vapour pressure.

\subsection{Inertia gravity waves in a saturated atmosphere without initial clouds}
\label{appendix:subsec.gravitywavesnoclouds}

The hydrostatic state is determined by solving the (1-dimensional)
problem: Find $\pbar, \rhobar_{d}, \rhobar_{vs}$ and $\Tbar$ such that
\begin{align*}
  \partial_z \pbar &= - (\rhobar_{d} + \rhobar_{vs}) g,\\
  \pbar &= (\rhobar_{d}R_d + R_v \rhobar_{vs})\Tbar,\\
  \rhobar_{vs} R_v \Tbar &=
  \eref \left(\frac{\Tbar}{\Tref}\right)^\frac{c_{pv} - c_l}{R_v}
    \exp \left(\frac{\Lref - (c_{pv} - c_l)\Tref }{R_v}\left(\frac{1}{\Tref}
      - \frac{1}{\Tbar}\right) \right),\\
  \theta_e &= \Tbar \left(\frac{\rhobar_{d} R_d \Tbar}{\pref}\right)^{-R_d/(c_{pd} + c_l (\rhobar_{vs}/\rhobar_d))}
    \exp\left(\frac{(\Lref + (c_{pv} - c_l)(\Tbar - \Tref))\rhobar_{vs}}{(c_{pd}\rhobar_d + c_{l}\rhobar_{vs})\Tbar}\right),
\end{align*}
and the pressure boundary condition $\pbar(0) = \pref$.
The perturbation of the hydrostatic state can then be computed
analogously to \eqref{eqn.appendix.gravitywaves_pertubation} by
replacing $q_w$ with $q_{vs} = \rho_{vs} / \rho_d$.

\subsection{Atmosphere at rest with a steep mountain}
\label{appendix.mountain}

With the temperature profile $\Tbar(z)$ given in
\eqref{eqn.hydrostatic_mountain:Tprofile} water vapour is at saturation
and there is no cloud water present. The hydrostatic state can then be
computed by solving the following (one-dimensional) non-linear problem:
Find $\pbar,\rhobar_d$, such that
\begin{align*}
  \partial_z \pbar = -(\rhobar_{d} + \rhobar_{vs})g,\qquad
  \pbar = (\rhobar_{d}R_d + \rhobar_{vs} R_v)\Tbar,
  \qquad\text{with}\qquad
  \rhobar_{vs}(\Tbar) = \frac{e_s(\Tbar)}{R_v \Tbar},
\end{align*}
the pressure boundary condition $\pbar(0) = \pref$, and
where the saturation-vapour pressure is computed from the Clausius–Clapeyron
relation \eqref{eqn.saturation-vapour-pre}.

\subsection{Rising thermal with rain}
\label{appendix:subsec.risingthermalrain}
The relative humidity is related to our variables by
\begin{equation*}\HC = \frac{q_{v}}{q_{vs}}\left(\frac{1 + q_{vs}/\epsilon}{1 + q_{v}/\epsilon} \right).
\end{equation*}
The hydrostatic base state is computed by solving the
(one-dimensional) problem: Find $\pbar,\rhobar_d, \rhobar_v, \Tbar$,
such that
\begin{align*}
  \partial_z \pbar &= -(\rhobar_{d} + \rhobar_{v})g,\\
  \pbar &= (\rhobar_{d} R_d + \rhobar_{v} R_v)\Tbar,\\
  \frac{\qbar_v}{\qbar_{vs}}\left(\frac{1 + \qbar_{vs}/\epsilon}{1 + \qbar_v/\epsilon} \right) &= 0.2,\\
  \Tbar \left(\frac{\pref}{\pbar} \right)^{R_d/ c_{pd}} & = T_\text{surf} \left(\frac{\pref}{8.5\times 10^4} \right)^{R_d/ c_{pd}}\exp(1.3\times 10^{-5}z),
\end{align*}
the pressure boundary condition $\pbar(0) = 8.5\times 10^4$, and
with $\qbar_v = \rhobar_v / \rhobar_d, \qbar_{vs} = \rhobar_{vs} / \rhobar_d$.

The initial relative humidity is given by
\begin{equation*}
  \tilde{\HC} =
  \begin{cases}
    \overline{\HC} & r > r_1,\\
    \overline{\HC} + (1 - \overline{\HC})\cos^2\left(\frac{\pi(r - r_1)}{2(r_1 - r_2)}\right) & r_2 \leq r < r_1,\\
    1 & r < r_2,
  \end{cases}
\end{equation*}
with the radius $r=\sqrt{(x - c_x)^2 + (z - c_z)^2}$, $c_x = L/2$,
$c_z = 800\,\unit{\metre}$, $r_1 = 300\,\unit{\metre}$ and
$r_2=200\,\unit{\metre}$. The initial perturbations are then obtained by
solving the problem
\begin{align*}
  \frac{q_v}{q_{vs}(T)}\left(\frac{1 + q_{vs}(T)/\epsilon}{1 + q_v/\epsilon} \right) = \tilde{\HC},\qquad
  T \left(\frac{\pref}{\pbar} \right)^{R_d/ c_{pd}} = \overline{\theta}_d,\qquad
  \rho_d R_d T + \rho_v R_v T = \pbar.
\end{align*}
for $\rho_v, \rho_d$ and $T$. The initial cloud and rain densities are
zero, and the initial velocity is at rest, i.e., $\rho_v=\rho_r=0,
u = (0, 0)^T$ and the initial energy density can be computed from these
variables.

\subsection{Squall Line}
\label{appendix:subsec.squall_line}

The hydrostatic base state is computed by solving the one-dimensional problem: Given $\overline{\theta}, \overline{q}_v$, find $(\pbar,\rhobar_d, \rhobar_v, \Tbar)$ such that
\begin{align*}
  \partial_z \pbar = -(\rhobar_{d} + \rhobar_{v})g,\quad
  \pbar = (\rhobar_{d} R_d + \rhobar_{v} R_v) \Tbar,\quad
  \rhobar_v = \rhobar_d \overline{q}_v, \quad 
  \overline{\theta} = T \left(\frac{\pref}{p} \right)^{R_d / c_{pd}},
\end{align*}
with pressure boundary condition $p(0) = \pref$. The initial temperature perturbation defined by
\begin{equation*}
  \theta' =
  \begin{cases}
    \theta_c\cos\left(\frac{\pi r}{2}\right) &\text{if } r\leq r_c,\\
    0 &\text{if } r > r_c,
  \end{cases}
\end{equation*}
with
\begin{equation*}
  r = \sqrt{\frac{(x - x_c)^2}{r_x^2} + \frac{(z - z_c)^2}{r_z^2}},\quad x_c = 75\,\unit{\km},\quad r_x = 10\,\unit{\km},\quad z_c = 2\,\unit{\km},\quad r_x = 1.5\,\unit{\km},\quad r_c = 1,\quad \theta_c = 3\,\unit{\K}.
\end{equation*}
We obtain the temperature perturbation by solving the constitutive equation for the potential temperature
\begin{equation*}
  \theta = T \left(\frac{\pref}{p} \right)^{R_d / c_{pd}} = T \left(\frac{\pref}{(\rho_d R_d + \rho_v R_v)T} \right)^{R_d / c_{pd}}.
\end{equation*}

\end{document}